\documentclass[a4paper,makeidx,reqno]{amsart}
\usepackage[utf8]{inputenc}
\usepackage[T1]{fontenc}
\usepackage{fixltx2e}
\usepackage{graphicx}
\usepackage{caption}
\usepackage[labelformat=simple]{subcaption}

\DeclareCaptionLabelFormat{carinasub}{#2.}
\makeatletter
\renewcommand{\p@subfigure}{}
\makeatother
\makeatletter
\RequirePackage{ifthen}

\provideboolean{omarfont}
\setboolean{omarfont}{false}

\provideboolean{usesvjour}
\setboolean{usesvjour}{false}%

\provideboolean{nohyperref}
\setboolean{nohyperref}{false}%

\provideboolean{usecleveref}
\setboolean{usecleveref}{false}%

\provideboolean{landscape}
\setboolean{landscape}{false}%

\makeatletter
      \RequirePackage{ifthen}
      \provideboolean{useutopia}
      \setboolean{useutopia}{false}%
      \provideboolean{usesvjour}
      \setboolean{usesvjour}{false}%
      \provideboolean{usebeamer}%
      \setboolean{usebeamer}{false}%
      \provideboolean{usinghyperref}
      \setboolean{usinghyperref}{false}%
\makeatletter
          \RequirePackage{xparse}
          \RequirePackage{ifthen}
          \RequirePackage{iftex}
          \provideboolean{usemathrsfs}%
          \provideboolean{useutopia}
          \setboolean{useutopia}{false}%
          \provideboolean{usebeamer}%
          \setboolean{usebeamer}{false}%
          \provideboolean{isthesis}%
          \provideboolean{isamsltex}%
          \provideboolean{issiamltex}%
          \ifXeTeX%
            \RequirePackage{xltxtra}%
            \ifthenelse{2=1}{
              \RequirePackage{polyglossia}
              \setdefaultlanguage[variant=british]{english}
              \setotherlanguages{french,vietnamese,russian}%
            }{
              \RequirePackage[american,main=british,french,vietnamese]{babel}
            }
            \DeclareSymbolFont{usualmathcal}{OMS}{cmsy}{m}{n}
            \DeclareSymbolFontAlphabet{\mathcalbf}{usualmathcal}
            \ifthenelse{\boolean{useutopia}}{
              \RequirePackage[utopia,euro]{mathdesign}%
              \ifthenelse{\boolean{usebeamer}}{}{
                \setmainfont{Baskerville}
              }
              \RequirePackage[OMLmathbf,OMLmathsfit]{isomath}
            }{
              \RequirePackage{amssymb}    %
              \RequirePackage{bbold}    %
            }
            \RequirePackage{mathrsfs} %
            \setboolean{usemathrsfs}{true}%
          \else%
          \RequirePackage[utf8]{inputenc}%
          \ifthenelse{\boolean{useutopia}}{%
            \RequirePackage[utopia,euro]{mathdesign}%
            \RequirePackage[OMLmathrm,OMLmathbf,OMLmathsfit]{isomath}
            \RequirePackage{bbold}    %
            \RequirePackage{baskervald}
            \DeclareSymbolFont{usualmathcal}{OMS}{cmsy}{m}{n}
            \DeclareSymbolFontAlphabet{\mathcalbf}{usualmathcal}
            \RequirePackage{stmaryrd} %
            \providecommand{\diracdelta}[1][]{\ensuremath{\deltaup_{#1}}}
            
            \providecommand{\lap}{\ensuremath{\Deltaup}}
            \providecommand{\pic}{\ensuremath{\mathrm\pi}}
            \providecommand{\measure}[1]{\ensuremath{\mathcalbf{\uppercase{#1}}}}
            \providecommand{\olnum}[1]{\ensuremath{\mathsfit{#1}}}
            \providecommand{\numsca}[1]{\olnum{#1}}
            \providecommand{\numvec}[1]{\ensuremath{\mathsfbfit{#1}}}%
          }{%
            \RequirePackage{amssymb}    %
            \RequirePackage{mathrsfs} %
            \RequirePackage{bbold}    %
            \RequirePackage{stmaryrd} %
            \setboolean{usemathrsfs}{true}%
            \providecommand{\mathcalbf}{\mathcal}
            \providecommand{\mathsfit}{\mathsf}
          }%
          \RequirePackage[american,main=british,vietnamese]{babel}
          \fi%
          \RequirePackage{xstring}%
          \RequirePackage{chngcntr}%
          \RequirePackage{mathtools}%
          \RequirePackage{nicefrac}
          \RequirePackage{stmaryrd} %
          \RequirePackage{amsthm}%
          \ifthenelse{\boolean{usebeamer}}{}{
            \RequirePackage[shortlabels]{enumitem}%
          }
          \RequirePackage{xspace}%
          \RequirePackage{verbatim}%
          \RequirePackage{fvextra}
          \RequirePackage{fancyvrb}
          \RequirePackage{listings}%
          \RequirePackage{xcolor}
          \RequirePackage{epsfig}
          \RequirePackage{graphicx}
          \RequirePackage[space]{grffile}%
          \RequirePackage{booktabs}%
          \RequirePackage{tikz}%
          \usetikzlibrary{calc}%
          \usetikzlibrary{fadings}
      \def\olprovideenvironment{\@star@or@long\provide@environment}
      \def\provide@environment#1{%
              \@ifundefined{#1}%
                      {\def\reserved@a{\newenvironment{#1}}}%
                      {\def\reserved@a{\renewenvironment{dummy@environ}}}%
              \reserved@a
      }
      \def\dummy@environ{}
      \colorlet{a}{magenta}
      \colorlet{b}{green!75!blue}
      \colorlet{c}{yellow!87.5!red}
      \colorlet{d}{cyan}
      \colorlet{e}{red}
      \colorlet{f}{blue}
      \colorlet{g}{white}
      \colorlet{i}{black}
      \colorlet{h}{i!50!g}
      \colorlet{j}{a!75!g}

      \ifthenelse{\boolean{usinghyperref}}{%
        \providecommand{\linkedurl}[1]{\url{1}}%
        \providecommand{\linkedemail}[1]{\href{mailto:#1}{#1}}%

      }{%
        \providecommand{\linkedurl}[1]{\texttt{#1}}%
        \providecommand{\linkedemail}[1]{\texttt{#1}}%
      }
      \providecommand{\email}[1]{{\linkedemail{#1}}}
      \providecommand{\Ignore}[1]{}
      \providecommand{\ignore}[1]{}
      \providecommand{\freeze}[1]{}%
      
      \providecommand{\crossout}[1]{{\color{i!20} #1}}
      \providecommand{\highlightcolor}{a}
      \providecommand{\highlight}[1]{{\color{\highlightcolor}#1}}

      \providecommand{\memo}[1]{%
        \ensuremath{%
          \framebox{\tiny\textbf{\kern-2pt\textsf{#1}}\kern-2pt}%
        }%
        \xspace}

      \RequirePackage{alphalph}
      \provideboolean{shownotes}
      \setboolean{shownotes}{true}%
      \newcounter{margnote}[page]
      \providecommand{\mgcolor}{a}%
      \providecommand{\mgcolorset}[1]{\renewcommand{\mgcolor}{\alphalph{#1}}}
      \providecommand{\mgcolorsetbycounter}[1]{%
        \newcount\@olmodn
        \@olmodn \value{#1}\relax
        \newcount\@olmodd
        \@olmodd 6\relax
        \newcount\@olmodq %
        \newcount\@olmodr %
        \newcount\@olmodc %
        \@olmodc\@olmodd\relax
        \multiply\@olmodc by 2\relax
        \advance\@olmodc-1\relax
        \@olmodq 0\relax
        \@olmodr\@olmodn\relax
        \ifnum\@olmodr>\@olmodc
        \loop
        \advance\@olmodr-\@olmodd\relax
        \advance\@olmodq1\relax
        \ifnum\@olmodr>\@olmodc
        \repeat
        \fi
        \setcounter{tmpcounter}{\the\@olmodr}
        \stepcounter{tmpcounter}
        \mgcolorset{\value{tmpcounter}}
      }
      \providecommand{\mgcolormake}{\mgcolorsetbycounter{margnote}}
      \providecommand{\margnotecolor}{\mgcolormake}
      \providecommand{\margnotemark}{{\colorbox{\mgcolor}{\tiny\color{g}\upshape\texttt{\arabic{page}.\arabic{margnote}}}}\,}
      \providecommand{\margnote}[2][]{%
        \ifthenelse{%
          \boolean{shownotes}%
        }{%
          \stepcounter{margnote}%
          \margnotecolor%
          \margnotemark %
          \marginpar{%
            \color{\mgcolor}%
            \texttt{\bfseries{%
              \begin{minipage}{2cm}%
                \raggedright\tiny%
                \margnotemark%
                #2%
                \\
                {\ifx|#1|{}\else{ - #1}\fi}%
              \end{minipage}%
              }%
            }%
          }%
        }{%
        }%
      }%
      \providecommand{\mathnote}[2][]{%
        \ifthenelse{%
          \boolean{shownotes}%
        }{%
          \stepcounter{margnote}%
          \margnotecolor%
          \text{%
            \colorbox{g}{%
              \color{\mgcolor}%
              \texttt{\bfseries{%
                \tiny%
                    \margnotemark\,%
                    #2%
                    \ifx|#1|{}\else{ - #1}\fi%
                }%
              }%
            }%
          }%
        }{%
        }%
      }%
      \providecommand{\textnote}[2][]{%
        \ifthenelse{%
          \boolean{shownotes}%
        }{%
          \stepcounter{margnote}%
          \margnotecolor%
          \ \\
          \text{%
              \begin{minipage}{\textwidth}
              \color{\mgcolor}%
              \texttt{%
                \margnotemark%
                #2%
                \ifx|#1|{}\else{ - #1}\fi%
              }%
              \end{minipage}
          }%
        }{}%
      }%

      \providecommand{\Todo}[2][To do:]{
        \ifthenelse{\boolean{shownotes}}{
          \begin{tikzpicture}
           \node[fill=a!17]{
             \begin{minipage}{\textwidth}
               \tiny
               \texttt{#1}
               \texttt{\bfseries{#2}}
             \end{minipage}
           };
          \end{tikzpicture}
        }{}}
      \provideboolean{showrevisions}
      \setboolean{showrevisions}{true}
      \provideboolean{emphrevisions}
      \setboolean{emphrevisions}{false}
      \newcommand{\revisionsheader}{\ \clearpage\Warning{the following part is under development/revision}}
      \newcommand{\revisionsfooter}{\ \newline\Warning{end of part under development/revision}\clearpage}

      \providecommand{\HighlightBox}[2][a!6.25]{
        \begin{center}
          \begin{tikzpicture}
            \node[fill=#1]{
              \begin{minipage}{\textwidth}
                #2
              \end{minipage}
            };
          \end{tikzpicture}
        \end{center}
      }
      \providecommand{\Warning}[1]{    
        \HighlightBox[b!25]{%
          \texttt{\bfseries{\small Warning: #1}}
        }
      }

      \provideboolean{showcomments}
      \providecommand{\margincomment}[1]{
      \ifthenelse{\boolean{showcomments}}{\marginpar{\tiny #1}}{}
      }
      \provideboolean{showchanges}
      \setboolean{showchanges}{false}
      \providecommand{\changes}[2][]{%
        \ifthenelse{\boolean{showchanges}}{{%
            \ifx|#1|{}\else\margnote{#1}\fi%
            \highlight{#2}%
        }}{%
          #2}}
      \providecommand{\mathchanges}[2][]{%
        \ifthenelse{\boolean{showchanges}}{{\ifx|#1|{}\else\mathnote{#1}\fi\highlight{#2}}}{#2}}

      \providecommand{\changefromto}[3][replace with]{%
        \ifthenelse{\boolean{showchanges}}{{%
            \crossout{#2}\margnote{#1}%
          }{%
            \highlight{#3}
          }%
        }{%
          #3\xspace%
        }%
      }
      \providecommand{\ChangePar}[3][]{%
        \ifthenelse{\boolean{showchanges}}{
          {\par\textcolor{i!20}{#2}\ifx|#1|\else\margnote{#1}\fi}{\par\textcolor{a}{#3}}
        }{%
          \par #3%
        }%
      }
      \providecommand{\InsertPar}[1]{
        \ifthenelse{\boolean{showchanges}}
        {{\par$\mapsto$ \textcolor{blue}{#1}}}
        {\par #1}
      }
      
      \providecommand{\mathchangefromto}[3][]{\crossout{#2}\ifx|#1|\else\mathnote{#1}\fi\highlight{#3}}

      \let\trueMakeUppercase\MakeUppercase
      \newcommand{\UCmath}[1]{%
        \begingroup
        \ucmathlist\trueMakeUppercase{#1}%
        \endgroup
      }
      \ifthenelse{\(\boolean{useutopia}\)\OR\(\boolean{usesvjour}\)}{
        \newcommand{\ucmathlist}{%
          \def\alpha{A}%
          \def\beta{B}%
          \let\gamma\Gamma
          \let\delta\Delta
          \def\epsilon{E}%
          \def\varepsilon{E}%
          \def\zeta{Z}%
          \def\eta{H}%
          \let\theta\Theta
          \let\vartheta\Theta
          \def\iota{I}%
          \def\kappa{K}%
          \let\lambda\Lambda
          \def\mu{M}%
          \def\nu{N}%
          \let\xi\Xi
          \def\omicron{O}
          \let\pi\Pi
          \let\varpi\Pi
          \def\rho{P}%
          \def\varrho{P}%
          \let\sigma\Sigma
          \def\varsigma{C}
          \def\tau{T}%
          \let\upsilon\Upsilon
          \let\phi\Phi
          \let\varphi\Phi
          \def\chi{X}%
          \let\psi\Psi
          \let\omega\Omega
      }}{
        \newcommand{\ucmathlist}{
          \def\alpha{\mathrm{A}}%
          \def\beta{\mathrm{B}}%
          \let\gamma\Gamma
          \let\delta\Delta
          \def\epsilon{\mathrm{E}}%
          \def\varepsilon{\mathrm{E}}%
          \def\zeta{\mathrm{Z}}%
          \def\eta{\mathrm{H}}%
          \let\theta\Theta
          \let\vartheta\Theta
          \def\iota{\mathrm{I}}%
          \def\kappa{\mathrm{K}}%
          \let\lambda\Lambda
          \def\mu{\mathrm{M}}%
          \def\nu{\mathrm{N}}%
          \let\xi\Xi
          \let\pi\Pi
          \let\varpi\Pi
          \def\rho{\mathrm{P}}%
          \def\varrho{\mathrm{P}}%
          \let\sigma\Sigma
          \def\tau{\mathrm{T}}%
          \let\upsilon\Upsilon
          \let\phi\Phi
          \let\varphi\Phi
          \def\chi{\mathrm{X}}%
          \let\psi\Psi
          \let\omega\Omega
        }
      }
      \providecommand{\mathscript}
      	       {\mathscr}

      \providecommand{\cD}{\ensuremath{\mathscript D}\xspace}
      \providecommand{\cE}{\ensuremath{\mathscript E}\xspace}

      \providecommand{\bbbold}{\mathbb}

      \providecommand{\rK}{\ensuremath{\bbbold K}\xspace}

      \providecommand{\rN}{\ensuremath{\bbbold N}\xspace}
      
      \providecommand{\rP}{\ensuremath{\bbbold P}\xspace}
      
      \providecommand{\rR}{\ensuremath{\bbbold R}\xspace}
      
      \providecommand{\rT}{\ensuremath{\bbbold T}\xspace}

      \providecommand{\Ae}[1][]{\ensuremath{\ifx|#1|{\ }\else{\:#1\text{-}}\fi\text{almost everywhere\xspace}}}
      \providecommand{\Aeon}[1][]{\ensuremath{\text{\Ae[#1]\ on }}}
      \providecommand{\Aein}[1][]{\ensuremath{\text{\Ae[#1]\ in }}}
      \providecommand{\Aa}[1][]{\ensuremath{\text{ for }\ifx|#1|{}\else{\:#1\text{-}}\fi\text{almost all }}}
      \providecommand{\as}[1][]{\ensuremath{\ifx|#1|{\ }\else{#1\text{-}}\fi\text{almost surely}}\xspace}
      \providecommand{\aposteriori}{aposteriori\xspace}
      \providecommand{\Aposteriori}{Aposteriori\xspace}
      
      \providecommand{\apriori}{apriori\xspace}

       \providecommand{\naturals}{\ensuremath{\rN}}
       
       \providecommand{\NO}[1][]{\ensuremath{\naturals_0\ifx|#1|{}\else^{#1}\fi}}

       \providecommand{\reals}{\rR}

       \providecommand{\R}[1]{{\reals^{#1}}}

       \providecommand{\fieldmats}[3][F]{\csname#1\endcsname{#2\times#3}}
       
       \providecommand{\fieldtens}[3][F]{\csname#1\endcsname{{#2}_1\times\dotsb\times{#2}_{#3}}}

       \providecommand{\RO}[1][]{{\reals_{0+}\ifx|#1|{}\else^{#1}\fi}}
       \providecommand{\RP}[1][]{{\reals_+\ifx|#1|\else^{#1}\fi}}
       \providecommand{\Rneg}[1][]{{\reals_-\ifx|#1|\else^{#1}\fi}}

       \providecommand{\ring}[1][A]{\csname r#1\endcsname}
       \providecommand{\field}[1][K]{\csname r#1\endcsname}

       \providecommand{\torus}[1]{\rT\ifthenelse{\equal{#1}1}{}{^#1}}
      
       \providecommand{\one}{\ensuremath{\bbbold 1}}
       \providecommand{\zerofun}{\ensuremath{\bbbold 0}}
       \providecommand{\ones}[1][]{\one\ifx|#1|\else_{#1}\fi}
       \providecommand{\zeros}[1][]{\zerofun\ifx|#1|\else_{#1}\fi}
       \providecommand{\charfun}[1]{\one_{#1}}
       \providecommand{\iverson}[1]{\one_{\qb{#1}}}

       \providecommand{\diracdelta}[1][]{\ensuremath{{\mathrm{\delta}}\ifx|#1|{}\else_{#1}\fi}}

       \providecommand{\pic}{\pi}%
       \providecommand{\pifracl}[2][]{\fracl{\ifx|#1|\else#1\fi\pic}{#2}}
       \providecommand{\pifrac}[2][]{\frac{\ifx|#1|\else#1\fi\pic}{#2}}

       \providecommand{\take}{\smallsetminus}
       \providecommand{\takesetof}[1]{\take\setof{#1}}
       \providecommand{\takeset}\takesetof
       \providecommand{\takeel}\oldneg%

       \providecommand{\closure}[2][]{\ifx|#1|\overline{#2}\else\operatorname{clos}_{#1}{#2}\fi}
       \providecommand{\inner}{\cdot}
       \providecommand{\vecprod}{\times}
       \providecommand{\outerp}{\wedge}

       \providecommand{\W}{\ensuremath{\varOmega}\xspace}
       \providecommand{\qgroup}[2][1]{{#2}}%

       \providecommand{\qp}[2][]{\ensuremath{\ifx|#1|\left(\else\csname#1\endcsname(\fi{#2}\ifx|#1|\right)\else\csname#1\endcsname)\fi}}

       \providecommand{\qpreg}[1]{\ensuremath{(#1)}}
       \providecommand{\qpbig}[1]{\qp[big]{#1}}%
       \providecommand{\qpBig}[1]{\ensuremath{\Big(#1\Big)}}
       \providecommand{\qpbigg}[1]{\ensuremath{\bigg(\!#1\!\bigg)}}
       \providecommand{\qpBigg}[1]{\ensuremath{\Bigg(\!#1\!\Bigg)}}
       \providecommand{\qb}[2][]{\ifx|#1|\left[\else\csname#1\endcsname[\fi{#2}\ifx|#1|\right]\else\csname#1\endcsname]\fi}
       \providecommand{\qc}[2][]{\ensuremath{\ifx|#1|\left\{\else\csname#1\endcsname\{%
           \fi{#2}\ifx|#1|\right\}\else\csname#1\endcsname\}\fi}}

       \providecommand{\qa}[2][]{\ifx|#1|\left\langle\else\csname#1\endcsname\langle%
         \fi{#2}\ifx|#1|\right\rangle\else\csname#1\endcsname\rangle\fi}%
       \providecommand{\qareg}[1]{\ensuremath{\langle#1\rangle}}
       \providecommand{\qabig}[1]{\ensuremath{\big\langle#1\big\rangle}}
       \providecommand{\qaBig}[1]{\ensuremath{\Big\langle#1\Big\rangle}}
       \providecommand{\qabigg}[1]{\ensuremath{\bigg\langle#1\bigg\rangle}}
       \providecommand{\qaBigg}[1]{\ensuremath{\Bigg\langle#1\Bigg\rangle}}
       \providecommand{\qv}[2][]{\ifx|#1|\left|\else\csname#1\endcsname|%
         \fi{#2}\ifx|#1|\right|\else\csname#1\endcsname|\fi}%

       \providecommand{\opinter}[2]{\ensuremath{\left(#1,#2\right)}\xspace}

       \providecommand{\clinter}[2]{\ensuremath{\left[#1,#2\right]}\xspace}
       
       \providecommand{\opclinter}[2]{\ensuremath{\left(#1,#2\right]}\xspace}
       \providecommand{\clopinter}[2]{\ensuremath{\left[#1,#2\right)}\xspace}   
       \providecommand{\opintertopinfty}[1]{\opinter{#1}\infty}
       \providecommand\optyinter\opintertopinfty
       \providecommand{\opinterbotinfty}[1]{\opinter{-\infty}{#1}}
       \providecommand\tyopinter\opinterbotinfty
       \providecommand{\clintertopinfty}[1]{\clopinter{#1}\infty}
       \providecommand{\cltyinter}\clintertopinfty
       \providecommand{\clinty}\clintertopinfty
       
       \providecommand{\clinterbotinfty}[1]{\opclinter{-\infty}{#1}}
       \providecommand{\tyclinter}\clinterbotinfty
       \providecommand{\expp}[1]{\ensuremath{\e^{#1}}}
       
       \providecommand{\compowqp}[2]{\ensuremath{\qp{\!#2\!\!}^{\kern -.4em #1}\!}}
       
       \providecommand{\powqpreg}[2]{\ensuremath{%
           \qpreg{#2}^{\kern 0em\lower .1ex\hbox{\scriptsize $#1$}}\kern-.3em}}
       \providecommand{\powqpbig}[2]{\ensuremath{%
           \qpbig{#2}^{\kern -.2em\lower .3ex\hbox{\scriptsize $#1$}}\kern-.3em}}
       \providecommand{\powqpBig}[2]{\ensuremath{%
           \qpBig{#2}^{\kern -.2em\lower .3ex\hbox{\scriptsize $#1$}}\kern-.3em}}
       \providecommand{\powqpbigg}[2]{\ensuremath{%
           \qpbigg{#2}^{\kern -.2em\lower .3ex\hbox{\scriptsize $#1$}}\kern-.3em}}
       \providecommand{\powqpBigg}[2]{\ensuremath{%
           \qpBigg{#2}^{\kern -.2em\lower .3ex\hbox{\scriptsize $#1$}}}}
       \providecommand{\powp}[3][]{{#3}\ifx|#1|^{#2}\else{#1}^{#2}\fi}%
       \providecommand{\pow}[2][]{\ifx|#1|\operatorname{pow}^{#2}\else\powp{#2}{#1}\fi}%
       \providecommand{\ppow}[3][]{\powp[#1]{#3}{#2}}

       \providecommand{\powqp}[3][]{\powp[#1]{#2}{\qp{#3}}}%
       \providecommand{\qppow}[3][]{\ppow[#1]{\qp{#2}}{#3}}%

       \providecommand{\powsqrt}[2][2]{\powp{\fracl1{#1}}{#2}}
       \providecommand{\powqpsqrt}[2][2]{\powsqrt[#1]{\qp{#2}}}

       \providecommand{\norm}[2][]{\ifx|#1|\left|\else\csname#1\endcsname|\fi#2\ifx|#1|\right|\else\csname#1\endcsname|\fi}
       \providecommand{\normon}[3][]{\norm[#1]{#2}_{#3}}

       \providecommand{\abs}[2][]{\ensuremath{\ifx|#1|{\left|#2\right|}\else{\csname#1\endcsname|{#2}\csname#1\endcsname|}\fi}}

       \providecommand{\Norm}[2][]{\ifx|#1|\left\|\else\csname#1\endcsname\|\fi{#2}\ifx|#1|\right\|\else\csname#1\endcsname\|\fi}

       \providecommand{\Normon}[3][]{\Norm[#1]{#2}_{#3}}
       \providecommand{\Normonspace}[3][]{\Norm[#1]{#2}_{\vecspace{#3}}}

       \providecommand{\normonsob}[5][]{\normon[#1]{#2}{\sob{#3}{#4}\if|#5|{}\else(#5)\fi}}
       \providecommand{\Normonsob}[5][]{\Normon[#1]{#2}{\sob{#3}{#4}\if|#5|{}\else(#5)\fi}}
       \providecommand{\normonsobh}[4][]{\normon[#1]{#2}{\sobh{#3}\if|#4|{}\else(#4)\fi}}
       \providecommand{\normonsobhz}[4][]{\normon[#1]{#2}{\sobhz{#3}\if|#4|{}\else(#4)\fi}}
       \providecommand{\Normonsobh}[4][]{\Normon[#1]{#2}{\sobh{#3}\if|#4|{}\else(#4)\fi}}
       \providecommand{\Normonsobhz}[4][]{\Normon[#1]{#2}{\sobhz{#3}\if|#4|{}\else(#4)\fi}}
       \providecommand{\Normonleb}[4][]{\Normon[#1]{#2}{\leb{#3}\if|#4|\else(#4)\fi}}
       \providecommand{\Normsupon}[3][]{\Normonleb[#1]{#2}\infty{#3}}
       \providecommand{\ltwop}[3][]{\ensuremath{\qa{#2,#3}\ifx|#1|\else_{#1}\fi}}
       \providecommand{\ltwopreg}[2]{\ensuremath{\qareg{#1,#2}\ifx|#1|\else_{#1}\fi}}
       \providecommand{\ltwopbig}[2]{\ensuremath{\qabig{#1,#2}\ifx|#1|\else_{#1}\fi}}
       \providecommand{\ltwopBig}[2]{\ensuremath{\qaBig{#1,#2}\ifx|#1|\else_{#1}\fi}}
       \providecommand{\ltwopbigg}[2]{\ensuremath{\qabigg{#1,#2}\ifx|#1|\else_{#1}\fi}}
       \providecommand{\ltwopBigg}[2]{\ensuremath{\qaBigg{#1,#2}\ifx|#1|\else_{#1}\fi}}
       \providecommand{\ltwopon}[3]{\ltwop{#1}{#2}_{#3}}

       \providecommand{\duality}[3][]{%
         \ifx#1\left\langle#2\middle|#3\right\rangle%
         \else%
         #1\langle%
         #2%
         #1|%
         #3%
         #1\rangle%
         \fi}%
       
       \providecommand{\average}[2][]{{\qa{#2}\ifx|#1|\else_{#1}\fi}}

       \providecommand{\ensemble}[2]{\ensuremath{\left\{ #1:\;#2 \right\}}}
       \providecommand{\setofsuch}{\ensemble}%

       \providecommand{\floor}[1]{\ensuremath{\left\lfloor{#1}\right\rfloor}}
       \providecommand{\ceil}[1]{\ensuremath{\left\lceil{#1}\right\rceil}}

       \providecommand{\setof}[1]{{\qc{#1}}}
       \providecommand{\pair}[2]{\qp{#1,#2}}

       \providecommand{\setpair}[2]{\setof{#1,#2}}

       \providecommand{\conditionalto}[1]{{\left|{#1}\right.}}

      \providecommand{\measure}[1]{\ensuremath{\mathcalbf{\MakeUppercase{#1}}}}
      
      \providecommand{\probmeasure}[2][]{{\measure{#2}}\ifx|#1|\else_{#1}\fi}
      \providecommand{\Prob}{}
      \renewcommand{\Prob}[1][]{\probmeasure[{#1}]{p}}

      \providecommand{\randvars}[1][\Prob]{\operatorname{RV}\ifx|#1|{}\else{(#1)}\fi}
      \providecommand{\discrandvars}[1][\Prob]{\operatorname{DRV}\ifx|#1|{}\else{({#1)}\fi}} 
      \providecommand{\contrandvars}[1][\Prob]{\ensuremath{\operatorname{CDRV}\ifx|#1|{}\else(#1)\fi}} 
       \def\env@matrix{\hskip -\arraycolsep
        \let\@ifnextchar\new@ifnextchar
        \array{*\c@MaxMatrixCols c}}
       \renewcommand*\env@matrix[1][c]{\hskip -\arraycolsep
         \let\@ifnextchar\new@ifnextchar
         \array{*\c@MaxMatrixCols #1}}
       \providecommand{\irow}[2]{#1_{#2}}%
       \providecommand{\icol}[2]{#1^{#2}}%

       \providecommand{\ijrowcol}[3]{\icol{\irow{#1}{#2}}{#3}}
       \providecommand{\entry}[1]{\qb{#1}}
       \providecommand{\vecentry}[2]{\irow{#1}{#2}}

       \providecommand{\colvecentry}\vecentry
       \providecommand{\covecentry}[2]{\icol{#1}{#2}}
       \providecommand\rowvecentry\covecentry

       \providecommand{\rowof}[1]{\qb{#1}}
       \providecommand{\disvecof}[2][r]{\begin{bmatrix}[#1]#2\end{bmatrix}}
       \providecommand{\vecof}[2][r]{\mathchoice{\disvecof[#1]{#2}}{\qp{#2}}{\qp{#2}}{\qp{#2}}}

       \providecommand{\getentryi}[2]{\irow{\entry{#1}}{#2}}
       \providecommand{\getcolentry}\getentryi
       
       \providecommand{\getvecentry}[2]{\getentryi{\vec #1}{#2}}

       \providecommand{\discolvecitwo}[1]{\discolvectwo{\vecentry{#1}1}{\vecentry{#1}2}}
       \providecommand{\discolvecitwoz}[1]{\discolvectwo{\vecentry{#1}0}{\vecentry{#1}1}}
       \providecommand{\discolvecintwo}\discolvecitwo%
       \providecommand{\discolvecitwoz}[1]{\discolvectwo{\vecentry{#1}0}{\vecentry{#1}1}}

       \providecommand{\vecdots}[2]{\qp{{#1},\dotsc,{#2}}}
       \providecommand{\vecdotsfromto}[3]{\vecdots{\irow{#1}{#2}}{\irow{#1}{#3}}}%

       \providecommand{\dismatof}[2][r]{%
         \mathchoice{%
           \begin{bmatrix}[#1]#2\end{bmatrix}%
         }{%
           \qb{
             \begin{smallmatrix}
               #2
             \end{smallmatrix}
           }
         }{%
           \qb{
             \begin{smallmatrix}
               #2
             \end{smallmatrix}
           }
         }{%
           \qb{
             \begin{smallmatrix}
               #2
             \end{smallmatrix}
           }
         }
       }

       \providecommand{\matentry}[3]{\ijrowcol{#1}{#2}{#3}}

       \providecommand{\block}[5]{\ijrowcol{#1}{\ifx#2#3{\rowof{#2}}\else\rowof{{#2}\dotsc{#3}}\fi}{\ifx#4#5{\rowof{#4}}\else\rowof{{#4}\dotsc{#5}}\fi}}
       \providecommand{\colblock}[3]{\getvecentry{#1}{\ifx#2#3{#2}\else\fromto{#2}{#3}\fi}}

       \providecommand{\dismatskeldots}[4]{
         \dismatof[c]{
           #1&\dotsc&#3
           \\
           \vdots & \ddots &\vdots
           \\
           #2&\dotsc&#4
         }
       }
       \providecommand{\dismatcommfromtofromto}[5]{
         \dismatskeldots{#1#2#4}{#1#3#4}{#1#2#5}{#1#3#5}
       }
       \providecommand{\dismatcustfromtofromto}[6][matentry]{
         \dismatcommfromtofromto{\csname#1\endcsname{#2}}#3#4#5#6
       }
       \providecommand{\dismatcustfromtofromto}[6][matentry]{
         \dismatskeldots{%
           \csname#1\endcsname{#2}{#3}{#4}%
         }{%
           \csname#1\endcsname{#2}{#3}{#6}%
         }{%
           \csname#1\endcsname{#2}{#5}{#4}%
         }{%
           \csname#1\endcsname{#2}{#5}{#6}%
         }%
       }%
       \providecommand{\dismatcustfromtofromto}[6][matentry]{
         \dismatof{
           \csname#1\endcsname{#2}{#3}{#4}&\dotsc&\csname#1\endcsname{#2}{#3}{#6}
           \\
           \vdots & \ddots &\vdots
           \\
           \csname#1\endcsname{#2}{#5}{#4}&\dotsc&\csname#1\endcsname{#2}{#5}{#6}
         }
       }

       \providecommand{\dissysaxbdotsnm}[5]{\begin{matrix}[r]%
           \matentry{#1}11\vecentry{#2}1&+\dotsb&+\matentry{#1}1{#5}\vecentry{#2}{#5}
           &
           =
           \ifx|#3|0\else{\vecentry {#3}1}\fi
           \\
           \dotsb
           \\
           \matentry{#1}{#4}1\vecentry{#2}1&+\dotsb&+\matentry{#1}{#4}{#5}\vecentry{#2}{#5}
           &
           =
           \ifx|#3|0\else{\vecentry {#3}{#4}}\fi
       \end{matrix}}
       \providecommand{\seqof}[1]{\qp{#1}}%
       \providecommand{\seqs}[2]{\seqof{#1}_{#2}}
       \providecommand{\sets}[2]{\setof{#1}_{#2}}%
       \providecommand{\seqi}[3][]{\seqs{#2_{#3}}{\ifx|#1|{#3}\else{{#3}\in{#1}}\fi}}%
       \providecommand{\sequ}[3][]{\seqs{#2^{#3}}{\ifx|#1|{#3}\else{{#3}\in{#1}}\fi}}%
       \providecommand{\subseqi}[4][]{\seqs{#2_{{#3}_{#4}}}{\ifx|#1|{#4}\else{{#4}\in{#1}}\fi}}%

       \providecommand{\seti}[3][]{\sets{#2_{#3}}{\ifx|#1|_{#3}\else_{{#3}\in{#1}}\fi}}%
       \providecommand{\setu}[3][]{\sets{#2^{#3}}{\ifx|#1|{#3}\else{{#3}\in{#1}}\fi}}%
       \let\liminf\relax
       \DeclareMathOperator*{\liminf}{liminf}
       \let\limsup\relax
       \DeclareMathOperator*{\limsup}{limsup}
       \providecommand{\limofat}[3][]{\ensuremath{\lim_{\ifx|#1|{}\else{#1\ni}\fi#3}{#2}}}
       \providecommand{\limsupofat}[3][]{\ensuremath{\limsup_{\ifx|#1|{}\else{#1\ni}\fi#3}{#2}}}
       \providecommand{\liminfofat}[3][]{\ensuremath{\liminf_{\ifx|#1|{}\else{#1\ni}\fi#3}{#2}}}

       \providecommand{\seqsfromto}[4]{\seqs{#1}{\rangefromto{#2}{#3}{#4}}}
       \providecommand{\setsfromto}[4]{\setofsuch{#1}{\rangefromto{#2}{#3}{#4}}}
       \providecommand{\setsifromto}[4]{\ensemble{{#1}_{#2}}{\rangefromto{#2}{#3}{#4}}}

       \providecommand{\stringdotsfrom}[3][]{\ensuremath{#2\ifx|#1|\else#1\fi\,#3\ifx|#1|\else#1\fi\,\dotsc}}
       \providecommand{\listdotsfrom}[3][]{\ensuremath{#2\ifx|#1|\else#1\fi,#3\ifx|#1|\else#1\fi,\dotsc}}
       \providecommand{\stringdotsfromto}[3][]{\ensuremath{#2\ifx|#1|\else#1\fi\,\dotsc\,#3\ifx|#1|\else#1\fi}}
       \providecommand{\listdotsfromto}[3][]{\ensuremath{#2\ifx|#1|\else#1\fi,\dotsc,#3\ifx|#1|\else#1\fi}}
       \providecommand{\listifromto}[5][]{\ensuremath{{#2}_{#3}\ifx|#1|\else#1\fi},\text{ for }\ensuremath{\rangefromto{#3}{#4}{#5}}\xspace}
       \providecommand{\listufromto}[5][]{\ensuremath{{#2}^{#3}\ifx|#1|\else#1\fi},\text{ for }\ensuremath{\rangefromto{#3}{#4}{#5}}\xspace}
       \providecommand{\listitwo}[2][]{\ensuremath{#2_1\ifx|#1|\else#1\fi,#2_2\ifx|#1|\else#1\fi}}
       \providecommand{\listutwo}[2][]{\ensuremath{#2^1\ifx|#1|\else#1\fi,#2^2\ifx|#1|\else#1\fi}}
       \providecommand{\listithree}[2][]{\ensuremath{#2_1\ifx|#1|\else#1\fi,#2_2\ifx|#1|\else#1\fi,#2_3\ifx|#1|\else#1\fi}}
       \providecommand{\listithreez}[2][]{\ensuremath{#2_0\ifx|#1|\else#1\fi,#2_1\ifx|#1|\else#1\fi,#2_2\ifx|#1|\else#1\fi}}
       \providecommand{\listifourz}[2][]{\ensuremath{#2_0\ifx|#1|\else#1\fi,#2_1\ifx|#1|\else#1\fi,#2_2\ifx|#1|\else#1\fi,#2_3\ifx|#1|\else#1\fi}}
       \providecommand{\listuthree}[2][]{\ensuremath{#2^1\ifx|#1|\else#1\fi,#2^2\ifx|#1|\else#1\fi,#2^3\ifx|#1|\else#1\fi}}

       \providecommand{\listudotsfromto}[4][]{\listdotsfromto[#1]{#2^{#3}}{#2^{#4}}}

       \providecommand{\ltidotsfromto}[3]{\ensuremath{#1_{#2}\lt\dotsb\lt#1_{#3}}}

       \providecommand{\sums}[2]{\ensuremath{\sum_{#1\in#2}}}

       \providecommand{\sumifromto}[3]{\ensuremath{\sum_{#1=#2}^{#3}}}

       \providecommand{\jump}[2][]{\ensuremath{\left\llbracket #2\right\rrbracket\ifx|#1|{}\else_{#1}\fi}}

       \providecommand{\fromto}[2]{\ensuremath{\setof{#1\dotsc#2}}}%

       \providecommand{\integerbetween}[2]{\ensuremath{={#1},\dotsc,{#2}}}
       \providecommand{\betweenonetwoend}[3]{\ensuremath{={#1},{#2},\dotsc,{#3}}}

       \providecommand{\rangefromto}[3]{\ensuremath{#1\integerbetween{#2}{#3}}}

       \providecommand{\e}{\ensuremath{\operatorname{e}\!}\xspace}

       \providecommand{\backdiff}{\overleftarrow\partial}

       \providecommand{\forediff}{\overrightarrow\partial}

       \providecommand{\d}{}
       \renewcommand{\d}[1][]{\ensuremath{\operatorname{d}\!\ifx|#1|\else{_{#1}}\fi}}
       
       \providecommand{\ds}[1][]{\d{\measure S}}
       \providecommand{\D}[1][]{\ensuremath{\operatorname{D}\!\ifx|#1|\else{_{#1}}\fi}}

      \providecommand{\registered}%
      {\ensuremath{^\text{\textregistered}}}

      \providecommand{\tand}{\ensuremath{\text{ and }}}
      \providecommand{\quand}{\ensuremath{\quad\tand\quad}}
      \providecommand{\tif}{\ensuremath{\text{ if }}}

      \providecommand{\tor}{\ensuremath{\text{ or }}}

      \providecommand{\constant}[1]{\ensuremath{C_{#1}}}
      \providecommand{\constext}[2][]{\constant{\textup{#2}{\ifx|#1|{}\else{,\ensuremath{#1}}\fi}}}            %
      \providecommand{\constref}[2][]{\ensuremath{\constant{\textup{\ref{#2}{\ifx|#1|{}\else{,\ensuremath{#1}}\fi}}}}}
      \providecommand{\constdef}[2][]{\label{#2}\ensuremath{\constant{\textup{\ref{#2}{\ifx|#1|{}\else{,\ensuremath{#1}}\fi}}}}}

      \providecommand{\funkref}[3][]{\ensuremath{{#3}_{\textup{\ref{#2}{\ifx|#1|{}\else{,\ensuremath{#1}}\fi}}}}}

      \providecommand{\diam}{\operatorname{diam}}
      
      \providecommand{\curl}{\operatorname{curl}}
      \renewcommand{\curl}[1][]{\nabla\ifx|#1|{}\else\kern-2pt_{#1}\fi\kern-2pt\vecprod}
      \renewcommand{\div}[1][]{\nabla\ifx|#1|{}\else\kern-2pt_{#1}\fi\kern-1pt\inner}
      \providecommand{\divof}[2][]{\div[#1]\ifx|#2|{}\else\qb{#2}\fi}%
      \providecommand{\divideabyb}[2]{\operatorname{div}(a,b)}
      \providecommand{\grad}{}
      \renewcommand{\grad}[1][]{\nabla\ifx|#1|\else_{#1}\fi}
      \providecommand{\gradof}[2][]{\grad[#1]\qb{#2}}
      \providecommand{\rot}[1][]{\nabla\ifx|#1|\else_{#1}\fi\outerp}
      
      \providecommand{\rowdiv}[1][]{\D\ifx|#1|{}\else\kern-1pt_{#1}\kern-2pt\fi\cdot}
      \providecommand{\rowdivof}[2][]{\rowdiv[#1]\ifx|#2|{}\else\qb{#2}\fi}

      \providecommand{\area}{\measure s}
      \providecommand{\areaof}[1]{\area\qp{#1}}

      \providecommand{\interior}{\operatorname{int}}
      \providecommand{\inv}[1][]{\operatorname{inv}\ifx|#1|\else^{#1}\fi}
      \providecommand{\ivt}[1]{\operatorname{ivt}\ifx|#1|\else^{#1}\fi}
      \providecommand\tensorinvariant\ivt
      
      \providecommand{\inverse}[2][]{\powp[#1]{-1}{#2}}

      \providecommand{\mod}{}
      \renewcommand{\mod}[1][]{\operatorname{mod}\ifx|#1|\else\kern-1pt_{#1}\fi}
      
      \let\oldfrac\frac
      \renewcommand{\frac}[3][]{\ifx|#1|\oldfrac{#2}{#3}\else\begin{array}{#1}{#2}\\\hline{#3}\end{array}\fi}
      \providecommand{\fracl}[3][]{\ifx|#1|\nicefrac{#2}{#3}\else{#2}#1/{#3}\fi}

      \providecommand{\half}[1]{\frac{#1}2}

      \providecommand{\qpfracl}[3][]{\qp{\ifx|#1|\fracl{#2}{#3}\else{#2}#1/{#3}\fi}}
      \providecommand{\qpfrac}[3][]{\qp{\ifx|#1|\frac{#2}{#3}\else{#2}#1/{#3}\fi}}

      \providecommand{\absfracl}[3][]{\abs{\ifx|#1|\fracl{#2}{#3}\else{#2}#1/{#3}\fi}}
      \providecommand{\absfrac}[3][]{\abs{\ifx|#1|\frac{#2}{#3}\else{#2}#1/{#3}\fi}}
      \providecommand{\fraclff}[3][]{\ifx|#1|{#2}/{#3}\else{#1}\fracl{#2}{#3}\fi}

      \providecommand{\eye}[1][]{\vec{\mathrm I}\ifx|#1|{}\else_{#1}\fi}%
      \providecommand{\numeye}[1][]{\boldsymbol{\mathsf{I}}\ifx|#1|{}\else_{#1}\fi}%
      \providecommand{\doteye}{{\ooalign{$\eye${\kern-0.41em\raise1.56ex\hbox{\tiny$\bullet$}}}}}%
      \providecommand{\Eye}[1]{
        \begin{bmatrix}
        \ifthenelse{#1>1}{
          \ifthenelse{#1>2}{
            \ifthenelse{#1>3}{
              \ifthenelse{#1>4}{
                1&\zeroentry&\dotso&\zeroentry
                \\
                \zeroentry&1&\dotso&\zeroentry
                \\
                \vdots&\vdots&\ddots&\vdots
                \\
                \zeroentry&\zeroentry&\dotso&1
              }{        
                1&\zeroentry&\zeroentry&\zeroentry
                \\
                \zeroentry&1&\zeroentry&\zeroentry
                \\
                \zeroentry&\zeroentry&1&\zeroentry
                \\
                \zeroentry&\zeroentry&\zeroentry&1
              }
            }{
              1&\zeroentry&\zeroentry
              \\
              \zeroentry&1&\zeroentry
              \\
              \zeroentry&\zeroentry&1
            }
          }{
            1&\zeroentry
            \\
            \zeroentry&1
          }
        }{
          1
        }
        \end{bmatrix}
      }
      \providecommand{\Id}{\operatorname{Id}}                   %

      \providecommand{\lebmeas}[1][]{\measure L^{#1}}     %
      \providecommand{\lebmeasof}[2][]{\ifx|#1|\left|#2\right|\else\lebmeas[#1]\qp{#2}\fi}         %
      \providecommand{\area}[1]{\operatorname{area}#1}          %

      \providecommand{\meshsize}[1][]{h\ifx|#1|\else_{#1}\fi}
      \providecommand{\meshsizemesh}[1]{\meshsize[\mesh{#1}]}
      \providecommand{\maxi}[2]{#1\vee#2}                       %

      \let\oldneg\neg
      \renewcommand{\neg}[1]{\left[#1\right]_-}
      \providecommand{\Oh}{\operatorname{O}}                   %
      \providecommand{\dash}[1][']{\ifthenelse{\equal{#1}{'}\OR\equal{#1}{''}}{#1}{^{(#1)}}}
      
      \providecommand{\pdfrac}[2][]{\ensuremath{\frac{\partial\ifx|#1|\phantom{#2}\else{#1}\fi}{\partial{#2}}}} %
      \providecommand{\pdfracpow}[3][]{\ensuremath{\frac{\partial^{#3}\ifx|#1|\phantom{#2}\else{#1}\fi}{\partial{#2}^{#3}}}} %

      \providecommand{\pd}[2][]{\ensuremath{\partial_{#2}}{\ifx|#1|{}\else{\qb{#1}}\fi}} %

      \providecommand{\pdt}[1][]{\pd[#1]t}                       %
      \providecommand{\pdtt}[1][]{\pd[#1]{tt}}                    %
      \providecommand{\dd}[2][]{\ensuremath{\ifx|#1|\frac{\d}{\d{#2}}\else\frac[l]{\d{#1}}{\d{#2}}\fi}}    %

      \providecommand{\ddt}[1][]{\dd[{#1}]t}    %

      \renewcommand{\Im}{\operatorname{im}}                 %
      \renewcommand{\Re}{\operatorname{re}}                 %
      \providecommand{\imaginpart}[1][]{\Im{\ifx|#1|{}\else\qp{#1}\fi}} %
      \providecommand{\realpart}[1][]{\Re{\ifx|#1|{}\else\qp{#1}\fi}} %
      \providecommand\determinant\det

      \providecommand{\transpose}{\intercal}%

      \providecommand{\transposed}{{}^\transpose}

      \providecommand{\orthogonalto}[1][]{\ensuremath{\perp\ifx|#1|{}\else{\!_{#1}\,}\fi}}
      
      \providecommand{\rowof}[1]{\ensuremath{\mathchoice{\vecof{#1}}{\qb{#1}}{\qb{#1}}{\qb{#1}}}}

      \providecommand{\rowvectwo}[2]{\ensuremath{\vecof{#1,#2}}}

      \providecommand{\colvectwo}[2]{\ensuremath{%
          \mathchoice%
              {\discolvectwo{#1}{#2}}%
              {\rowvectwo{#1}{#2}\transposed}%
              {\rowvectwo{#1,#2}\transposed}%
              {\rowvectwo{#1,#2}\transposed}%
        }
      }
      \providecommand{\coltwovec}\colvectwo

      \providecommand{\colvecitwoz}[1]{\colvectwo{\vecentry{#1}0}{\vecentry{#1}1}}

      \providecommand{\discolvec}[2][r]{\ensuremath{\begin{bmatrix}[#1]#2\end{bmatrix}}}

      \providecommand{\discolvectwo}[3][r]{\ensuremath{\discolvec[#1]{#2\\#3}}}

      \providecommand{\discolvecitwo}[1]{\discolvectwo{\vecentry{#1}1}{\vecentry{#1}2}}

      \provideboolean{showzeroentries}%
      \setboolean{showzeroentries}{true}%
      \providecommand{\zeroentry}{\ifthenelse{\boolean{showzeroentries}}{{0}}{\phantom0}}
      \providecommand{\zeroentrywarning}{\ifthenelse{\boolean{showzeroentries}}{}{%
          \ensuremath{\text{($0$ entries omitted)\xspace}}}}

      \providecommand{\dismatrix}[2][r]{\ensuremath{\dismatof[#1]{#2}}}
      \providecommand{\dismattwo}[5][r]{\dismatrix[#1]{#2&#3\\ #4&#5}}

      \providecommand{\smint}{\ensuremath{{\text{\textbf{/}}}\kern-.75em\smallint}}
      \renewcommand{\smint}[1][]{\lower12.3pt\hbox{\begin{tikzpicture}\draw[line width=.75pt] (-3pt,-0.5)--(1pt,-0.5) node[pos=0.6]{$\int$};\path (3pt,-24pt)node {\scriptsize $#1$};\end{tikzpicture}}}

      \providecommand{\lap}{\ensuremath{\mathrm\Delta}}
      \providecommand{\lapin}[1][]{\lap\ifx|#1|\else_{#1}\fi}
      \providecommand{\normalsymbol}{\operatorname{\mathbf{n}}}
      \renewcommand{\normalsymbol}{\vec{\operatorname{n}}}
      \providecommand{\normal}[1][]{\normalsymbol\ifx|#1|\else_{#1}\fi}%
      \providecommand{\norm@l}[1][]{\normalsymbol\ifx|#1|\else_{#1}\fi}%
      \providecommand{\normalto}[2][]{\ensuremath{\norm@l[#2]\ifx|#1|\else\qp{#1}\fi}}
      \providecommand{\normalder}[1][]{\ensuremath{\norm@l\ifx|#1|\else\qp{#1}\fi{\inner\grad}}}
      \providecommand{\normalderto}[2][]{\ensuremath{\normalto[#1]{#2}{\inner\grad}}}

      \providecommand{\tangentialsymbol}{\operatorname{\textbf{t}}}
      \providecommand{\tangentialto}[2][]{\tangentialsymbol\ifx|#1|\else^{#1}\fi\ifx|#2|\else_{#2}\fi}

      \providecommand{\intersected}{\ensuremath{\cap}}
      \providecommand{\united}{\ensuremath{\cup}}
      \providecommand{\meet}{\intersected}
      
      \providecommand{\join}{\united}
      
      \providecommand{\union}[1]{\ensuremath{\bigcup\nolimits_{#1}}}
      \providecommand{\intersection}[1]{\ensuremath{\bigcap\nolimits_{#1}}}

      \providecommand{\unions}[3][]{\union{#2\in{#3}\ifx|#1|\else:#1\fi}}

      \ifthenelse{\boolean{usesvjour}}{}{
        \let\vec\undefined
        \providecommand{\vec}[1]{\ensuremath{\boldsymbol{#1}}}
        \renewcommand{\vec}[1]{\ensuremath{\boldsymbol{#1}}}
      }

      \providecommand{\hatmat}[1]{\hat{\mat{#1}}}

      \providecommand{\mat}[1]{\geomat{#1}} %

      \providecommand{\olnum}[1]{\mathsfit{#1}} %
      \providecommand{\numsca}[1]{\olnum{#1}} %
      \providecommand{\numvec}[1]{{\vec{\numsca{#1}}}} %
      \providecommand{\numvecentry}[2]{\irow{\numsca{#1}}{#2}} %
      \providecommand{\numvecdotsfromto}[3]{\vecdotsfromto{\numsca{#1}}{#2}{#3}}

      \providecommand{\Prob}[1][]{\ensuremath{\operatorname{Prob}\ifx|#1|{}\else_{#1}\fi}}

      \providecommand{\pdf}[2][]{\ensuremath{\operatorname{pdf}_{#2\ifx|#1|{}\else{\conditionalto{#1}}\fi}}\xspace}

      \providecommand{\expectation}{\ensuremath{\operatorname{E}}}
      \providecommand{\EX}[1][]{\ensuremath{\expectation\ifx|#1|{}\else_{#1}\fi}}

      \providecommand{\gausskernel}[3][x]{%
        \ensuremath{
          \exp\frac{-\if#20{#1}\else(#1-\mu)\fi^2}{%
            2\if#31{}\else\powp2{#3}\fi}%
        }%
      }
      \providecommand{\gaussdistribution}[3][x]{%
        \ensuremath{\frac1{\sqrt{2\pic}\if#31{}\else#3\fi}%
          \gausskernel[#1]{#2}{#3}
        }%
      }%

      \providecommand{\boundary}{\partial}
      \providecommand{\PD}[1]{\operatorname{PD}\qpreg{#1}}
      \providecommand{\pdspace}[1]{\PD{\linspace v}}
      \providecommand{\pdmats}[2][F]{\PD{\csname#1\endcsname{#2}}}
      
      \providecommand{\SPD}{\operatorname{SPD}}
      \providecommand{\spdmats}[2][F]{\SPD(\csname#1\endcsname{#2})}

       \providecommand{\Continuous}{\ensuremath{\operatorname C}\xspace}%
       \providecommand{\Hspace}{\ensuremath{\operatorname H}\xspace}
       \providecommand{\Lebesgue}{\ensuremath{\operatorname L}\xspace}
       \providecommand{\Besovspace}{\ensuremath{\operatorname B}\xspace}

       \providecommand{\Weaklyder}{\ensuremath{\operatorname W}\xspace}
       
       \providecommand{\dual}[1]{\ensuremath{{#1}'}}
       
       \providecommand{\dualspace}[2][]{\dual{\linspace{#2}\ifx|#1|\else{_{#1}}\fi}}
       \providecommand{\bidual}[1]{\ensuremath{{#1}''}}
       \providecommand{\bidualspace}[2][]{\bidual{\linspace{#2}\ifx|#1|\else{_{#1}}\fi}}

       \providecommand{\cont}[1]{\ensuremath{\Continuous^{#1}}}

       \providecommand{\diffable}[2][]{\ensuremath{\cD\ifx|#1|\else^{#1}\fi(#2)}}
       
       \providecommand{\BV}[1]{\ensuremath{\operatorname{BV}}}

       \providecommand{\leb}[1]{\ensuremath{\Lebesgue_{#1}}}

       \providecommand{\lebloc}[1]{\ensuremath{{{\Lebesgue}^{\kern-.20em\lower .1ex\hbox{\tiny\textrm{\textup{loc}}}}_{#1}}}}
       \providecommand{\lebnorm}[3][]{\ensuremath{\Norm{#2}_{\leb{#3}\ifx|#1|{}\else(#1)\fi}}}
       \providecommand{\bes}[3][]{\ensuremath{\Besovspace^{#2}_{#3\ifx|#1|\else,#1\fi}}}

       \providecommand{\sob}[2]{\ensuremath{{\smash\Weaklyder}^{#1}_{#2}}}

       \providecommand{\sobh}[1]{\ensuremath{\Hspace^{#1}}}
       \providecommand{\vecsobh}[1]{\ensuremath{\vec\Hspace^{#1}}}
       \ProvideDocumentCommand{\hdiv}{ O{} O{}}{\vecsobh{\operatorname{div}}\ifx+#1+\else_{0|#1}\fi\ifx|#2|\else(#2)\fi}
       \providecommand{\hcurl}[1][]{\vecsobh{\operatorname{curl}}\ifx|#1|\else(#1)\fi}
       
       \providecommand{\sobhz}[2][]{\sobh{#2}_{0\ifx+#1+\else|#1\fi}}

       \providecommand{\Lip}[1][]{\ensuremath{\operatorname{Lip}}\ifx|#1|{}\else{\qp{#1}}\fi}

       \ProvideDocumentCommand{\polyring}{ O{X} O{A} }{\ring[#2][#1]}
       \ProvideDocumentCommand{\polyfield}{ O{X} O{} }{\field[#2][#1]}
       \providecommand{\polyreals}[1][]{\polyfield[][R]\ifx|#1|\else^{#1}\fi}
       \providecommand{\poly}[2][]{\ensuremath{\rP\ifx#1\else_{#1}\fi^{#2}}}

       \providecommand{\Symmatrices}[2][R]{\ensuremath{\operatorname{Sym}{(\csname#1\endcsname{#2})}}}
       
       \providecommand{\SAmatrices}[2][F]{\ensuremath{\operatorname{SA}{(\csname#1\endcsname{#2})}}}
       \providecommand{\mesh}[2][]{\ensuremath{\mathcalbf{\MakeUppercase{#2}}\ifx|#1|\else_{#1}\fi}}

      \providecommand{\convexhull}{\operatorname{Cnvx}}
      \providecommand{\simplex}{\convexhull}
      \providecommand{\crouzeixraviart}[1][1]{\operatorname{CR}\ifx|#1|{}\else{^{#1}}\fi}

      \providecommand{\sides}{\operatorname{Sides}}
      
      \providecommand{\sidesofmesh}[2][]{\sides\mesh[#1]{#2}}

      \providecommand{\linspace}[1]{\mathscript{\MakeUppercase{#1}}}
      
      \providecommand{\vecspace}{\linspace}
      \providecommand{\linop}[1]{\mathcalbf{\MakeUppercase{#1}}}

      \providecommand{\Lin}{\operatorname{Lin}}
      \providecommand{\CL}{\operatorname{CL}}
      \providecommand{\linops}[3][]{\ensuremath{\Lin\ifx|#1|\else^{#1}\fi\qp{{#2}\to{#3}}}}

      \providecommand{\clinops}[3][]{\ensuremath{\CL\ifx|#1|\else^{#1}\fi\qp{{#2}\to{#3}}}}
      
      \providecommand{\fepartition}[2][]{\mathscript{\MakeUppercase{#2}}\ifx|#1|{}\else_{#1}\fi}
      \providecommand{\fespace}[2][]{\mathbb{\MakeUppercase{#2}}\ifx|#1|{}\else_{#1}\fi}
      \providecommand{\hatfespace}[2][]{\widehat{\mathbb{\MakeUppercase{#2}}}\ifx|#1|{}\else_{#1}\fi}

      \providecommand{\vespace}[1][]{\fespace v\ifx|#1|\else_{#1}\fi}
      \providecommand{\hatvespace}[1][]{\hatfespace v\ifx|#1|\else_{#1}\fi}
      \providecommand{\fezerospace}[2][]{\ensuremath{\mathring{\fespace{#2}}\ifx|#1|{}\else_{#1}\fi}}

      \providecommand{\fe}[2][]{\ensuremath{\UCmath{#2}\ifx|#1|\else_{#1}\fi}}%

      \providecommand{\vecfe}[2][]{\ensuremath{\vec{\fe{#2}}\ifx|#1|{}\else{_{#1}}\fi}}%
      
      \providecommand{\matfe}[2][]{\ensuremath{\mat{\fe{#2}}\ifx|#1|{}\else{_{#1}}\fi}}%
      
      \providecommand{\hatmatfe}[2][]{\ensuremath{\hatmat{\UCmath{#2}}\ifx|#1|{}\else{_{#1}}\fi}}%

      \providecommand{\Forall}{\:\forall\:}
      
      \providecommand{\tfor}{\text{ for }}
      \providecommand{\Foreach}[1][]{\text{\; for each \ifx|#1|\else#1\ \fi}}%
      \providecommand{\Forsome}{\text{ for some }}

      \RequirePackage{amscd}
      \providecommand{\lt}{<}

      \providecommand{\ideq}{\equiv}
      
      \providecommand{\funk}[4][]{\ensuremath{#2:#3\ifx|#1|\else(\subseteq #1)\fi\to#4}}

      \providecommand{\dfunkmapsto}[6][]{\ensuremath{
          \begin{array}{rrcl}
            {#2}: & {#4} &  \to   & {#6}
            \\
                  & {#3} &\mapsto & {#5\text{\ #1}}
          \end{array}\quad}}

      \providecommand{\isomorphicto}{\leftrightarrows}
      \providecommand\isomorphic\isomorphicto

      \providecommand{\implies}{\ensuremath{\:\Rightarrow\:}\xspace}
      \renewcommand{\implies}{\ensuremath{\:\Rightarrow\:}\xspace}

      \providecommand{\imbedded}{{\ensuremath{\,\hookrightarrow\,}}}
      
      \providecommand{\embedsin}{\imbedded}

      \providecommand{\restriction}[2]{\left.#1\right|_{#2}}
      \renewcommand{\restriction}[2]{\left.#1\right|_{#2}}

      \providecommand{\evalat}[3][]{\qb{#2}_{\ifx|#1|{}\else#1=\fi#3}}
      \providecommand{\evaldiff}[4][]{\qb{#2}^{\ifx|#1|{}\else#1=\fi#3}_{\ifx|#1|{}\else#1=\fi#4}}

      \providecommand{\aka}[1]{(also known as {#1})\xspace}

      \providecommand{\akaindexen}[2][]{\aka{\indexen[#1]{#2}}}

      \providecommand{\CBS}{Cauchy--Bunyakovsky--Schwarz inequality\xspace}
      \providecommand{\bs}{\char '134}   %

      \providecommand{\Program}[1]{\textsf{#1}\xspace}%
      \providecommand{\Source}[1]{\nolinkurl{#1}\xspace}

      \providecommand{\matlabplot}[2][]{%
        \begin{center}
          \includegraphics[width=0.9375\linewidth,trim=64 200 64 200,clip]{#2}%
          \ifx|#1|\else\\#1\fi
        \end{center}%
      }

      \providecommand{\texcommand}[1]{\texttt{\bs{\nolinkurl{#1}}}\xspace}
      
      \providecommand{\codename}[1]{\nolinkurl{#1}\xspace}
      \providecommand\colorvar[2][a]{\colorbox{#1!6.25}{#2}}%
      \providecommand{\colorvarname}[2][a]{\colorvar[#1]{\Verb{#2}}}
      
      \providecommand{\codevarname}[1]{\colorvarname[a]{#1}}
      
      \providecommand\olco\codevarname

      \ifthenelse{\boolean{useutopia}}{%
        
      }{%
        
      }

      \providecommand{\matlab}{{\footnotesize\Program{MATLAB}}\xspace}%
      \providecommand\MATLAB\matlab

      \ProvideDocumentCommand{\codesnip}{ O{.} O{1.0} m}{%
        \newline
        \begin{minipage}{#2\linewidth}
          \lstinputlisting{#1/#3}
        \end{minipage}
      }
      \providecommand{\codeprint}[2][.]{
        \ \newline
        \begin{minipage}{\linewidth}
          \lstinputlisting{#1/#2}
          \framebox{Contents of file %
            \ifthenelse{\isundefined\pickuppath}{%
             \codename{#2}%
            }{%
              \providecommand{\fullpickuppath}{}%
              \renewcommand{\fullpickuppath}{\pickuppath/\ifx|#1|\else#1/\fi#2}%
              \href{\fullpickuppath}{\codename{#2}}%
          }}
        \end{minipage}
      }
      \providecommand{\codenoprint}[2][.]{
              \providecommand{\fullpickuppath}{}%
              \renewcommand{\fullpickuppath}{\pickuppath/\ifx|#1|\else#1/\fi#2}%
              \href{\fullpickuppath}{\codename{#2}}%
      }

      \providecommand{\indexen}[2][]{{\ifthenelse{\boolean{shownotes}}{\color b}{}#2\ifx|#1|\index{#2}\else\index{#1}\fi}}
      \providecommand{\indexemph}[2][]{\emph{\indexen[#1]{#2}}}
      \providecommand{\indexma}[2][]{{\ifthenelse{\boolean{shownotes}}{\color b}{}#2\ifx|#1|\index{\(#2\)}\else\index{<#1@\(#2\)}\fi}}

      \newenvironment{Hint}{\par\textit{Hint}.}{\par}

      \providecommand{\ListParameters}{}
      \renewcommand{\ListParameters}%
      {
      	 \setlength{\topsep}{0pt}
      	 \setlength{\leftmargin}{0pt}
               \setlength{\itemsep}{0pt}
      	 \setlength{\parsep}{0pt}
      	 \setlength{\parskip}{0pt}
               \setlength{\labelsep}{0pt}
      	 \setlength{\itemindent}{0pt}
      }
      {%
        \begin{list}%
          {}%
          {\ListParameters%
          
      }}%
      {\end{list}}
      \newcounter{tmpcounter}
      \newcounter{LetterListItem}
      \renewcommand{\theLetterListItem}{(\alph{LetterListItem})}

      \newcounter{CapitalListItem}
      \renewcommand{\theCapitalListItem}{\Alph{CapitalListItem}.}

      \newcounter{NumberListItem}
      \renewcommand{\theNumberListItem}{\arabic{NumberListItem}}
      {
      	\begin{list}%
      	{\theNumberListItem.\ }%
      	{\usecounter{NumberListItem}%
      	 \ListParameters
      	}
      }%
      {\end{list}}
      \newcounter{QuestionListItem}
      \renewcommand{\theQuestionListItem}{\textbf{Question \arabic{QuestionListItem}}}
      {
      	\begin{list}%
      	{\theQuestionListItem.\ }%
      	{\usecounter{QuestionListItem}%
      	 \ListParameters
      	}
      }%
      {\end{list}}
      \newcounter{RomanListItem}
      \renewcommand{\theRomanListItem}{(\roman{RomanListItem})}
      {
      	\begin{list}%
      	{\theRomanListItem\ }%
      	{\usecounter{RomanListItem}
      	 \ListParameters
      	}
      }%
      {\end{list}}
      \newcounter{StepsItem}
      {
      	\begin{list}%
      	{Step \theStepsItem.\ }%
      	{\usecounter{StepsItem}%
      	 \ListParameters
      	}
      }%
      {\end{list}}
      \newcounter{CasesListItem}
      \renewcommand{\theCasesListItem}{\Alph{CasesListItem}}
      {
      	\begin{list}%
      	{\emph{Case \theCasesListItem.}\ }%
      	{\usecounter{CasesListItem}%
      	 \ListParameters
      	}
      }%
      {\end{list}}
      \newcounter{NumCasesListItem}
      \renewcommand{\theNumCasesListItem}{\arabic{NumCasesListItem}}
      {
      	\begin{list}%
      	{\emph{Case \theNumCasesListItem.}\ }%
      	{\usecounter{NumCasesListItem}%
      	 \ListParameters
      	}
      }%
      {\end{list}}
      \newcounter{QAListItem}
      \renewcommand{\theQAListItem}{Q\arabic{QAListItem}:}
      {
      	\begin{list}%
      	{\theQAListItem}%
      	{\usecounter{QAListItem}
      	 \ListParameters
      	}
      }%
      {\end{list}}

      \ifthenelse{\boolean{isthesis}}{%
        \setcounter{secnumdepth}{1}%
      }{
        \setcounter{secnumdepth}{2} %
      }
      \providecommand{\ListParameters}{}
      \renewcommand{\ListParameters}
      {
      	 \setlength{\topsep}{0em}
      	 \setlength{\leftmargin}{0em}
               \setlength{\itemsep}{0ex}
      	 \setlength{\parsep}{.5ex}
      	 \setlength{\itemindent}{\labelsep}
      	 \addtolength{\itemindent}{\labelwidth}
      }

        \providecommand{\ObsName}{Remark}%
        \providecommand{\RemName}{Remark}%
        \providecommand{\NotName}{Notation}%
        \providecommand{\BFNName}{Big~Fantastic~Note}%
        \providecommand{\DefName}{Definition}%
        \providecommand{\ExaName}{Example}%
        \providecommand{\TheName}{Theorem}%
        \providecommand{\LemName}{Lemma}%
        \providecommand{\ProName}{Proposition}%
        \providecommand{\CorName}{Corollary}%
        \providecommand{\PbmName}{Problem}%
        \providecommand{\HypName}{Hypothesis}%
        \providecommand{\AlgName}{Algorithm}%
        \providecommand{\ExeName}{Exercise}%
        \providecommand{\SolName}{Solution}%
        \providecommand{\ClaName}{Claim}%
        \providecommand{\EsyName}{Essay}%
        \providecommand{\Proofname}{Proof}%
        \providecommand{\Derivename}{Derivation}%
      
      \ifthenelse{\boolean{isthesis}}{%
        \providecommand{\Thecounter}{The}
      }{%
        \providecommand{\Thecounter}{subsection}
      }
      \newcommand{\oltikzgetxy}[3]{%
        \tikz@scan@one@point\pgfutil@firstofone#1\relax
        \edef#2{\the\pgf@x}%
        \edef#3{\the\pgf@y}%
      }
      \providecommand{\pdfformat}[1]{
         \provideboolean{pdfoutput}
         \setboolean{pdfoutput}{#1}%
        \ifthenelse{\boolean{pdfoutput}}{
          \typeout{using pdf}
\makeatletter
\usepackage{pdfsync}
          \providecommand{\graphext}{pdf}
          \renewcommand{\graphext}{pdf}
          \providecommand{\graphextex}{pdf_t}
          \renewcommand{\graphextex}{pdf_t}
        }{
          \typeout{using eps}
          \RequirePackage[dvips]{graphicx,xcolor}
          \providecommand{\graphext}{eps}
          \renewcommand{\graphext}{eps}
          \providecommand{\graphextex}{eps_t}
          \renewcommand{\graphextex}{eps_t}
        }
        \RequirePackage{epsfig}
        \RequirePackage{tikz}
        \RequirePackage{rotating}
\makeatletter
        \RequirePackage{graphicx}
        \RequirePackage{xcolor}
        \provideboolean{darkcolortheme}
        \definecolor{SussexFlint}{rgb}{.00,.19,.21}
        \definecolor{SussexGrey}{rgb}{.51,.58,.49}
        \definecolor{SussexOrange}{rgb}{.94,.29,.00}
        \definecolor{SussexYellow}{rgb}{1.00,.73,.00}
        \definecolor{SussexRed}{rgb}{.94,.01,.49}
        \definecolor{SussexPurple}{rgb}{.48,.06,.44}
        \definecolor{SussexGreen}{rgb}{.00,.58,.46}
        \definecolor{OmarGreen}{rgb}{.00,.68,.36}
        \definecolor{SussexBlue}{rgb}{.00,.58,.65}
        \definecolor{OmarBlue}{rgb}{.00,.38,.65}
        \colorlet{olgreen}{SussexGreen}
        \colorlet{olblue}{OmarBlue}
        \colorlet{a}{OmarBlue}%
        \colorlet{b}{SussexOrange}
        \colorlet{c}{SussexGreen}
        \colorlet{d}{SussexPurple}%
        \colorlet{e}{SussexRed}
        \colorlet{f}{SussexYellow}
        \colorlet{g}{white}%
        \colorlet{h}{SussexGrey}%
        \colorlet{i}{black}%
        \colorlet{j}{SussexFlint}
        \colorlet{colora}{a}
        \colorlet{colorb}{b}
        \colorlet{colorc}{c}
        \colorlet{colord}{d}
        \colorlet{colore}{e}
        \colorlet{colorf}{f}
        \colorlet{colorg}{g}
        \colorlet{colorh}{h}
        \colorlet{colori}{i}
        \colorlet{colorj}{j}
        \newcommand{\mausDarkColorTheme}{
          \colorlet{a}{SussexYellow!50!yellow}
          \colorlet{b}{SussexBlue}%
          \colorlet{c}{SussexRed!50!red}
          \colorlet{d}{SussexOrange!50!yellow}
          \colorlet{e}{SussexGreen!50!green}
          \colorlet{f}{SussexPurple!50!magenta}
          \colorlet{g}{black}%
          \colorlet{h}{SussexFlint!50!black}
          \colorlet{i}{white}%
          \colorlet{j}{SussexGrey}
        }
        \ifthenelse{\boolean{darkcolortheme}}{\mausDarkColorTheme}{}
\makeatletter
      }
      \providecommand{\solution}{\textbf{\SolName.}\xspace}

      \newcounter{phantomedinput}
      \newcounter{phantombox}
      \counterwithin*{phantombox}{phantomedinput}%
      \provideboolean{showphantoms}
      \renewcommand{\thephantombox}{\Alph{phantombox}}%
      \providecommand{\phantombox}[2][]{\stepcounter{phantombox}%
        \ensuremath{
          \boxed{%
            \ifthenelse{\boolean{showphantoms}}{#2}{\phantom{#2}}
            \texttt{\tiny\ \colorbox{i!50}{\color g\thephantombox}}
          }%
        }%
      }%

      \provideboolean{hidesolution}
      \newcommand{\consolution}[2][]{
        \ifthenelse{\boolean{hidesolution}}{#1\setboolean{showphantoms}{false}}{%
          {\setboolean{showphantoms}{true}\color{i!50}\par \small {\solution}\ #2\par\ \\[5pt]}}
      }
      \provideboolean{showmarks}
      \providecommand{\showmarks}[1]{%
        \ifthenelse{%
          \boolean{showmarks}}{%
          \!\,\marginpar{%
            \tiny [$#1$ mark\ifthenelse{\equal{#1}1}{\phantom{s}}s]}%
        }{}}%

      \newcommand{\condibreak}{\ifthenelse{\boolean{hidesolution}}{\clearpage}{}}
      
      \newcommand{\solutibreak}{\ifthenelse{\boolean{hidesolution}}{}{\clearpage}}
      \newcommand{\questionly}[1]{\ifthenelse{\boolean{hidesolution}}{#1}{}}
      \newcommand{\solutionly}[1]{\ifthenelse{\boolean{hidesolution}}{}{#1}}

       \providecommand{\qeyword}[1]{\index{#1}\ifthenelse{\boolean{shownotes}}{{\tiny\color e\colorbox{e!6.25}{#1}}}{}}
       \providecommand{\pathwordbase}{file:///Users/ol20/Sokratex}
       \providecommand{\pathword}[2][]{%
         \label{#2}%
         \ifthenelse{\boolean{shownotes}}{%
           \ \\\index{#2@\tiny\codevarname{#2}}{%
             \tiny{\color f\href{\pathwordbase/#2}{\colorvarname[f]{#2}%
           }}}\\
         }{}%
       }
       \providecommand{\targword}[2][]{%
         \label{#2}%
         \ifthenelse{\boolean{shownotes}}{%
           \index{#2@\tiny\codevarname{#2}}{\ensuremath{\tiny\color d-> \href{\pathwordbase/#2}{\colorvarname[d]{#2}\ifx|#1|\else\colorvarname[d]{[#1]}\fi}}}\\
         }{}%
       }
       \providecommand{\sourceurl}[2][]{%
         \ifthenelse{\boolean{shownotes}}{{\ \\\tiny\colorbox{d!6.25}{\color d\texttt{source: \ifx|#1|\href{#2}{#1}\else\url{#2}\fi}}}}}
       \providecommand{\sourcecite}[2][]{\ifthenelse{\boolean{shownotes}}{{\ \\\tiny\colorbox{d!6.25}{\color d\texttt{source: \citet[#1]{#2}}}}}{%
       }}
       \providecommand{\conword}[2][]{\ifthenelse{\boolean{shownotes}}{#2}{#1}}
       \providecommand{\solword}[2][]{\ifthenelse{\boolean{hidesolution}}{#1}{#2}}
       \providecommand{\solghost}[1]{\ifthenelse{\boolean{showphantoms}}{#1}{\phantom{#1}}}

      \RequirePackage{lineno}
      \ifthenelse{\boolean{showchanges}}{
        \newcommand{\llabel}[1]{\hypertarget{llineno:#1}{\linelabel{#1}}}
        \newcommand{\lref}[1]{\hyperlink{llineno:#1}{\ref*{#1}}}
      }{
        \newcommand\llabel[1]{}
        \newcommand\lref[1]{}
      }
      \provideboolean{includeresponses}
      \setboolean{includeresponses}{false}
      \providecommand{\mailto}[1]{\href{mailto:#1}{\nolinkurl{#1}}}
      \provideboolean{showoldetails}
      \setboolean{showoldetails}{true}
      \providecommand{\oldetails}[2]{\ifthenelse{\boolean{showoldetails}}{#1}{#2}}
      \RequirePackage{hyphenat}
   \ifthenelse{\boolean{nohyperref}}{
     \PackageWarning{omarstyle}{package hyperref not loaded, please load manually if needed}
   }{
     \RequirePackage[obeyspaces,hyphens]{url}
     \RequirePackage[bookmarks,hypertexnames=false,debug,pdfpagelabels]{hyperref}
     \RequirePackage{bookmark}
   }

   \newtheoremstyle{plain}%
     {}%
     {}%
     {\mdseries\slshape}%
     {\parindent}%
     {\bfseries}%
     {.}%
     {.5em}%
     {}%
   
   \newtheoremstyle{note}%
     {}%
     {}%
     {}%
     {\parindent}%
     {\bfseries}%
     {.}%
     {.5em}%
     {}%
   
   \newtheoremstyle{claim}%
     {}%
     {}%
     {\mdseries\slshape}%
     {}%
     {\bfseries}%
     {}%
     {.5em}%
     {}%
   
   \newtheoremstyle{exercise}%
     {}%
     {}%
     {}%
     {}%
     {\bfseries}%
     {.}%
     {1em}%
     {}%
   
   \newtheoremstyle{break}%
     {}%
     {}%
     {}%
     {}%
     {\bfseries}%
     {.}%
     {\newline}%
     {}%
   
   \swapnumbers{
     \theoremstyle{plain}
     \ifthenelse
         {\boolean{isthesis}}
         {\newtheorem{The}{\TheName}[section]}%
         {
         }%
   {
      \theoremstyle{plain}
   
      \renewcommand{\Thecounter}{subsection}

      \newtheorem*{The*}{\TheName}
      \newtheorem*{Lem*}{\LemName}
      \newtheorem*{Pro*}{\ProName}
      \newtheorem*{Cor*}{\CorName}
      \newtheorem*{Pbm*}{\PbmName}
      \newtheorem*{Hyp*}{\HypName}
      \newtheorem*{Exe*}{\ExeName}
      \newtheorem*{Txx*}{\ExeName} %
      \newtheorem*{Con*}{Conclusion}
      \newtheorem*{Sum*}{Summary}
    }
    {
      \theoremstyle{claim}

    }
    {
      \theoremstyle{note}

      \newtheorem*{Obs*}{\ObsName}

      \newtheorem*{Def*}{\DefName}
      \newtheorem*{Exa*}{\ExaName}
      \newtheorem*{Alg*}{\AlgName}
    }
   
    {
      \theoremstyle{break}
    }
   }

   \ifthenelse{\boolean{usecleveref}}{
     \usepackage{aliascnt}
     \usepackage{cleveref}
     \providecommand{\@oltocline}{tocline}
     \newtheorem{OLThe}[subsection]{Theorem}%
     \crefname{OLThe}{theorem}{theorems}
     \newenvironment{The}[1][]{%
       \begin{OLThe}[#1]%
         \addcontentsline{toc}{subsection}{\thesubsection. \TheName\ifx|#1|\else\ - #1\fi}
     }{\end{OLThe}}
     \newaliascnt{Lem}{subsection}%
     \newtheorem{OLLem}[Lem]{Lemma}
     \crefname{OLLem}{lemma}{lemmata}
     \newenvironment{Lem}[1][]{%
       \begin{OLLem}[#1]%
         \addcontentsline{toc}{subsection}{\thesubsection. \LemName\ifx|#1|\else\ - #1\fi}
     }{\end{OLLem}}  
     \newaliascnt{Pro}{subsection}
     \newtheorem{OLPro}[Pro]{Proposition}
     \crefname{OLPro}{proposition}{propositions}
     \newenvironment{Pro}[1][]{%
       \begin{OLPro}[#1]%
         \addcontentsline{toc}{subsection}{\thesubsection. \ProName\ifx|#1|\else\ - #1\fi}
     }{\end{OLPro}}  
     \newaliascnt{Cor}{subsection}
     \newtheorem{OLCor}[Cor]{Corollary}
     \crefname{OLCor}{corollary}{corollaries}
     \newenvironment{Cor}[1][]{%
       \begin{OLCor}[#1]%
         \addcontentsline{toc}{subsection}{\thesubsection. \CorName\ifx|#1|\else\ - #1\fi}
     }{\end{OLCor}}  
     \newaliascnt{Def}{subsection}
     \theoremstyle{note}
     \newtheorem{OLDef}[Def]{Definition}
     \crefname{OLDef}{definition}{definitions}
     \newenvironment{Def}[1][]{%
       \begin{OLDef}[#1]%
         \addcontentsline{toc}{subsection}{\thesubsection. \DefName\ifx|#1|\else\ of #1\fi}
     }{\end{OLDef}}
     \newaliascnt{Obs}{subsection}
     \theoremstyle{note}
     \newtheorem{OLObs}[Obs]{Remark}
     \crefname{OLObs}{remark}{remarks}
     \newenvironment{Obs}[1][]{%
       \begin{OLObs}[#1]%
         \addcontentsline{toc}{subsection}{\thesubsection. \ObsName\ifx|#1|\else\ on #1\fi}
     }{\end{OLObs}}
     \newaliascnt{Sol}{subsection}
     \theoremstyle{note}
     \newtheorem{OLSol}[Sol]{Solution}
     \crefname{OLSol}{solution}{solutions}
     \newenvironment{Sol}[1][]{%
       \begin{OLSol}[#1]%
         \addcontentsline{toc}{subsection}{\thesubsection. \SolName\ifx|#1|\else\ on #1\fi}
     }{\end{OLSol}}
     \newaliascnt{Txx}{subsection}
     \theoremstyle{note}
     \newtheorem{OLTxx}[Txx]{\ExeName}
     \crefname{OLTxx}{exercise}{exercises}
     \newenvironment{Txx}[1][]{%
       \begin{OLTxx}[#1]%
         \addcontentsline{toc}{subsection}{\thesubsection. \ExeName\ifx|#1|\else\ - #1\fi}
     }{\end{OLTxx}}
     \newaliascnt{Exa}{subsection}
     \theoremstyle{note}
     \newtheorem{OLExa}[Exa]{\ExaName}
     \crefname{OLExa}{example}{examples}
     \newenvironment{Exa}[1][]{%
       \begin{OLExa}[#1]%
         \ifx|#1|%
         \addcontentsline{toc}{subsection}{\thesubsection. \ExaName}%
         \else%
         \addcontentsline{toc}{subsection}{\thesubsection. \ExaName\ of #1}%
         \fi
     }{\end{OLExa}}
     \newaliascnt{Alg}{subsection}
     \theoremstyle{note}
     \newtheorem{OLAlg}[Alg]{\AlgName}
     \crefname{OLExe}{exercise}{exercises}
     \newenvironment{Alg}[1][]{%
       \begin{OLAlg}[#1]%
         \ifx|#1|%
         \addcontentsline{toc}{subsection}{\thesubsection. \AlgName}%
         \else%
         \addcontentsline{toc}{subsection}{\thesubsection. \AlgName: #1}%
         \fi
     }{\end{OLAlg}}
     \newaliascnt{Pbm}{subsection}
     \theoremstyle{note}
     \newtheorem{OLPbm}[Pbm]{\PbmName}
     \crefname{OLPbm}{problem}{problems}
     \newenvironment{Pbm}[1][]{%
       \begin{OLPbm}[#1]%
         \addcontentsline{toc}{subsection}{\thesubsection. \PbmName\ifx|#1|\else\ - #1\fi}
     }{\end{OLPbm}}
     \aliascntresetthe{Lem}
     \aliascntresetthe{Pro}
     \aliascntresetthe{Def}
     \aliascntresetthe{Obs}
     \aliascntresetthe{Sol}
     \aliascntresetthe{Txx}
     \aliascntresetthe{Exa}
     \aliascntresetthe{Pbm}
     \crefname{Lem}{lemma}{lemmata}
     \crefname{Pro}{proposition}{propositions}
     \crefname{Def}{definition}{definitions}
     \crefname{Obs}{remark}{remarks}
     \crefname{Sol}{solution}{solutions}
     \crefname{Txx}{exercise}{exercises}
     \crefname{Exa}{example}{examples}
     \crefname{Pbm}{problem}{problems}
   }{
     \newenvironment{The}[1][]{%
       \ifx&#1&%
       \subsection{\TheName\xspace}%
       \else%
       \subsection[\MakeUppercase#1 theorem]{\TheName\ (#1)}%
       \fi%
       \slshape}{%
       \upshape}
     \newenvironment{Pro}[1][]{\subsection{\ProName\xspace{\ifx&#1&{}\else{ (#1)}\fi}}\slshape}{\upshape}
     \newenvironment{Lem}[1][]{\subsection{\LemName\xspace{\ifx&#1&{}\else{ (#1)}\fi}}\slshape}{\upshape}
     \newenvironment{Cor}[1][]{\subsection{\CorName\xspace{\ifx&#1&{}\else{ (#1)}\fi}}\slshape}{\upshape}
     \newenvironment{Txx}[1][]{\subsection{\ExeName\xspace{\ifx&#1&{}\else{ (#1)}\fi}}\slshape}{\upshape}
     \newenvironment{Pbm}[1][]{\subsection{\PbmName\xspace{\ifx&#1&{}\else{ (#1)}\fi}}\slshape}{\upshape}
     \newenvironment{Def}[1][]{\subsection{\DefName\xspace{\ifx&#1&{}\else{ of \indexen{#1}}\fi}}}{}
     \newenvironment{Obs}[1][]{\subsection{\ObsName\xspace{\ifx&#1&{}\else{ (#1)}\fi}}}{}
     \newenvironment{Exa}[1][]{\subsection{\ExaName\xspace{\ifx&#1&{}\else{ (#1)}\fi}}}{}
     \newenvironment{Alg}[1][]{\subsection{\AlgName\xspace{\ifx&#1&{}\else{ (#1)}\fi}}}{}
   }
   \providecommand{\qed}{\vrule height 5pt depth 0pt width 3pt}
   \providecommand{\qqed}{{\raggedright{\ \hfill\qed}}}
   
   \newcounter{passo}

   \newenvironment{Proof}[1][]%
   {\par\noindent{\bf \Proofname\ifx|#1|.\ \else\ #1.\ \fi}\setcounter{passo}{0}}%
   {\qqed\par}
   {\par\noindent{\bf \Derivename\ #1}\setcounter{passo}{0}}%
   {\qqed\par}
   \newenvironment{Proof*}[1][{}]%
   {\subsection{\Proofname\ #1}\setcounter{passo}{0}}
   {\qqed\par}

\usepackage{calc}
\usepackage{tikz}
\usetikzlibrary{calc}
\newlength\figureheight
\newlength\figurewidth
\usepackage{pgfplots}
\usepackage[authoryear]{natbib}
\usepackage{cleveref}
\crefname{equation}{}{}
\crefname{section}{\S}{\SS}
\usepackage{hyperref}
\makeatletter
\DeclareRobustCommand{\crefnosort}[1]{%
  \begingroup\@cref@sortfalse\cref{#1}\endgroup
}
\makeatother
\providecommand{\mathsfit}{\mathsf}
\usepackage{pgf}
\usepackage{pgfplots}
\usepackage{grffile}
\pgfplotsset{compat=newest}
\usetikzlibrary{plotmarks}
\usetikzlibrary{arrows.meta}
\usepgfplotslibrary{patchplots}
\makeatletter
\newboolean{longversion}
\setboolean{longversion}{true}%
\newcommand{\longonly}[2][]{%
  \ifthenelse{\boolean{longversion}}{#2}{#1}%
}
\providecommand{\field}{\rF}
\renewcommand{\field}{\rF}

\providecommand{\OR}{\tor}

\providecommand{\abs}[1]{\norm{#1}}

\providecommand{\div}{}
\renewcommand{\div}{\operatorname{div}}

\renewcommand{\mat}[1]{\vec{\MakeUppercase{#1}}}
\renewcommand{\maxi}[2]{\operatorname{max}\setpair{#1}{#2}}

\setboolean{showzeroentries}{true}

\providecommand{\drawcard}[1]{}
\renewcommand{\drawcard}[1]{\tikz{\draw[fill=a!20] (-5pt,-8pt) rectangle ++(10pt,16pt); \path (0,0) node{\textbf{#1}};}\:}

\renewcommand{\mat}[1]{\ensuremath{\vec{\MakeUppercase{#1}}}}
\providecommand{\matspace}[3][R]{\mathbb #1^{#2\times#3}}

\providecommand{\fieldmats}[2]{\matspace[F]{#1}{#2}}

\providecommand{\linop}[1]{\ensuremath{\mathcalbf{\MakeUppercase{#1}}}}
\providecommand{\opeye}[1]{\linop I}

\providecommand{\vecspace}[1]{\ensuremath{\mathscript{\MakeUppercase{#1}}}}

\renewcommand{\div}[1][]{\nabla\ifx|#1|{}\else\kern-2pt_{#1}\fi\kern-2pt\inner}

\providecommand{\inner}{\cdot}
\renewcommand{\inner}{\cdot}

\setboolean{showzeroentries}{false}
\RequirePackage{siunitx}
\providecommand{\iverson}{}
\renewcommand{\iverson}[1]{\left\llbracket #1\right\rrbracket}

\renewcommand{\d}{\operatorname d\!}

\providecommand{\attime}[1]{^{#1}}
\providecommand{\ta}[2][]{\ensuremath{\mathsfit{#2}\ifx|#1|\else\attime{#1}\fi}}
\providecommand{\tavec}[2][]{\ensuremath{\mathsfbfit{#2}\ifx|#1|\else\attime{#1}\fi}}
\providecommand{\feop}[2][\vespace]{{\uppercase{#2}}\ifx|#1|\else_{#1}\fi}%
\providecommand{\tildefeop}[2][n]{\widetilde{\fe[#1]{#2}}}
\ProvideDocumentCommand{\febas}{O{\Phi} m O{}}{\fe{#1}^{#3}_{#2}}%

\ProvideDocumentCommand{\listfebasndotsfromto}{O{\Phi} mm O{}}{%
  \listdotsfromto{\febas[#1]{#2}[#4]}{\febas[#1]{#3}[#4]}}%
\providecommand{\listfebasndots}[2][\Phi]{\listfebasndotsfromto[#1]1{#2}}%

\ProvideDocumentCommand{\febasndots}{O{\Phi} m O{}}{\qb{\listfebasndots[#1]{#2}[#3]}}

\providecommand{\iverson}[1]{\left\llbracket #1\right\rrbracket}
\renewcommand{\iverson}[1]{\left\llbracket #1\right\rrbracket}

\renewcommand{\field}{\rK}

\providecommand{\Hmat}[1]{\ensuremath{\mat{\operatorname H}(#1)}}
\providecommand{\Hmatvec}[2][]{\Hmat{\vec #2\ifx|#1|{}\else{^{#1}}\fi}}

\providecommand{\rayleigh}[2][]{\ensuremath{R\ifx&#1&{}\else_{\mat #1}\fi}\qp{\vec #2}}
\providecommand{\inner}{\cinner}
\renewcommand{\inner}{\cinner}

\renewcommand{\d}{\operatorname d}

\providecommand{\lmm}[1][]{linear \ifx|#1|multistep\else$#1$-step\fi\ method\xspace}

\providecommand{\y}[2][]{\ensuremath{y^{#2}\ifx|#1|\else_{#1}\fi}}

\providecommand{\inew}{\operatorname{new}}
\providecommand{\iold}{\operatorname{old}}
\providecommand{\imid}{\operatorname{mid}}
\providecommand{\ynew}[1][]{\y[\inew\ifx|#1|\else,#1\fi]}
\providecommand{\yold}[1][]{\y[\iold\ifx|#1|\else,#1\fi]}
\providecommand{\ymid}[1][]{\y[\imid\ifx|#1|\else,#1\fi]}

\providecommand{\vecy}[2][]{\ensuremath{\vec y^{#2}\ifx|#1|\else_{#1}\fi}}
\providecommand{\vecynew}[1][]{\vecy{\inew\ifx|#1|\else,#1\fi}}
\providecommand{\vecyold}[1][]{\vecy{\iold\ifx|#1|\else,#1\fi}}
\providecommand{\vecymid}[1][]{\vecy{\imid\ifx|#1|\else,#1\fi}}

\providecommand{\vecynplus}[1]{\ensuremath{\vec y^{n+{#1}}}}
\providecommand{\vecynp}\vecynplus

\providecommand{\tn}[1][n]{t_{#1}}

\providecommand{\tnplus}[2][n]{\tn[#1+#2]}
\providecommand{\tnp}[2][n]{\tnplus[#1]{#2}}%
\providecommand{\tnph}[1][n]{\ifx|#1|\tn[\fracl12]\else\tnp[#1]{\fracl12}\fi}
\providecommand{\tnminus}[2][n]{\tn[#1-#2]}
\providecommand{\tnm}[2][n]{\tnminus[#1]{#2}}%
\providecommand{\tnmh}[1][n]{\ifx|#1|\tn[-\fracl12]\else\tnm[#1]{\fracl12}\fi}
\providecommand{\x}{}
\renewcommand{\x}[2][]{\bgroup x_{#2}\ifx|#1|\else^{#1}\fi\egroup}

\providecommand{\xnew}[1]{\x{\inew}}

\providecommand{\timestep}[1][]{\ensuremath{h\ifx|#1|\else_{#1}\fi}}

\providecommand{\trunc}{\ensuremath{\tau}}
\providecommand{\truncation}[2][]{\ensuremath{\trunc_{#2}\ifx|#1|\else^{\text{#1}}\fi}}
\providecommand{\vectruncation}[2][]{\ensuremath{\vec\trunc_{#2}\ifx|#1|\else^{\text{#1}}\fi}}

\providecommand{\OR}{\ensuremath{\quad\operatorname{or}\quad}}

\providecommand{\magfac}{M}
\providecommand{\magnification}[2][]{\ensuremath{\magfac_{#2}\ifx|#1|\else^{\text{#1}}\fi}}

\renewcommand{\timestep}[1][]{\ensuremath{\tau\ifx|#1|\else_{#1}\fi}}
\providecommand{\tat}[1]{t_{#1}}
\renewcommand{\tat}[1]{t_{#1}}
\providecommand{\tn}[1][n]{\tat{#1}}
\providecommand{\tnpm}[2][n]{\ifx|#1|\tat{\pm#2}\else\tat{#1\pm#2}\fi}%

\providecommand{\tnmo}[1][n]{\tnminus[#1]1}
\providecommand{\tnpmo}[1][n]{\tnpm[#1]1}
\providecommand{\tnph}[1][n]{\tnplus[#1]{\fracl12}}%
\providecommand{\tnmh}[1][n]{\tnminus[#1]{\fracl12}}%
\providecommand{\tnpmh}[1][n]{\tnpm[#1]{\fracl12}}

\providecommand{\feta}[2][]{\ta[#1]{\fe{#2}}}
\providecommand{\ellop}{\linop a}
\providecommand{\abil}[3][a]{\ensuremath{{#1}\qb{#2,#3}}}
\providecommand{\aqua}[2][a]{\abil[#1]{#2}{#2}}
\providecommand{\lproj}[1][\vespace]{\feop[#1]P}

\providecommand{\discellop}[1][]{\feop[#1]a}
\providecommand{\pivot}{\linspace h}
\providecommand{\elldom}{\linspace v}
\providecommand{\ellran}{\dualspace v}

\providecommand{\inner}{\cdot}
\renewcommand{\inner}{\cdot}

\providecommand{\starof}[2][]{\ifx|#1|\omega_{#2}\else\omega_{\mesh{#1}}(#2)\fi}

\providecommand{\patchof}[2][]{\ifx|#1|\tilde\omega_{#2}\else\tilde\omega_{\mesh{#1}}(#2)\fi}
\providecommand{\cosineschemecquafn}[6][-2p]{
  \backdiff^2\ta[{#6}]{#3}
  +
  {#2}\linop{#4}\ta[{#6}]{#3}
  \ifx|#1|{-2p}\else{#1}\fi
  \linop{#4}\ta[{#6}-1]{#3}
  +
  {#2}\linop{#4}\ta[{#6}-2]{#3}
  =
  {#2}{#5}\attime {#6}
  \ifx|#1|{-2p}%
  \else%
    {#1}{#5}\attime{{#6}-1}
    +
    {#2}{#5}\attime{{#6}-2}
  \fi
}

\renewcommand{\iverson}[1]{\one_{\qb{#1}}}
\providecommand{\ellop}{}
\renewcommand{\ellop}[1][]{\linop A\ifx|#1|\else_{#1}\fi}
\providecommand{\ellopmesh}[1]{\ellop[\mesh #1]}
\renewcommand{\abil}[3][]{\duality[#1]{\ellop #2}{#3}}
\renewcommand{\aqua}[2][]{\abil[#1]{#2}{#2}}
\renewcommand{\backdiff}{\partial^-}
\providecommand{\centdiff}{\partial}
\renewcommand{\centdiff}{\partial}
\renewcommand{\forediff}{\partial^+}
\providecommand{\extforediff}{}
\renewcommand{\extforediff}[2][n]{\frac{{#2}\attime{#1+1}-{#2}\attime{#1}}{\timestep}}
\providecommand{\extcentdiff}[2][n]{}
\renewcommand{\extcentdiff}[2][n]{\frac{{#2}\attime{#1+1}-{#2}\attime{#1-1}}{2\timestep}}
\providecommand{\secodiff}{}
\renewcommand{\secodiff}[1][]{\partial^2}
\providecommand{\extsecodiff}{}
\renewcommand{\extsecodiff}[2][n]{\frac{{#2}\attime{#1+1}-2{#2}\attime{#1}+{#2}\attime{#1-1}}{\timestep^2}}

\renewcommand{\timestep}[1][]{\Delta\tn[#1]}
\providecommand{\polydeg}{}
\renewcommand{\polydeg}{k}

\providecommand{\nmth}[1][n]{#1-\fracl32}
\providecommand{\nmh}[1][n]{#1-\fracl12}
\providecommand{\nph}[1][n]{#1+\fracl12}
\providecommand{\npth}[1][n]{#1+\fracl32}

\RenewDocumentCommand{\meshsize}{O{n} O{}}{h^{#1}\ifx|#2|\else_{#2}\fi}
\providecommand{\meshsizeof}[1]{\meshsize[][#1]}
\renewcommand{\meshsizemesh}[1]{\meshsize[][\mesh{#1}]}

\providecommand{\coarse}{{\operatorname{c}}}
\providecommand{\fine}{{\operatorname{f}}}
\providecommand{\coarsemesh}[2][n]{\mesh[#1]{#2}^{\coarse}}
\providecommand{\coarsemeshsize}[1][]{\meshsize[\coarse][#1]}
\providecommand{\finemesh}[2][n]{\mesh[#1]{#2}^{\fine}}
\providecommand{\finemeshsize}[1][]{\meshsize[\fine][#1]}

\providecommand{\finespace}[1][n]{\vespace[#1]^{\fine}}
\providecommand{\coarsespace}[1][n]{\vespace[#1]^{\coarse}}

\providecommand{\passop}[1][n]{\feop[#1]\Pi}%

\providecommand{\finepass}[1][n]{\passop[#1]^{\fine}}
\providecommand{\fineproj}[1][n]{\lproj[#1]^{\fine}}

\providecommand{\projf}[1]{\feta[#1]f}
\providecommand{\extprojf}{\overline{\projf{}}}

\renewcommand{\discellop}[1][n]{\feop[#1]a}
\renewcommand{\vespace}[1][n]{\fespace[#1]v}
\providecommand{\hatvespace}[1][n-]{\hatfespace[#1\fracl12]v}

\providecommand{\ellrecop}[1][n]{\linop R\ifx|#1|\else_{#1}\fi}

\providecommand{\clementop}[1][]{\linop Q\ifx|#1|\else_{#1}\fi}

\providecommand{\ta}[2][]{{#2}\if|#1|\else^{#1}\fi}%
\renewcommand{\ta}[2][]{{#2}\if|#1|\else^{#1}\fi}%
\renewcommand{\feta}[2][n]{\ta[#1]{\fe{#2}}}

\providecommand{\timeinter}[1][n]{I_{#1}}
\providecommand{\timestaginter}[1][n+]{I_{#1\fracl12}}
\providecommand{\timehalfinter}[1][\nu]{I'_{#1}}

\providecommand{\recu}{\omega}
\renewcommand{\recu}{\omega}
\providecommand{\recv}{\psi}%
\renewcommand{\recv}{\psi}%

\providecommand{\tarecu}[1][n]{\ta[#1]\recu}
\renewcommand{\tarecu}[1][n]{\ta[#1]\recu}
\providecommand{\tarecvbase}[1][n-\frac12]{\ta[#1]\recv}
\renewcommand{\tarecvbase}[1][n-\frac12]{\ta[#1]\recv}
\providecommand{\tarecv}[1][n-]{\ta[#1\fracl12]\recv}
\renewcommand{\tarecv}[1][n-]{\ta[#1\fracl12]\recv}
\providecommand{\fetav}[1][n-]{\feta[#1\fracl12]v}

\providecommand{\error}{e}

\providecommand{\erroru}[1][]{\error_0\ifx|#1|\else^{#1}\fi}
\providecommand{\errorv}[1][]{\error_1\ifx|#1|\else^{#1\fracl12}\fi}

\providecommand{\errvec}[1][]{\vec\error\ifx|#1|\else^{#1}\fi}

\providecommand{\erc}[1][]{\sigma\ifx|#1|\else_{#1}\fi}
\providecommand{\errrecu}{\erc[0]}%
\providecommand{\errrecv}{\erc[1]}%
\providecommand{\errrec}{\vec\erc}

\providecommand{\errel}{\epsilon}
\providecommand{\errellu}{\errel_0}
\providecommand{\errellv}{\errel_1}
\providecommand{\errell}{\vec\errel}

\providecommand{\qes}{}%
\renewcommand{\qes}[2][n]{\ta[#1]\rho_1}%
\providecommand{\res}{}
\renewcommand{\res}[2][n]{\ta[#1]\rho_0}%

\providecommand{\reszero}[1][]{r_0}
\providecommand{\resone}[1][]{r_1}
\providecommand{\resvec}{\vec r}
\providecommand{\resrecu}[1][n]{\res[#1]{\recu,\recv}}

\providecommand{\resrecuexttime}[1][n]{\frac\ellop4\qb{\tarecu[#1+1]-2\tarecu[#1]+\tarecu[#1-1]}}
\providecommand{\resrecuextoperators}[1][n]{\qb{\discellop[n]-\ellrecop[n+1]\passop[n+1]\tildefeop[n]a}\feta[n]u}
\providecommand{\resrecuextmeshchange}[1][n]{\qb{\ellrecop[#1+1]\passop[#1+1]\fetav[#1-]-\tarecv[#1-]}\inverse\timestep}

\providecommand{\resrecut}{\overline{\resrecu[]}}
\providecommand{\resrecvbase}[1][n-\frac12]{\qes[#1]{\recv,\recu}}%
\providecommand{\resv}[1][n-]{%
  \ifx|#1|\qes[]{\recv,\recu}\else\qes[#1\fracl12]{\recv,\recu}\fi}
\providecommand{\resvexttime}[1][n]{%
  -\frac14\qp{\tarecv[#1+]-2\tarecv[#1-]+\tarecv[#1-1]}}

\providecommand{\resvextmeshchange}[1][n]{%
  \qb{\ellrecop[#1]\passop[#1]-\ellrecop[#1-1]}\feta[#1-1]u\inverse\timestep%
}

\providecommand{\resvt}{\overline{\resv[]}}%

\providecommand{\lintrec}[2][]{#2}
\providecommand{\slintrec}[2][]{\hat{#2}}
\providecommand{\quatrec}[2][]{\breve{#2}\ifx|#1|\else^{#1}\fi}

\providecommand{\lintrecu}[1][]{\lintrec{\recu}\ifx|#1|\else\attime{#1}\fi}
\providecommand{\slintrecu}[1][]{\slintrec{\recu}\ifx|#1|\else\attime{#1}\fi}
\providecommand{\lintrecv}[1][]{\lintrec{\recv}\ifx|#1|\else\attime{#1}\fi}
\providecommand{\slintrecv}[1][]{\slintrec{\recv}\ifx|#1|\else\attime{#1}\fi}
\providecommand{\quatrecu}[1][]{\quatrec[#1]{\recu}}
\providecommand{\quatrecv}[1][]{\ifx|#1|\quatrec[]{\recv}\else\quatrec[#1\fracl12]{\recv}\fi}

\providecommand{\potenorm}[2][]{\Norm[#1]{#2}_{\ellop}}
\DeclareMathOperator{\erg}{erg}

\renewcommand{\pivot}{\leb2(\W)}
\newcommand{\pivotnorm}[2][]{\Normon[#1]{#2}{\pivot}}
\providecommand{\ergprod}[3][\erg,\linop A]{\ltwop[#1]{#2}{#3}}
\providecommand{\ergnorm}[2][]{\Normon[#1]{#2}{\erg,\linop A}}
\providecommand{\supergnorm}[4][]{\Normsupon[#1]{#2}{#3,#4;\erg,\linop A}}
\providecommand{\lebergnorm}[5][]{\Normonleb[#1]{#2}{#3}{#4,#5;\erg,\linop A}}

\providecommand{\indellzero}[1][n]{\varepsilon_0^{#1}}
\providecommand{\indellone}[1][n]{\varepsilon_1^{#1}}
\providecommand{\indoperator}[1][n]{\alpha^{#1}}

\providecommand{\indmeshchange}[1][n]{\mu\ifx|#1|\else^{#1}\fi}
\providecommand{\indmeshchangezero}[1][n]{\indmeshchange[#1]_0}
\providecommand{\indmeshchangeone}[1][n]{\indmeshchange[#1]_1}
\providecommand{\indmeshchangelts}[1][n]{\indmeshchange[#1]_2}
\providecommand{\inddatafun}[1][n]{\delta^{#1}}

\providecommand{\indtimefunzero}[1][n]{\vartheta_0^{#1}}

\providecommand{\indtimefunone}[1][n]{\vartheta_1^{#1}}

\providecommand{\indtotnospace}[1][m]{\zeta^{#1}}
\providecommand{\abil}[3][a]{\ensuremath{{#1}\qb{#2,#3}}}
\providecommand{\aqua}[2][a]{\abil[#1]{#2}{#2}}

\renewcommand{\areaof}[2][]{\qv[#1]{#2}}
\makeindex
\author{Marcus J. Grote}
\thanks{This work was partly supported by the Swiss National Science Foundation.}
\address{Marcus J. Grote
  ---
  Department of Mathematics and Informatics,
  University of Basel,
  Spiegelgasse 1,
  4051 Basel, Switzerland}
\email{marcus.grote@unibas.ch}
\author{Omar Lakkis}
\thanks{OL likes to thank the EU Marie Curie ITN ``ModCompShock'' and the Dr Perry James (Jim) Browne Research Centre on Mathematics and its Applications at the University of Sussex for their support.}
\address{Omar Lakkis
  ---
  Department of Mathematics,
  University of Sussex,
  Brighton,
  BN1 9QH,
  England UK}
\email{lakkis.o.maths@gmail.com}
\author{Carina S. Santos}
\address{Carina Santos
  ---
  Department of Mathematics and Informatics,
  University of Basel,
  Spiegelgasse 1,
  4051 Basel, Switzerland}
\email{carina.santos@unibas.ch}
\date{\today}
\title[A posteriori error analysis of a variable mesh leapfrog method]{A posteriori error estimates for the wave equation with mesh change in the leapfrog method}
\hypersetup{
  pdfauthor={Marcus Grote and Omar Lakkis},
  pdfkeywords={wave, partial differential equation, computational, leapfrog},
  pdfsubject={numerical analysis},
  pdfcreator={Emacs 25.3.1 (Org mode 8.2.10)}}
\setboolean{showrevisions}{false}%
\setboolean{emphrevisions}{false}%
\setboolean{shownotes}{false}%
\setboolean{showchanges}{false}%
\setboolean{longversion}{true}%
\makeatletter
\@ifpackageloaded{biblatex}{
  \newbibmacro*{bbx:parunit}{%
    \ifbibliography
        {\setunit{\bibpagerefpunct}\newblock
          \usebibmacro{pageref}%
          \clearlist{pageref}%
          \setunit{\adddot\newline\nobreak}}
        {}}
  
  \renewbibmacro*{doi+eprint+url}{%
    \usebibmacro{bbx:parunit}%
    \iftoggle{bbx:doi}
             {\printfield{doi}}
             {}%
             \iftoggle{bbx:eprint}
                      {\usebibmacro{eprint}}
                      {}%
                      \iftoggle{bbx:url}
                               {\usebibmacro{url}}
                               {}}
  
  \renewbibmacro*{eprint}{%
    \usebibmacro{bbx:parunit}%
    \iffieldundef{eprinttype}
                 {\printfield{eprint}}
                 {\printfield[eprint:\strfield{eprinttype}]{eprint}}}
  
  \renewbibmacro*{url+urldate}{%
    \usebibmacro{bbx:parunit}%
    \printfield{url}%
    \iffieldundef{urlyear}
                 {}
                 {\setunit*{\addspace}%
                   \printtext[urldate]{\printurldate}}}
  
  \renewbibmacro*{url}{%
    \usebibmacro{bbx:parunit}%
    \iffieldundef{doi}{%
      \printfield{url}%
    }
  }

}{
  \RequirePackage{doi}
}

\numberwithin{equation}{section}
\makeatother
\begin{document}
\definecolor{mycolor1}{rgb}{0.00000,0.44700,0.74100}%
\definecolor{mycolor2}{rgb}{0.85000,0.32500,0.09800}%
\definecolor{mycolor3}{rgb}{0.92900,0.69400,0.12500}%
\maketitle
\begin{abstract}
  We derive a fully computable \aposteriori error estimator for a
  Galerkin finite element solution of the wave equation with explicit
  leapfrog time-stepping. Our discrete formulation accommodates both
  time evolving meshes and leapfrog based local time-stepping
  \citet{DiazGrote:09:article:Energy-Conserving}, which overcomes the
  stringent stability restriction on the time-step due to local mesh
  refinement. Thus we account for adaptive time-stepping with mesh
  change in a fully explicit time integration while retaining its
  efficiency.
  The error analysis relies on elliptic reconstructors and abstract
  grid transfer operators, which allows for use-defined elliptic error
  estimators.
  Numerical results using the elliptic Babuška--Rheinboldt estimators
  illustrate the optimal rate of convergence with mesh size of the
  \aposteriori error estimator.
\end{abstract}
\pathword{sec-introduction.tex}%
\conword\par%
\section{Introduction}
Adaptive time-step and mesh refinement are key strategies in many efficient
numerical simulation of partial differential equations. \Aposteriori
error estimates are the cornerstone of any adaptive strategy that
relies on mathematically rigorous and computable error bounds. For
elliptic problems, standard residual based \aposteriori error
estimates yield elementwise error indicators used to steer the mesh
adaptation process \citep[][and references
  therein]{AinsworthOden:00:book:A-posteriori,Verfurth:13:book:A-posteriori}.
For time-dependent, e.g., parabolic problems, \aposteriori error
estimates naturally involve a time-discretization part
\citep[e.g.,][]{ErikssonJohnson:91:article:Adaptive%
  ,Picasso:98:article:Adaptive%
  ,ChenJia:04:article:An-adaptive%
  ,LakkisMakridakis:06:article:Elliptic%
  ,AkrivisMakridakisNochetto:06:article:A-posteriori%
  ,GaspozSiebertKreuzerZiegler:19:article:A-convergent}.

For second order hyperbolic problems, such as the wave equation,
\aposteriori error estimation is less developed than in the elliptic
or parabolic case. In \citet{Johnson:93:article:Discontinuous}
\aposteriori estimates were derived for a space-time discretization of
the second order wave equation with continuous finite elements (FEM)
in space and a discontinuous Galerkin (DG) discretization in time
\citep{HulbertHughes:90:article:Space-Time-Finite}. Goal oriented
adaptive methods based on duality and hence on the solution of adjoint
problems was proposed in
\citet{BangerthRannacher:01:article:Adaptive%
  ,BangerthGroteHohenegger:04:article:Finite-Element}.
Residual
based \aposteriori error estimates with first-order implicit
time-stepping were developed in \citet{BernardiSuli:05:article:Time},
and also in \citet{Adjerid:02:article:A-posteriori,
  Adjerid:06:article:A-posteriori} using spatial bi-$p$ FEM on
rectangular grids. More recently, \aposteriori error estimates in the
$\leb2(0,T;\sobh1(\W))$ norm were derived for semi-discrete
formulations with anisotropic mesh refinement using elliptic
reconstructions \citep{Picasso:10:article:Numerical,%
  GoryninaLozinskiPicasso:19:article:An-Easily-Computable}.
\changes{%
  These works consider either semi-discrete formulations (continuous
  in time), or fully discrete formulations based on implicit rather
  than explicit time integration.

  To the best of our knowledge,
  \citet{GeorgoulisLakkisMakridakisVirtanen:16:article:A-Posteriori}
  derived
}%
the first \aposteriori error estimate for semi-discrete formulations
(continuous in space) for second-order wave equations, discretized in
time using two-step Newmark method (also known as cosine method)
family, which includes the explicit leapfrog method herein addressed.
\changes{%
  In \cite{GoryninaLozinskiPicasso:19:article:An-Easily-Computable}
  provide a time-reconstruction with a fixed spatial mesh with
  continuous FEM in space and Newmark (including leapfrog) methods in
  time.  More recently,
  \citet{Chaumont-FreletErn:24:techreport:Damped-Energy-Norm}
  established error estimates for the fully discrete leapfrog method
  in time and continuous FEM in space in ``damped energy norms'',
  previously introduced in a semidiscrete setting by
  \citet{Chaumont-Frelet:23:article:Asymptotically-Constant-Free}; the
  work is developed under fixed mesh and fixed time-step assumptions.
}%

\changes{%
  In the design of adaptive methods, \aposteriori error estimates of
  fully discrete formulations in time-dependent problems should
  include
}%
the added effect on the error due to mesh change from one time-step to the next;
in fact, some of the above cited
works address that particular issue for parabolic problems \citep[see
  also][]{Dupont:82:article:Mesh,LakkisPryer:12:article:Gradient}.  
Both for accuracy and efficiency, it is indeed expected, often even required, 
from any adaptive method to locally adapt and change
the mesh repeatedly during the entire simulation. 
Although quantification of mesh-change error in second order
hyperbolic problems is less studied, a notable result in this
direction was provided by
\citet{KarakashianMakridakis:05:article:Convergence} in an \apriori
setting. 

While local mesh refinement is certainly key to any efficient
numerical method, it also hampers any explicit time-stepping method
due to the stringent CFL stability condition which imposes a tiny
time-step across the entire computational domain. By taking smaller
time-steps, but only inside those smaller elements due to local mesh
refinement, local time-stepping (LTS) methods overcome that major
bottleneck without sacrificing explicit time-stepping. For
  this reason the main objective of this paper is the derivation of
fully discrete \aposteriori error estimates in the presence of locally
refined meshes that may vary in time together with the associated
local time-stepping procedure
\citep{DiazGrote:09:article:Energy-Conserving}.

\changes{%
  We focus on conforming methods (continuous FEM) in space and on
  the \aposteriori estimation of time-maximum potential and kinetic
  energy in space of the
  discretization error.  We stress that \emph{mesh change} and
  \emph{local time-stepping} are rigorously taken into account with a
  view to developing adaptive explicit wave
  solvers. }%
A difficulty in establishing error bounds for the leapfrog method is
related to its symplectic nature where the velocity and the state are
intimately related; the error can be analyzed by considering these
quantities on two \indexen{staggered} time-grids. That is, a primal
time grid for the position variable, $u$, and a grid offset by half a
time-step for its time derivative $v=\pdt u$.  The spatial
discretization uses (continuous) $\sobh1(\W)$-conforming finite
elements of arbitrary polynomial degree. Moreover, our estimates allow
for a changing mesh and also accommodate the use of leapfrog based LTS
methods as proposed by \citet{DiazGrote:09:article:Energy-Conserving}
and \citet{GroteMitkova:10:article:Explicit-Local} \citep[see
  also][and the references therein]{%
  GroteMehlinSauter:18:article:Convergence,%
  GroteMichelSauter:21:article:Stabilized,%
  CarleHochbruck:22:article:Error-Analysis}.
\changes{%
  Our fully discrete \aposteriori error-estimates for the wave
  equation thus pave the way to adaptive space and time (with mesh
  change) solvers while retaining efficiency of the fully explicit
  nature of leapfrog.  Note that, the leapfrog's property of
  preserving the discrete energy \citep[as described
    in][IX.8]{HairerLubichWanner:10:book:Geometric} on a time-constant
  spatial mesh is impossible to maintain exactly for time-varying
  spatial meshes. Nonetheless, thanks to the error estimators bounding
  the error's full energy norm, an adaptive algorithm that uses them
  will be able to quantify the error in the energy norm, and thus
  approximate the \emph{exact energy} within an arbitrarily small
  tolerance over exponentially long integration times.
}%

The rest of our paper is structured as follows. In
\cref{sec:notation-and-scheme}, we present the problem, introduce
notation and state the fully discrete Galerkin formulation of the wave
equation using $\sobh1$-conforming finite elements and the the
leapfrog method in time. With a careful choice of finite element
spaces and their bases as to make degrees of freedom coincide with
certain quadrature nodes, these methods allow for high-order
\indexen{mass lumping} in space, which means that the numerical method
is fully explicit, efficient and easily parallelizable
\citep{CohenJolyRobertsTordjman:01:article:Higher-Order}.  The
proposed approach accommodates for both time evolving meshes (under a
reasonable \indexen{mesh compatibility} condition, briefly discussed
in \cref{sec:residual-esimators} and leapfrog based local
time-stepping \cite{DiazGrote:09:article:Energy-Conserving}.  Starting
from the time discrete numerical solutions in possibly varying FE
spaces, in \cref{sec:reconstruction} we recall the corresponding
elliptic and time reconstructions together with the associated
residuals. In \cref{sec:aposteriori-error-analysis} those space-time
reconstructions lead to a continuous error equation akin to the wave
equation reformulated as a first-order system.
\changes{%
  The energy-based estimators are fully computable energy-norm error
  bounds, formed as the sum in time of mean-square sum of local error
  indicators accounting for spatial discretization, time
  discretization and local time-stepping and mesh-change indicators.
  Finally, in \cref{sec:numerical-results}, we consider a
  one-dimensional Gaussian pulse on a \emph{locally refined} and
  \emph{time-varying} mesh and compare the true error with the
  \aposteriori estimates, as we progressively refine the mesh.
}%
\pathword{sec-notation-and-scheme.tex}%
\conword\par%
\section{The wave equation and its discrete counterpart}
\label{sec:notation-and-scheme}
Here we define the model problem and functional analytic framework
(\crefrange{sec:wave-equation}{sec:energy-norms}), the leapfrog
discretization in time and space
(\crefrange{sec:time-discretization}{sec:finite-element-spaces}), and
  the associated local time-stepping on variable meshes
  (\crefrange{sec:fine-coarse-Du-splitting}{sec:time-varying-mesh}).
\subsection{The wave equation}
\label{sec:wave-equation}
We consider the \indexen{wave equation} for the unknown $u(\vec x,t)$ with
$\vec x\in\W$ (on a Lipschitz domain in $\R d$) and time $0\leq t\leq T$,
with forcing $f(\vec x,t)$
\begin{equation}
  \label{eq:def:wave-equation:concrete}
  \pdtt u(\vec x,t)-\divof{c(\vec x)^2\grad u(\vec x,t)}=f(\vec x,t)
  \text{ for }
  \vec x\in\W
  \tand
  t\in\opclinter0T
\end{equation}
coupled with \indexen{Dirichlet--Neumann boundary conditions}
\begin{equation}
  \restriction u{\Gamma_0}(t)=0
  \tand
  \restriction{\normalderto\W u(t)}{\boundary\W\take\Gamma_0}=0
  \text{ for }
  t\in\opclinter0T
\end{equation}
(where $u(t)$ is short for $u(\cdot,t)$)
and the initial conditions
\begin{equation}
  u(0)=u_0\tand \pdt u(0)=v_0
\end{equation}
for given functions $u_0,v_0$.  We assume the \indexen{Dirichlet
  boundary} $\Gamma_0\subseteq\boundary\W$ to have a strictly positive
measure%
. The \indexen{wave velocity field} $c$ is a constant or a
function in $\leb\infty(\W)$ which satisfies
\begin{equation}
  0<c_\flat\leq c(\vec x)\leq c_\sharp\text{ in $\W$}
\end{equation}
for two constants $c_\flat,c_\sharp\in\reals$. The \indexen{forcing
  term} $f$ is a space-time function.
\subsection{Functional spaces and PDE in abstract form}
\label{sec:functional-space-and-pde}
\changes{%
  We will denote by $\elldom:=\sobhz[\Gamma_0]1(\W)$, the space of
  Sobolev square-summable functions of order one which vanish (in the
  sense of traces) on $\Gamma_0\subseteq\boundary\W$, a strictly
  positive measure, $\areaof{\Gamma_0}>0$.
  We also write $\ellran$ for the (topological) dual space of $\elldom$, with
  $\pivot$ being the pivot space to obtain the Gelfand triple structure
  \begin{equation}
    \elldom\embedsin\pivot\embedsin\ellran.
  \end{equation}
  The inner products of two elements, say $\phi$
  and $\psi$ in $\pivot$ and $\elldom$ are respectively indicated by
  \begin{equation}
    \indexma[<,>]{\ltwop\phi\psi}
    :=\int_\W\phi\psi
    \tand
    \indexma[<|>]{\ltwopon\phi\psi\elldom}
    :=\int_\W\grad\phi\inner\grad\psi,
  \end{equation}
  for any $\phi,\psi$ for which the integrals (and gradients in the
  second case) make sense.  The duality pairing, thought of as a
  bilinear form on $\ellran\times\elldom$ is indicated with
  \begin{equation}
    \indexma[<>]{\duality g\phi}\Foreach g\in\ellran,\phi\in\elldom.
  \end{equation}
  
  The spatial differential \indexen{operator} appearing in
  \cref{eq:def:wave-equation:concrete} will be denoted by
  \begin{equation}
    \ellop\phi(\vec x)
    =
    \divof{c(\vec x)^2\grad\phi(\vec x)},
  \end{equation}
  for twice differentiable functions $\phi$, and we consider its
  extension $\funk{\indexma[A]\ellop}{\linspace v}{\dualspace v}$ via
  its \indexen{bilinear form} representation, whereby for each
  $\phi\in\linspace v$, $\ellop\phi$ is the unique member of $\dualspace v$ such that
  \begin{equation}
    \duality{\ellop\phi}\psi
    =
    \ltwop{c\grad\phi}{\grad\psi}
    =
    \int_\W c(\vec x)^2\grad\phi(\vec x)\psi(\vec x)\d\vec x
    \Foreach
    \psi\in\elldom.
  \end{equation}
  The operator $\linop A$ and, equivalently, the associated bilinear form
  $\abil\cdot\cdot$, is symmetric,
  \begin{equation}
    \abil\phi\psi=\abil\psi\phi,
  \end{equation}
}%
and satisfies the Lax--Milgram theorem
assumptions,
\begin{equation}
  \label{eqn:ellop-Lax-Milgram}
  c_\flat^2
  \Normon\phi\elldom^2
  \leq
  \abil\phi\phi
  \tand
  \abil\phi\psi
  \leq
  c_\sharp^2
  \Normon\phi\elldom\Normon\psi\elldom
\end{equation}
for all $\phi,\psi\in\elldom$%
.
\changes{%
  With this notation we rewrite the wave problem
  \cref{eq:def:wave-equation:concrete},
  as that of finding $\funk u{\opclinter0T}\linspace v$ such that
  $\pdtt u\in\leb2(0,T;\ellran)$ and
}%
\begin{equation}
  \label{eq:def:wave-equation:abstract}
  \begin{gathered}
    \pdtt u+\ellop u=f
    \text{ on }\opclinter0T
    ,
    \\
    u(0)=u_0
    \tand
    \pdt u(0)=v_0
    .
  \end{gathered}
\end{equation}
\changes{%
  We often rewrite equation \cref{eq:def:wave-equation:abstract}
  as a first order system:
}%
\begin{equation}
  \label{eq:def:wave-equation:abstract-system}
  \pdt\discolvectwo{u}v
  +
  \dismattwo 0{-1}{\ellop}0
  \discolvectwo {u}v
  =
  \discolvectwo 0f.
\end{equation}
\changes{%
  \begin{Obs}[regularity of data]
    Our analysis applies under general conditions; but for
    simplicity's sake we deal with more specific ones.  In particular,
    we observe the following:
    \begin{enumerate}[(i)\ ]
    \item
      The boundary conditions in \cref{eq:def:wave-equation:abstract-system}
      do not need to be homogeneous.
    \item
      The source term $f$ may be taken in $\leb2(0,T;\pivot)$, or
      even in $\leb2(0,T;\ellran)$.
    \item
      It is possible to have a time-dependent wave-velocity
      $c(\vec{x},t)$ instead of the time-constant one $c(\vec x)$.  Since
      time-varying spatial meshes are allowed, a necessary requirement
      for a fully adaptive method, in our analysis the discrete
      elliptic operator defined below
      (\cref{sec:discrete-elliptic-operator}) could in fact be
      time-dependent, even for a time-constant $c$.
    \end{enumerate}
  \end{Obs}
}%
\subsection{Energy norms}
\label{sec:energy-norms}
The function $u$ satisfying the wave
equation \cref{eq:def:wave-equation:concrete} has, associated to it,
the \indexemph{wave-energy} \akaindexen{total energy} which is the sum
of its \indexen{kinetic energy} and \indexen{potential energy}:
\begin{equation}
  \frac12
  \Normonleb{\pdt u(t)}2\W^2
  +
  \frac12
  \Normonleb{c\grad u(t)}2\W^2
  .
\end{equation}
\changes{%
  The associated \indexen{potential energy norm} for any $\phi\in\elldom$
  is given by
  \begin{equation}
    \label{eq:def:potential-energy-norm}
    \indexma[(A]{\potenorm\phi}
    :=
    \powsqrt{\aqua\phi}
    =
    \Normonleb{c\grad\phi}2\W
    ,
  \end{equation}
  which thanks to the boundary conditions in
  \cref{eq:def:wave-equation:concrete}, or the assumptions on $\ellop$
  in \cref{sec:functional-space-and-pde}, is equivalent as a norm to
  the norm of $\elldom$.
  In %
  the special case of $c\ideq1$, $\ellop$ coincides with the
  (positive) Laplace operator, $-\div\grad$ and the potential energy
  norm coincide with the seminorm $\Normonleb{\grad\phi}2\W$, for
  $\phi$ members of $\sobhz[\Gamma_0]1(\W)$; this in fact a norm owing to
  the Poincaré--Friedrichs inequality and $\areaof{\Gamma_0}>0$.
}%
Introduce the \indexen{wave-energy scalar product}, as the bilinear
form
\begin{equation}
  \ergprod{\vec\phi}{\vec\chi}
  :=
  \duality{\linop A\phi_0}{\chi_0}
  +
  \ltwop{\phi_1}{\chi_1}
  \Foreach
  \vec\phi=\colvecitwoz\phi,
  \vec\chi=\colvecitwoz\chi
  \in
  \elldom\times\pivot
  .
\end{equation}
The corresponding full
\indexen{wave-energy norm} will be denoted by
\begin{equation}
  \indexma[(erg]{\ergnorm{\vec\phi}}
  :=
  \ergprod{\vec\phi}{\vec\phi}^{\fracl12}.
\end{equation}
In terms of the elliptic and mean-square norms we have
\begin{equation}
  \ergnorm{\vec\phi}^2
  =
  \potenorm{\vecentry\phi0}^2
  +
  \Normon{\vecentry\phi1}\pivot^2
  \sim
  \Normon{\vecentry\phi0}\elldom^2
  +
  \Normon{\vecentry\phi1}\pivot^2
  .
\end{equation}
\subsection{Time discretization}
\label{sec:time-discretization}
We discretize time with a \indexen{global time grid}
which a standard uniform partition of the time interval with integer indices
defined as
\begin{equation}
  0=\tn[0]<\tn[1]<\dotsb<\tn[N]=T,
  \text{ where }
  \tn:=n\timestep.
\end{equation}
We will use also the corresponding \indexen{staggered time grid},
whose nodes are the midpoints of the global time grid's nodes,
\begin{equation}
  \tnmh[]<\tnph[]<\dotsb<\tnmh[N]<T<\tnph[N]
  \text{ where }
  \tnpmh:=\frac{\tnpmo+\tn}2=\tn\pm\half\timestep.
\end{equation}
The corresponding \indexen{time intervals} are denoted by
\begin{equation}
  \timeinter :=
  \clinter{\tnmo}{\tn}
  \tand
  \timestaginter
  :=
  \clinter{\tnmh}{\tnph}.
\end{equation}
These two mutually ``dual'' grids play a central role in the analysis
and we will use piecewise polynomial time-basis-functions defined on them.

The simplest such time-basis-functions are \emph{two families}
of piecewise linear (i.e., piecewise affine) functions
\begin{equation}
  \setsifromto\ell n0N
  \tand
  \setsfromto{\ell_{\nmh}}n0{N+1}
\end{equation}
where for each integer or half-integer time index
\begin{equation}
  \nu=-\fracl12,0,\fracl12,1,\dotsc,N,N+\fracl12,
\end{equation}
$\indexma[lnu]{\ell_\nu(t)}$ is the piecewise linear (in fact, affine)
function in $t$ satisfying
\begin{equation}
  \ell_\nu(\tat \nu)=1
  \tand
  \ell_\nu(\tat \nu+k\timestep)=0\Foreach\text{ integer }k\neq0.
\end{equation}
We will occasionally use the \indexen{time half-intervals}
\begin{equation}
  \timehalfinter:=\clinter{\tnmh[\nu]}{\tn[\nu]},
  \text{ for }
  \nu\integerbetween{-\fracl12}{N+\fracl12}.
\end{equation}

Note that the integer-indexed $\sets{\ell_n}n$, constitute a partition
of unity on $\clinter 0T$ while the half-integer-indexed
$\sets{\ell_{\nmh}}n$, constitute a partition of
unity on the interval $\clinter{-\timestep/2}{T+\timestep/2}$.

\begin{figure}
  \begin{tikzpicture}[scale=.75]
    \providecommand{\oltikzyunit}{4}
    \providecommand{\oltikzyminoffset}{-1}
    \providecommand{\namenmto}{n-2}
    \providecommand{\namenmth}{n-\fracl32}
    \providecommand{\namenmo}{n-1}
    \providecommand{\namenmh}{n-\fracl12}
    \providecommand{\namen}{n}
    \providecommand{\namenph}{n+\fracl12}
    \providecommand{\namenpo}{n+1}
    \providecommand{\namenpth}{n+\fracl32}
    \providecommand{\namenpto}{n+2}
    \providecommand{\oltikzponmto}{0}
    \providecommand{\oltikzponmth}{1}
    \providecommand{\oltikzponmo}{2}
    \providecommand{\oltikzponmh}{3}
    \providecommand{\oltikzpon}{4}
    \providecommand{\oltikzponph}{5}
    \providecommand{\oltikzponpo}{6}
    \providecommand{\oltikzponpth}{7}
    \providecommand{\oltikzponpto}{8}
    \providecommand{\oltikztimestep}{1}
    \draw[-stealth] (-2,0) node[left]{$0$}
    --(10,0) node[above]{$t$};
    \draw[color=i!25] (-2,\oltikzyunit) node[left]{\color i$1$}
    --++(12,0);
    \draw[color=i!25] (-2,.125*\oltikzyunit) node[left]{\color i$\fracl18$}
    --++(12,0);
    \draw[-stealth] (-2,0)--++(0,1.25*\oltikzyunit);
    \foreach \p in {nmto,nmo,n,npo,npto}{
      \coordinate (\p) at (\csname oltikzpo\p\endcsname,0);
      \draw[line width=.125pt,color=a!12.5!g] ($(\p)+(0,\oltikzyminoffset)$)--++(0,7);
      \draw[line width=.5pt] ($(\p)+(0,\oltikzyminoffset)$)--++(0,-.125)
       node[below]{$\color a t_{\csname name\p\endcsname}$}
       ;
    }
    \foreach \p in {nmth,nmh,nph,npth}{%
      \coordinate (\p) at (\csname oltikzpo\p\endcsname,0);
      \draw[line width=.125pt,color=b!12.5!g] ($(\p)+(0,\oltikzyminoffset)$)--++(0,7);
      \draw[line width=.5pt] ($(\p)+(0,\oltikzyminoffset)$)--++(0,-.125)
      node[below]{$\color{b}t_{\csname name\p\endcsname}$}
      ;
    }
    \providecommand{\drawoneell}[1][a]{
      \fill[color=#1!50,opacity=.2]
      (\csname oltikzpo\p\endcsname,0)
      coordinate (A)
      --
      (\csname oltikzpo\q\endcsname,\oltikzyunit)
      coordinate (B)
      --(\csname oltikzpo\r\endcsname,0)
      coordinate (C)
      ;
      \draw[color=#1,line width=1pt]
      (-2,0)
      --(A)
      --(B)
      node[above]{$\ell_{\csname name\q\endcsname}$}
      --(C)
      --(9.75,0);
    }
    \foreach \p/\q/\r in {nmto/nmo/n,nmo/n/npo,n/npo/npto}{
      \drawoneell[a]
    }
    \foreach \p/\q/\r in {nmth/nmh/nph,nmh/nph/npth}{
      \drawoneell[b]
    }
    \providecommand{\drawonequad}[1][e]{
      \coordinate (B) at (\csname oltikzpo\q\endcsname,.16*\oltikzyunit);
      \fill[color=#1!50,opacity=.2]
      (\csname oltikzpo\p\endcsname,0) coordinate (A)
      ..
      controls(B)%
      ..
      (\csname oltikzpo\r\endcsname,0) coordinate(C);
      \draw[line width=2pt,color=#1]
      (A)
      ..
      controls (B)
      ..
      (C)
      node[pos=0.5,above] {$q_{\csname name\q\endcsname}$}
      ;
    }
    \foreach \p/\q/\r in {nmth/nmo/nmh,nmh/n/nph,nph/npo/npth}{
      \drawonequad[f]
    }
    \foreach \p/\q/\r in {nmo/nmh/n,n/nph/npo}{
      \drawonequad[d]
    }
  \end{tikzpicture}
  \caption{Schematic description of the linear and quadratic time basis functions, $\ell_\nu$ and $q_\nu$, for some values of $\nu$.}
  \label{fig:time-basis-functions}
\end{figure}
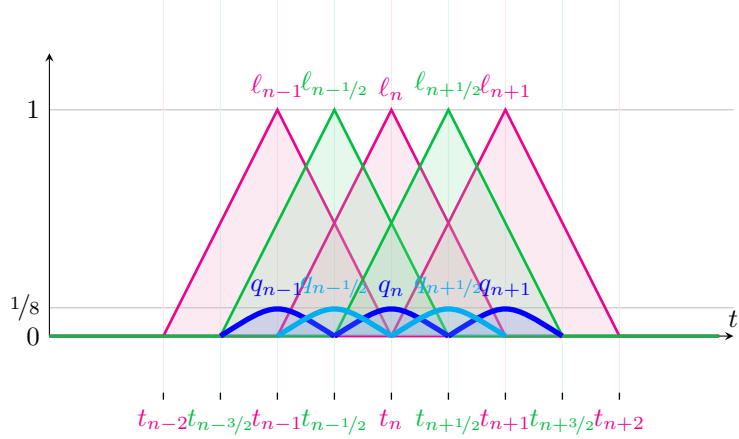

We will also use the following \indexen{quadratic bubble}
$q_\nu(t)$, defined as the positive part of the degree
$2$ polynomial in $t$ which vanishes at $\tat{\nu\pm\fracl12}$
and takes maximum $\fracl18$ at $\tat\nu$:
\begin{equation}
  \label{eq:def:quadratic-perturbation-polynomial}
  \indexma[q n]{q_\nu(t)}
  :=
  \frac{\qp{t-\tn[\nu-\fracl12]}\qp{\tn[\nu+\fracl12]-t}}{2\qppow{\timestep}2}
  \iverson{2\abs{t-\tat\nu}>\timestep}
  \text{ for }
  \nu=0,\fracl12,1,\dotsc,N-\fracl12,N.
\end{equation}
A graphic description of these functions is reported in
\cref{fig:time-basis-functions}.

For all pointwise functions $\funk\varphi J\reals$, for some interval $J$
containing time-grid points we use the shorthand
\begin{equation}
  \label{eq:def:time-grid-evaluation-shorthand}
  \varphi\attime \nu:=\varphi(\tat \nu)
  \Foreach
  \nu\betweenonetwoend{-\fracl12}0{N+\fracl12}
  .
\end{equation}
Conversely and consistently, given a sequence
$\seqsfromto{\ta[n]\phi}n0N$ defined on the integer-time grid
(respectively on the staggered time grid $\seqsfromto{\ta[\nmh]\phi}n0{N+1}$)
we will denote by $\phi(t)$ its
\indexen{continuous piecewise linear interpolation} in time
whereby
\begin{equation}
  \phi(t)
  :=
  \sum_\nu\ta[\nu]{\phi}\ell_\nu(t),
  \quad
  \text{ i.e., }
  \phi(t)
  :=
  \ta[\nu-1]\phi
  \ell_{\nu-1}(t)
  +
  \ta[\nu]\phi
  \ell_{\nu}(t)
  \text{ for }\tat{\nu-1}\leq t\leq\tat{\nu}
  .
\end{equation}
Furthermore, the \indexen{forward difference in time} of such sequence at $\tat\nu$
with
\begin{equation}
  {\indexma[d +]\forediff}\ta[\nu]\phi
  :=
  \extforediff[\nu]{\ta\phi}
\end{equation}
the \indexen{centered difference in time} at $\tat\nu$ with
\begin{equation}
  {\indexma[d 0]\centdiff}\ta[\nu]\phi
  :=
  \extcentdiff[\nu]{\ta\phi}
\end{equation}
and
the \indexen{centered second difference in time} at $\tat\nu$ with
\begin{equation}
  \label{eq:def:centered-second-difference}
        {\indexma[d 2]\secodiff}{\ta[\nu]\phi}
        :=
        \extsecodiff[\nu]{\ta\phi}
        .
\end{equation}
These difference operators need sequences defined on only one (or
both) of the two grids.
\subsection{Finite element spaces}
\label{sec:finite-element-spaces}
To each $\tn$, $n\integerbetween0N$, we associate a spatial mesh
$\indexma[M]{\mesh[n]m}$ made up of polytopal finite elements $K\in\mesh[n]m$
with flat sides grouped in a set $\indexma[Sides Mn]{\sidesofmesh[n]m}$.
The corresponding piecewise constant
\indexen{mesh-size} function
\begin{equation}
  \indexma[hn]{\meshsize[n]}(\vec x)
  :=
  \diam\qgroup{
    \intersection{\vec x\in K\in\mesh[n]m}\closure K
  }
  ;
\end{equation}
and we write
\begin{equation}
  \indexma[hE]{\meshsize[][E]}\text{ for the constant }
  \restriction{\meshsize}E
  \Foreach E\in\mesh[n]m\join\sidesofmesh[n]m.
\end{equation}
For some fixed \indexen{polynomial degree} $\indexma[k]{\polydeg}\in\naturals$ and
each $n\integerbetween0N$, we associate to the
mesh $\mesh[n]m$ the \indexen{finite element space}
\begin{equation}
  \indexma[Vn]{\vespace[n]}:=\poly\polydeg(\mesh[n]m)\meet\cont0(\W)
\end{equation}
and a corresponding finite element basis of \indexen{degrees of freedom}
\index{<Phi n m@$\febas m[n]$} 
\begin{equation}
  \febasndots{M_n}[n]\text{ where }M_n:=\dim\vespace[n].
\end{equation}
We will also use the corresponding \indexen{finite element nodes} $\vec z^n_m$
for $m\integerbetween1{M_n}$.
With this notation in mind, we can introduce the \indexen{space-pass operators}
\begin{equation}
  \funk{\indexma[Pi n]{\passop[n]}}{\cont0(\W)}{\vespace[n]}
  \text{ such that }
  \passop[n] v(\vec x)
  :=
  \sumifromto m1{M_n}\febas m[n] v(\vec z_m).
\end{equation}
Note that the choice of $\passop$ is user dependent, it could be
the Lagrange interpolator or a $\leb2(\W)$ projection, for example.

We also use the $\leb2$-projector
\begin{equation}
  \dfunkmapsto{\indexma[P n]{\lproj[n]}}g{\ellran}{\lproj[n]g}{\vespace[n]}
  \text{ where }
  \ltwop{\lproj[n] g}{\fe\phi}
  =
  \duality{g}{\fe\phi}
  \Foreach\fe\phi\in\vespace[n].
\end{equation}
\subsection{Fine and coarse degrees of freedom splitting}
\label{sec:fine-coarse-dof-splitting}
Each mesh $\mesh[n]m$ has two types of elements
\indexen{coarse} and \indexen{fine},
$\mesh[n]m=\coarsemesh m\join\finemesh m$, where
\begin{equation}
  K\in\coarsemesh m
  \implies
  \meshsize[][K]>\theta\max_{L\in\mesh[n]m}\meshsize[][L]
  \tand
  \finemesh m:=\mesh[n]m\take\coarsemesh m.
\end{equation}
for a ``user defined'' \indexen{fine--coarse threshold}
$\theta\in\opinter01$.  For example, in the simplest situation where a
generic element has size either $h$ or $h/2$, the coarse mesh has all
its elements of size $h$, while the fine mesh contains all of those of
size $h/2$ and their neighbors of size $h$.

We define a degree of freedom $\febas m[n]$ to be \indexen{fine} if
and only if its support intersects at least one element in the fine
mesh $\finemesh m$, otherwise it is \indexen{coarse} and let
$\finespace$ and $\coarsespace$ respectively be the subspaces
respectively spanned by the fine and coarse degrees of freedom. We
have thus that $\vespace[n]=\finespace\oplus\coarsespace$ and assuming
the indices are ordered into fine-first from $\listdotsfromto
1{M_n^{\fine}}$, for some integer $M_n^{\fine}\leq M_n$, and
coarse-last $\listdotsfromto{M_n^{\fine}+1}{M_n}$ every finite element
function $\fe v$ in $\vespace[n]$ can be written as
\begin{equation}
    {\fe V}{(\vec x)}
    =
    \qp{
      \sumifromto m1{M_n^{\fine}}
      +
      \sumifromto m{M_n^{\fine}+1}{M_n}
    }
    \febas m[n](\vec x)
  {\numvecentry vm}
\end{equation}
for a suitable vector
${\numvec v}={\numvecdotsfromto v1{M_n}}\in\R{M_n}.$

Similarly to \cref{sec:finite-element-spaces}
we define the \indexen{fine-mesh interpolator}
\begin{equation}
  \funk{\indexma[Pi n f]{\finepass}}{\cont0(\W)}{\finespace[n]}
\end{equation}
through relation
\begin{equation}
  \finepass\fe v
  :=
  \sumifromto m1{M_n^{\fine}}\febas m[n]\numvecentry vm
  \Foreach
  \fe v\in\vespace[n],
\end{equation}
and the \indexen{fine-mesh $\leb2$-projector}
\begin{equation}
  \funk{\indexma[P n f]{\fineproj}}{\ellran}{\finespace[n]}
\end{equation}
through relation
\begin{equation}
  \ltwop{\fineproj g}{\fe\phi}
  :=
  \duality g{\fe\phi}
  \Foreach
  \fe\phi\in\finespace[n]
  .
\end{equation}
With adaptive methods  in mind, we allow for the case where $\mesh[n]m$ (and
thus $\vespace$) changes with time, under the \indexen{mesh
  compatibility conditions}, which implies that at each point of the
domain either $\mesh[n-1]m$ is a refinement of $\mesh[n]m$ or
conversely, as explained in
\citet{LakkisMakridakis:06:article:Elliptic,LakkisPryer:12:article:Gradient}.
\subsection{Discrete elliptic operators and source approximation}
\label{sec:discrete-elliptic-operator}
For each $n$ we introduce the corresponding \indexen{discrete elliptic operator}
\begin{equation}
  \dfunkmapsto{%
    \indexma[A n]{\discellop[n]}%
  }{%
    \phi%
  }{%
    \elldom%
  }{%
    \discellop\phi:
    \ltwop{\discellop\phi}{\fe\chi}
    =
    \duality{\ellop\phi}{\fe\chi}
    \Forall\fe\chi\in\vespace[n]%
  }{%
    \vespace[n]}
  ,
\end{equation}
\indexen{local time-stepping discrete elliptic operator}
\begin{equation}
\label{eq:ltsoperator}
  \indexma[A n tilde]{\tildefeop[n]A}
  :=
  \feop[n]A-\frac{\timestep^2}{16}
  \feop[n]A\finepass\feop[n]A
\end{equation}
and the \indexen{source approximation}
\begin{equation}
  \indexma[Fn]{\projf n}
  :=
  \begin{cases}
    \lproj[n]f(\tn)
    &
    \text{ if %
      $f\in\cont0(\timestaginter;\ellran)$}
    \\
    \frac1{\timestep}\int_{\tnmh}^{\tnph}\lproj[n] f(t)\d t
    &
    \text{ if $f$ is discontinuous in time but in $\leb2(\timestaginter;\ellran)$}.
  \end{cases}
\end{equation}

The particular instance of $\tildefeop[n]A$ in \eqref{eq:ltsoperator}
corresponds to the simplest situation with two local time-steps of
size $\timestep/2$ each for each global time-step of size $\timestep$.
By letting $\tildefeop[n]A$ denote a generic perturbed bilinear
form induced by local time-stepping, our analysis inherently
encompasses situations with different coarse-to-fine time-step ratios,
too, which may even change from one locally refined region to
another across a single mesh.
\changes{%
  It also includes the ``stabilized'' version of
  LTS \citep{GroteMichelSauter:21:article:Stabilized,CarleHochbruck:22:article:Error-Analysis},
  an even the more general situation of
}%
a hierarchy of locally refined regions, each associated with its own
local time-step \citep{DiazGrote:15:article:Multi-Level-Explicit}.
\subsection{Local time-stepping}
The \indexen{leapfrog-based local time-stepping}
for time-invariant finite
element spaces, i.e.,
$\vespace[n]=\vespace[]$, $\lproj[n]=\lproj[]$, $\passop[n]=\passop[]$
and $\tildefeop[n]a=\tildefeop[]a$ for all $n$, consists in
finding a sequence $\listudotsfromto{\feta[]u}0N$ such
that
\begin{equation}
  \label{eq:leapfrog-scheme}
  \begin{aligned}
    \feta[0]u
    &
    :=\lproj[] u_0,
    \\
    \feta[1]u
    &
    :=\feta[0]u+\lproj[] v_0\timestep
    +
    \qp{\projf 0-\tildefeop[]A\feta[0]u}
    \frac{\timestep^2}2
    ,
    \\
    \feta[n+1]u
    &
    :=
    2\feta[n]u
    -
    \feta[n-1]u
    +
    \qp{
      \projf n-\tildefeop[]a\feta[n]u
    }
    \timestep^2
    \Foreach n\geq1
    ,
  \end{aligned}
\end{equation}
where the latter is equivalent to $\feta[n+1]u$
satisfying
\begin{equation}
  \secodiff{\feta u}
  +
  \tildefeop[]A\feta[n]u
  =
  \projf n
  \Foreach n\geq1.
\end{equation}

The two-step method \cref{eq:leapfrog-scheme} may be rewritten as
single-step method in the system form by introducing an auxiliary's
\indexen{discrete velocity}
\begin{equation}
  \indexma[Vn+]{\feta[n+\fracl12]v}
  :=
  \forediff\feta[n]u
  =
  \frac{\feta[n+1]u-\feta[n]u}{\timestep}
  \tfor 0\leq n<N
  ,
\end{equation}
which implies
\begin{equation}
  \forediff\feta[\nmh]v
  =
  \projf n-\tildefeop[]a\feta[n]u
  \tfor 0\leq n<N
  .
\end{equation}
This is equivalent to
\begin{align}
  \feta[n+\fracl12]v
  -
  \feta[n-\fracl12]v
  &
  =
  \qp{\projf n-\tildefeop a\feta[n]u}\timestep
  \tfor 0\leq n<N
  .
  \intertext{%
    By requiring the discrete velocities to average to the projected initial velocity,
  }
  \feta[\fracl12]v
  +
  \feta[-\fracl12]v
  &
  =
  2\lproj[]v_0,
\end{align}
we deduce the following \indexen{local time-stepping leapfrog scheme in system form}
on a fixed mesh:
\begin{equation}
  \label{eq:LTSLF:system-form:invariant-fespace}
  \begin{aligned}
    \feta[-\frac12]v
    &
    :=
    \lproj[]v_0
    -
    \qp{\projf 0-\tildefeop[]a\feta[0]u}
    \frac\timestep2
    ,
    &
    \feta[0]u
    &
    :=
    \lproj[]u_0
    ,
    &&
    \text{initially}
    \\
    \feta[\nph]v
    &
    :=
    \feta[\nmh]v
    +
    \qp{\projf n-\tildefeop[]A\feta[n]u}
    \timestep
    ,
    &
    \feta[n+1]u
    &
    :=
    \feta[n]u
    +
    \feta[\nph]v\timestep
    &&
    \tfor0\leq n<N
    .
  \end{aligned}
\end{equation}
\subsection{Time-varying mesh}
\label{sec:time-varying-mesh}
We now extend system~\cref{eq:LTSLF:system-form:invariant-fespace} to
cover the case of time-varying meshes and the corresponding finite
element spaces. So $\vespace[n-1]$ and $\vespace[n]$ may differ for
some (or all) $n\integerbetween1N$.  It is important to take care of
this aspect in an \aposteriori analysis as the associated adaptive
strategies may require time-varying meshes and thus time-varying
spaces.  In this case, looking at the case of a system first, we look
for a double sequence
$\pair{\feta[n]u}{\feta[\nmh]v}\in{\vespace[n]\times\vespace[n]}$, for
$n\integerbetween0N$ such that
\begin{equation}
  \label{eq:LTSLF:system-form:variable-fespace}
    \begin{aligned}
      \indexma[V-1/2]{\feta[-\frac12]v}
      &
      :=
      \lproj[0]v_0
      -
      \qp{\projf 0-\tildefeop[0]a\feta[0]u}
      \frac\timestep2
      \\
      \indexma[U0]{\feta[0]u}
      &
      :=
      \lproj[0]u_0
      ,
      \\
      \indexma[Vn+1/2]{\feta[\nph]v}
      &
      :=
      \passop[n+1]\qb{%
        \feta[\nmh]v
        +
        \qp{\projf n-\tildefeop[n]A\feta[n]u}
        \timestep
      }
      ,
      \\
      \indexma[Un+1]{\feta[n+1]u}
      &
      :=
      \passop[n+1]
      \feta[n]u
      +
      \feta[\nph]v\timestep
      \tfor0\leq n<N
      .
    \end{aligned}
\end{equation}
The equivalent time-varying finite element space two-step leapfrog scheme is
\begin{equation}
  \label{eq:LTSLF:two-step-form:variable-fespace}
  \begin{split}
    \feta[0]u
    &
    :=\lproj[0]u_0
    \\
    \feta[1]u
    &
    :=
    \passop[1]
    \qb{
      \feta[0]u
      +
      \qp{
        \lproj[0]v_0
        +
        \qp{
          \projf 0
          -
          \tildefeop[0]a
          \feta[0]u
        }\timestep
      }\timestep
    }
    \\
    \feta[n+1]u
    &
    :=
    \passop[n+1]\qb{
      2\feta[n]u
      -
      \passop[n]\feta[n-1]u
      +
      \qp{
        \projf n
        -
        \tildefeop[n]a
        \feta[n]u
      }\timestep^2
    }
    \tfor
    n\integerbetween1N
   .
  \end{split}
\end{equation}
\pathword{sec-reconstruction.tex}%
\conword\par%
\section{Reconstruction}
\label{sec:reconstruction}
Here we recall the concepts of elliptic reconstruction in
\cref{def:elliptic-reconstruction} and the associated elliptic error
estimator functionals in \cref{def:elliptic-estimators}. In
\cref{def:wave-residuals} we then introduce the residuals associated
with the discrete time-dependent wave equation
\eqref{eq:LTSLF:two-step-form:variable-fespace}.  In \cref{def:time-reconstructions}
we recall the time-reconstruction tools
from \citet{GeorgoulisLakkisMakridakisVirtanen:16:article:A-Posteriori},
which play a central role in our analysis, and outline their main properties
in
\cref{obs:staggered-interpolation},
lemmata
\ref{lem:quadratic-time-reconstruction-interpolates-at-nodes}%
--
\ref{lem:piecewise-linear-time-reconstruction-residual}
and \cref{the:full-time-recontstruction-residual}.%
\begin{Def}[elliptic reconstruction]
  \label{def:elliptic-reconstruction}
  For each $n\integerbetween0N$, recalling the definition of introduce
  the associated \indexemph{elliptic reconstructor} $\ellrecop$
  associated to the corresponding discrete elliptic operator $\discellop$
  (and finite element space $\vespace$)
  as follows
  \begin{equation}
    \label{eq:def:elliptic-reconstructor}
    \dfunkmapsto[.]{
      \indexma[Rn]{\ellrecop}
    }{
      \phi
    }{
      \elldom
    }{
      \ellrecop \phi:=\inverse{\ellop}\discellop{\phi}}
    \elldom
  \end{equation}
  We consider, throughout the paper, the following elliptic reconstructions
  \begin{equation}
    \tarecu:=\ellrecop\feta[n]u
    \tand
    \tarecv:=\ellrecop\fetav.
  \end{equation}

  In other words
  $\tarecu$ is the unique solution in $\elldom$
  of the elliptic BVP
  \begin{equation}
    \label{eq:elliptic-reconstruction-BVP:tarecu}
    \ellop\tarecu = \feop[n]A \feta[n]u.
  \end{equation}
  The same holds for $\tarecv$ with $\feop[n]A\fetav$ on the right-hand side
  of \cref{eq:elliptic-reconstruction-BVP:tarecu}.
  \begin{Def}[elliptic error estimators]
    \label{def:elliptic-estimators}
    We will \emph{assume throughout the analysis} in
    \crefrange{sec:reconstruction}{sec:aposteriori-error-analysis}, and we shall give
    concrete examples in, that we
    have access to \indexemph{\aposteriori error estimator functional}
    $\indexma[E]{\cE}$ such that
    \begin{equation}
      \label{eq:def:elliptic-esimators-assumption}
      \Normonspace{\tarecu-\feta u}z
      \leq
      \cE[\feta u,\vespace,\linspace z]
    \end{equation}
    where $\linspace z$ is one of $\elldom$, $\ellran$, $\ellop$ or
    $\pivot$. In \cref{sec:residual-esimators}, we describe regarding
    the estimator functionals $\cE$ in the context of residual
    Babuška--Rheinboldt estimators, and for the details we refer to
    specialized texts, such as
    \citet{Verfurth:13:book:A-posteriori%
      ,AinsworthOden:00:book:A-posteriori%
      ,Braess:07:book:Finite,BraessPillweinSchoberl:09:article:Equilibrated}.
  \end{Def}
\end{Def}
\begin{Def}[residuals]
  \label{def:wave-residuals}
  Define the following \indexemph{residuals}
  \begin{equation}
    \label{eq:def:residuals}
    \begin{split}
      \indexma[R]{\resrecu}
      &
      := %
      \resrecuexttime
      +
      \resrecuextoperators
      \\
      &\phantom{:=}
      +
      \resrecuextmeshchange
      ,
      \\
      \indexma[Q]{\resv}
      &
      :=
      \resvexttime
      +
      \resvextmeshchange
    \end{split}
  \end{equation}
  foreach $n\integerbetween1{N-1}$,
  and their (discontinuous) piecewise constant extensions:%
  \begin{equation}
    \label{eq:def:time-extended-residuals}
    \indexma[Rb]{\resrecut(t)}
    :=
    \sumifromto n0N \resrecu\charfun{\timestaginter}(t)
    \tand
    \indexma[Qb]{\resvt(t)}
    :=
    \sumifromto n0N \resv\charfun{\timeinter}(t).
  \end{equation}
  We will see that both residuals are either fully computable discrete objects
  or bounded by \aposteriori estimators of elliptic type.  In particular,
  we note the alternative expression
  \begin{equation}
    \begin{split}
      \label{eq:def:alt:displacement-residual}
      \resrecu
      &
      =
      \frac14\qp{
        \discellop[n+1]\feta[n+1]u
        -
        2\discellop[n]\feta[n]u
        +
        \discellop[n-1]\feta[n-1]u
      }
      +
      \resrecuextoperators
      \\
      &\phantom=
      +
      \resrecuextmeshchange
      ,
    \end{split}
  \end{equation}
  which means that this residual is in fact fully computable.
\end{Def}
\begin{Def}[time-reconstructions]
  \label{def:time-reconstructions}  
  Respectively define the \indexemph{primal piecewise linear time-reconstructions}
  of $\seqsfromto{\lintrecu[n]}n0N$ and $\seqsfromto{\lintrecv[\nmh]}n0{N+1}$
    with
  \begin{equation}
    \label{eq:def:linear-time-reconstructions}
    \indexma[omega (t)]{\lintrecu(t)}
    :=
    \sumifromto n0N\tarecu\ell_n(t)
    \quand
    \indexma[V (t)]{\lintrecv(t)}
    :=
    \sumifromto n0N\fetav\ell_{n-\fracl12}(t),
  \end{equation}
  where the functions $\ell_\nu$, $\nu=-\fracl12,0,\dotsc,N,N+\fracl12$,
  are defined in \cref{sec:time-discretization}.
  
  Next interpolate both these time-reconstructions, again as
  piecewise linear functions, albeit on the opposite time-grid (with a
  ``hat'' accent as mnemonic)
  \begin{equation}
    \label{eq:def:staggered-linear-time-reconstructions}
    \indexma[omega hat]{\slintrecu(t)}
    :=
    \sumifromto n0N
    \lintrecu[n-\fracl12]
    \ell_{n-\fracl12}(t)
    \quand
    \indexma[V hat]{\slintrecv(t)}
    :=
    \sumifromto n0N
    \lintrecv[n]
    \ell_n(t).
  \end{equation}
  As a result (and recalling our convention $\phi\attime\nu:=\phi(\tat\nu)$ for any
  $\phi(t)$ continuous in $t$) we have
  \begin{equation}
    \label{obs:staggered-interpolation}
    \begin{aligned}
      \slintrecu[n]
      &
      =
      \frac12\qp{\lintrecu[n-\fracl12]+\lintrecu[n+\fracl12]}
      \\
      \slintrecv[n-\fracl12]
      &
      =
      \frac12\qp{\lintrecv[n-1]+\lintrecv[n]}
      =
      \frac14\qp{\tarecvbase[\nmth]+2\tarecvbase[\nmh]+\tarecvbase[\nph]}
      .
    \end{aligned}
  \end{equation}
  \begin{figure}
   \providecommand{\scalefactor}{1.0}
   \begin{tikzpicture}[scale=\scalefactor]%
     \providecommand{\oltikzyminoffset}{-4}
     \providecommand{\namenmto}{n-2}
     \providecommand{\namenmth}{n-\fracl32}
     \providecommand{\namenmo}{n-1}
     \providecommand{\namenmh}{n-\fracl12}
     \providecommand{\namen}{n}
     \providecommand{\namenph}{n+\fracl12}
     \providecommand{\namenpo}{n+1}
     \providecommand{\namenpth}{n+\fracl32}
     \providecommand{\namenpto}{n+2}
     \providecommand{\oltikzponmto}{0}
     \providecommand{\oltikzponmth}{1}
     \providecommand{\oltikzponmo}{2}
     \providecommand{\oltikzponmh}{3}
     \providecommand{\oltikzpon}{4}
     \providecommand{\oltikzponph}{5}
     \providecommand{\oltikzponpo}{6}
     \providecommand{\oltikzponpth}{7}
     \providecommand{\oltikzponpto}{8}
     \providecommand{\oltikztimestep}{1}
     \providecommand{\oltikzwnmto}{4.75}
     \providecommand{\oltikzwnmo}{4.50}
     \providecommand{\oltikzwn}{4.00}
     \providecommand{\oltikzwnpo}{3.00}
     \providecommand{\oltikzwnpto}{1.00}
     \providecommand{\oltikzvnmth}{-0.25}
     \providecommand{\oltikzvnmh}{-0.5}
     \providecommand{\oltikzvnph}{-1.00}
     \providecommand{\oltikzvnpth}{-2.00}
     \providecommand{\oltikzrvnmth}{0.25}
     \providecommand{\oltikzrvnmh}{0.375}
     \providecommand{\oltikzrvnph}{0.5}
     \providecommand{\oltikzrvnpth}{-0.125}
     \draw[-stealth] (-2,0)--(10,0) node[above]{$t$};
     \foreach \p in {nmto,nmo,n,npo,npto}{
       \coordinate (\p) at (\csname oltikzpo\p\endcsname,0);
       \coordinate (pow\p) at (%
       \csname oltikzpo\p\endcsname,%
       \csname oltikzw\p\endcsname);
       \draw[line width=.125pt,color=a!12.5!g] ($(\p)+(0,\oltikzyminoffset)$)--++(0,7);
       \draw[line width=.5pt] ($(\p)+(0,\oltikzyminoffset)$)--++(0,-.125)
       node[below]{$\color a t_{\csname name\p\endcsname}$}
       ;
     }
     \foreach \p in {nmth,nmh,nph,npth}{
       \coordinate (\p) at (\csname oltikzpo\p\endcsname,0);
       \coordinate (pov\p) at (
       \csname oltikzpo\p\endcsname,
       \csname oltikzv\p\endcsname);
       \draw[line width=.125pt,color=b!12.5!g] ($(\p)+(0,\oltikzyminoffset)$)--++(0,7);
       \draw[line width=.5pt] ($(\p)+(0,\oltikzyminoffset)$)--++(0,-.125)
       node[below]{$\color{b}t_{\csname name\p\endcsname}$}
       ;
       \draw[color=a!50!d] ($(\p)+(-\oltikztimestep,\csname oltikzrv\p\endcsname)$)--++(2*\oltikztimestep,0) node[pos=0.5,above]{$\resrecvbase[\csname name\p\endcsname]$};
     }
     \foreach \p in {nmto,nmo,n,npo,npto}{
       \draw[fill=a] (pow\p) circle(1pt) node[anchor=south west]{%
         $\color a\tarecu[\csname name\p\endcsname]$};
     }
     \draw[color=a] (pownmto)--(pownmo)--(pown)--(pownpo)--(pownpto) node[below]{$\lintrecu$};
     \foreach \p in {nmth,nmh,nph,npth}{
       \draw[fill=i] (pov\p) circle(1pt) node[anchor=south west]{\color b$\tarecvbase[\csname name\p\endcsname]$};
     }
     \draw[color=b] (povnmth)--(povnmh)--(povnph)--(povnpth) node[below]{$\lintrecv$};
     \foreach \pbeg\pend\pmid in {%
       nmto/nmo/nmth,%
       nmo/n/nmh,%
       n/npo/nph,%
       npo/npto/npth%
     }{%
       \coordinate (pow\pmid) at ($(pow\pbeg)!0.50!(pow\pend)$);
     }
     \foreach \p in {nmth,nmh,nph,npth}{
       \draw[fill=c!50!a] (pow\p) circle(1pt);
     }
     \draw[color=c!50!a,line width=.75pt] (pownmth)--(pownmh)--(pownph)--(pownpth) %
     node[above]{$\slintrecu$};
     \foreach \p\q\r in {nmo/nmh/n,n/nph/npo}{
       \coordinate (tempo) at ($(pow\q)-.5*(0,\csname oltikzv\q\endcsname)$);
       \draw[color=c,line width=.75pt] (pow\p)..controls(tempo) ..(pow\r)
       coordinate[pos=0.75] (tmp)
       ;
     }
     \node[color=c,above] at (tmp) {$\quatrecu$};
     \foreach \pbeg\pend\pmid in {%
       nmth/nmh/nmo,%
       nmh/nph/n,%
       nph/npth/npo%
     }{
       \coordinate (pov\pmid) at ($(pov\pbeg)!0.50!(pov\pend)$);
     }
     \foreach \p in {nmo,n,npo}{
       \draw[fill=d!50!b] (pov\p) circle(1pt);
     }
     \draw[color=d!50!b,line width=.75pt] (povnmo)--(povn)--(povnpo) node[below]{$\slintrecv$};
   \end{tikzpicture}
    \caption{
      \label{fig:time-reconstructions}
      A schematic illustration of the time-reconstructions and cognate
      time-functions.  The values are only for graphing purposes and
      do not reflect actual ones.
    }
  \end{figure}
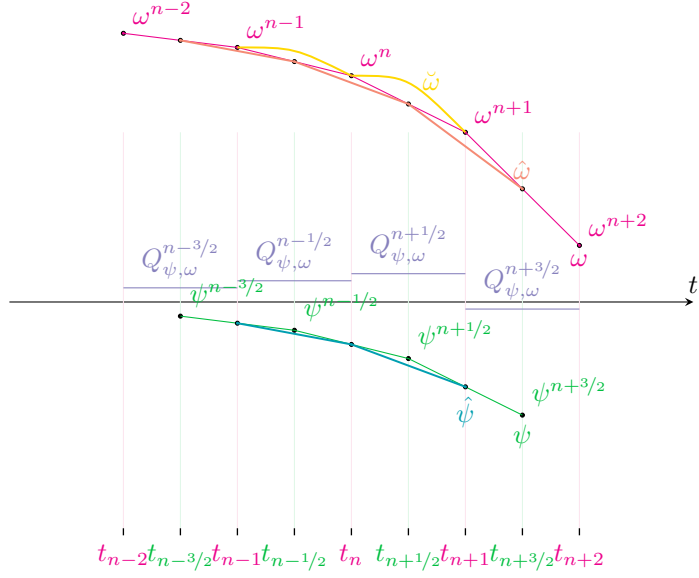
  For each $n\integerbetween1N$ we now can define the following
  \indexemph{quadratic time-reconstructions}
  \begin{equation}
    \label{eq:def:quadratic-time-reconstructions}
    \begin{split}
      \indexma[omega breve]{\quatrecu(t)}
      &
      :=
      \tarecu[n-1]%
      +\int_{\tnmo}^t\slintrecv(s)\d s
      +(t-\tnmo)\resv%
      ,%
      t\in\timeinter
      \\
      \indexma[V breve]{\quatrecv(t)}
      &
      :=
      \tarecvbase[\nmh]
      -
      \int_{\tnmh}^t
      \ellop \slintrecu(s)
      \d s
      \\
      &\phantom{:=}
      +(t-\tnmh)
      \qp{
        \ellrecop[n+1] \passop[n+1]
        \projf n
        +
        \resrecu
      }
      ,%
      t\in\timeinter[n+\fracl12]
    \end{split}
  \end{equation}
\end{Def}
\begin{Lem}[quadratic time-reconstructions interpolate at nodes]
  \label{lem:quadratic-time-reconstruction-interpolates-at-nodes}
  The quadratic displacement reconstrucion $\quatrecu$ defined in
  \cref{eq:def:quadratic-time-reconstructions} is a continuous
  piecewise quadratic in time funtion that interpolates the original
  values $\tarecu$ at the points $\ltidotsfromto t0N$. Similarly
  the quadratic velocity time-reconstruction $\quatrecv$ interpolates the
  values $\fetav$ at the staggered points $\ltidotsfromto
  t{-\fracl12}{N-\fracl12}$.
\end{Lem}
\begin{Proof}
  To see this, note that $\quatrecu(t_n^+)=\tarecu[n]$ follows immediately from the
  definition, while
  \begin{equation}
    \begin{split}
      \quatrecu(t_n^-)
      &
      =
      \tarecu[n-1]
      +
      \slintrecv[n-1]
      \int_{\tnmo}^{\tn}
      \ell_{n-1}(s)\d s
      +
      \slintrecv[n]
      \int_{\tnmo}^{\tn}
      \ell_n(s)\d s
      +
      \resv\timestep
      \\
      &
      =
      \tarecu[n-1]
      +
      \frac14\qp{
        \qp{\tarecvbase[\nmth]+\tarecvbase[\nmh]}
        +
        \qp{\tarecvbase[\nmh]+\tarecvbase[\nph]}
      }\timestep
      \\
      &
      \phantom=
      +\timestep\bigg(
      \resvexttime
      \\
      &\phantom{=+\timestep\bigg(}
      +
      \left[ \ellrecop[n-\frac{1}{2}] \fetav[n-] - \tarecvbase[\nmh] \right]
      +
      \resvextmeshchange
      \bigg)
      \\
      &
      =
      \tarecu[n-1]
      +
      \ellrecop[n]
      \qb{\fetav[n-]
        \timestep
        +
        \passop[n]
        \feta[n-1]u
      }
      -
      \ellrecop[n-1]
      \feta[n-1]u
      \\
      &
      =
      \ellrecop[n]\feta[n]u
      =
      \tarecu[n]
      .
    \end{split}
  \end{equation}
  Similarly $\quatrecv(\tnph^+)= \tarecvbase[\nph] $ is immediate for the
  integral in definition \cref{eq:def:quadratic-time-reconstructions} is $0$,
  while the same definition also implies
  \begin{equation}
    \begin{split}
      \quatrecv(\tnph^-)
      &
      =
     \tarecvbase[\nmh]
      -
      \ellop\slintrecu[\nmh]
      \int_{\tnmh}^{\tnph}
      \ell_{\nmh}(s)
      \d s
      \\
      &
      \phantom=
      -
      \ellop\slintrecu[\nph]
      \int_{\tnmh}^{\tnph}
      \ell_{\nph}(s)
      \d s
      +
      \qp{
        \ellrecop[n+1] \passop[n+1]
        \projf n
        +
        \resrecu
      }
      \timestep
      \\
      &
      =
      \tarecvbase[\nmh]
      -
      \frac\ellop4
      \qp{
        \tarecu[n-1]+\tarecu[n]
        +
        \tarecu[n]+\tarecu[n+1]       
      }\timestep
      \\
      &
      \phantom=
      +
      \ellrecop[n+1] \passop[n+1]
      \projf n\timestep
      +
      \frac14\qp{
        \discellop[n+1]\feta[n+1]u
        -
        2\discellop[n]\feta[n]u
        +
        \discellop[n-1]\feta[n-1]u
      }
      \timestep
      \\
      &
      \phantom=      
      +
      \resrecuextoperators\timestep
      +
      \resrecuextmeshchange[n]
      \\
      &
      =
      \ellrecop[n+1] \passop[n+1]
      \qb{
        \feta[\nmh]v
        +\qp{
          \projf n
          -
          \tildefeop[n]a\feta[n]u
        }
        \timestep
      }
      =
      \tarecvbase[\nph]
      .
    \end{split}
  \end{equation}
\end{Proof}
\begin{Lem}[quadratic time-reconstruction residual]
  \label{lem:quadratic-time-reconstruction-residual}
  Recalling the quadratic time-functions $q_\nu$ defined in
  \cref{sec:time-discretization}, let $n\integerbetween1{N-1}$, if
  $\tnmh\leq t\leq\tnph$ then
  \begin{equation}
    \label{eq:quadratic-time-reconstruction-residual:V}
    \quatrecv(t)
    -
    \lintrecv(t)
    =
    \frac{\discellop[n+1]\feta[n+1]u-\discellop[n-1]\feta[n-1]u}2
    q_n(t)\timestep
    \\
    =
    \centdiff\qb{\discellop[n]\feta[n]u}q_n(t)
    \qppow{\timestep}2    
  \end{equation}
  and if $\tnmo\leq t\leq\tn$
  then
  \begin{equation}
    \label{eq:quadratic-time-reconstruction-residual:recu}
    \quatrecu(t)
    -
    \lintrecu(t)
    =
    \frac{\tarecvbase[\nmth]-\tarecvbase[\nph]}2
    q_{\nmh}(t)\timestep
    =
    -
    \centdiff{\tarecvbase[\nmh]}q_{\nmh}(t)
    \qppow{\timestep}2
    .
  \end{equation}
\end{Lem}
\begin{Proof}
  Suppose $\tnmh\leq t\leq\tnph$, then by definition
  \cref{eq:def:quadratic-time-reconstructions} we have
  \begin{equation}
    {\quatrecv(t)}
      :=
      \tarecvbase[\nmh]
      -
      \int_{\tnmh}^t
      \ellop \slintrecu(s)
      \d s
      +(t-\tnmh)
      \qp{
        \ellrecop[n+1] \passop[n+1]
        \projf n
        +
        \resrecu
      }    
  \end{equation}
  where by \cref{eq:def:staggered-linear-time-reconstructions}
  and the fact that {$\ell_{\nmh}(s)+\ell_{\nph}(s)=1$},
  we obtain, for $\tnmh\leq s\leq t$, that
  \begin{equation}
    \begin{split}
      \slintrecu(s)
      &
      =
      \lintrecu[\nmh]\ell_{\nmh}(s)+\lintrecu[\nph]\ell_{\nph}(s)
      \\
      &
      =
      \frac{\tarecu[n-1]+\tarecu[n]}2
      \ell_{\nmh}(s)
      +
      \frac{\tarecu[n]+\tarecu[n+1]}2
      \ell_{\nph}(s)
      \\
      &
      =
      \frac{\tarecu[n-1]}2
      \ell_{\nmh}(s)
      +
      \frac{\tarecu[n]}2
      +
      \frac{\tarecu[n+1]}2
      \ell_{\nph}(s)
      ,
    \end{split}
  \end{equation}
  and thus, recalling \cref{eq:def:residuals}, we get
  \begin{equation}
    \begin{split}
      \resrecu
      -
      \ellop\slintrecu(s)
      &
      =
      \resrecuextmeshchange
      \\
      &\phantom=
      +
      \resrecuextoperators
      \\
      &
      \phantom=
      +
      \resrecuexttime
      \\
      &
      \phantom=
      -
      \ellop\qb{%
        \frac{\tarecu[n-1]}2
        \ell_{\nmh}(s)
        +
        \frac{\tarecu[n]}2
        +
        \frac{\tarecu[n+1]}2
        \ell_{\nph}(s)
      }.
    \end{split}
  \end{equation}
  Noting that 
  $\ellop\tarecu[n]=\discellop[n]\feta[n]u$
  we see that
  \begin{equation}
    \begin{split}
      \resrecu
      -
      \ellop\slintrecu(s)
      &
      =
      \resrecuextmeshchange
      -
      \ellrecop[n+1]\passop[n+1]
      \tildefeop[n]a\feta[n]u
      \\
      &
      \phantom=
      +
      \frac\ellop2
      \qb{%
        \tarecu[n+1]
        \qp{\frac12-\ell_{\nph}(s)}
        +
        \tarecu[n-1]
        \qp{\frac12-\ell_{\nmh}(s)}
      }
      .
    \end{split}
  \end{equation}
  To simplify further, we see that for our choice of $s$ we have
  \begin{equation}
    \ell_{\nmh}(s)+\ell_{\nph}(s)=1
  \end{equation}
  and thus
  \begin{equation}
    \qgroup{\frac12-\ell_{\nph}(s)}
    =
    -
    \qp{\frac12-\ell_{\nmh}(s)}
    =
    \tilde\ell_n(s),
  \end{equation}
  where for $\tnmh\leq s\leq\tnph$ we define
  \begin{equation}
    \indexma[l n]{\tilde\ell_n(s)}
    :=
    \frac{\tn-s}{\timestep}
    =
    \begin{cases}
      \ell_{n-1}(s)
      &
      \tfor
      \tnmh\leq s\leq\tn
      ,
      \\
      -\ell_{n+1}(s)
      &
      \tfor
      \tn\leq s\leq\tnph
      .
    \end{cases}
  \end{equation}
  Therefore we may write
  \begin{equation}
    \begin{split}
      \int_{\tnmh}^t
      \ellrecop[n+1]
      &
      \passop[n+1]
      \projf n
      +\resrecu[n]-\ellop\slintrecu(s)
      \d s
      \\
      &
      =
      \int_{\tnmh}^t
      \ellrecop[n+1] \passop[n+1]
      \projf n
      +
      \resrecuextmeshchange
      \\
      &
      \phantom{=\int_{\tnmh}}
      -
      \ellrecop[n+1] \passop[n+1]
      \tildefeop[n]a\feta[n]u
      +
      \frac\ellop2
      \qb{\tarecu[n+1]-\tarecu[n-1]}
      \tilde\ell_n(s)
      \:\d s
    \end{split}
  \end{equation}
  Definition \cref{eq:LTSLF:system-form:variable-fespace}
  and $\int_{\tnmh}^t\d s=(t-\tnmh)=\ell_{\nph}(t)\timestep$
  reveal that
  \begin{equation}
    \begin{split}
      \int_{\tnmh}^t
      &
      \ellrecop[n+1] \passop[n+1]
      \projf n
      +\resrecu[n]-\ellop\slintrecu(s)
      \d s
      \\
      &
      =
      \qp{
        \ellrecop[n+1] \fetav[n+]
        -
        \tarecvbase[\nmh]
      }
      \ell_{\nph}(t)
      +
      \frac\ellop2
      \qb{\tarecu[n+1]-\tarecu[n-1]}
      \tilde q_n(t)
\end{split}
\end{equation}
where we introduce
\begin{equation}
  \indexma[qtilden]{\tilde q_n(t)}
    :=
    \int_{\tnmh}^t
    \tilde\ell_n(s)
    \d s
  \end{equation}
  is the unique quadratic that equals $0$ at $\tnmh$, $\tnph$ and satisfies
  \begin{equation}
    \tilde q_n(\tn)
    =
    \frac\timestep{8}.
  \end{equation}
  It can be written the form $\tilde q_n(t)=q_n(t)\timestep$, with $q_n(t)$
  given by \cref{eq:def:quadratic-perturbation-polynomial}.

  To conclude note that for $\tnmh\leq t\leq\tnph$ we have
  \begin{multline}
    \tarecvbase[\nmh]
    +
    \qp{
      \tarecvbase[\nph]
      -
      \tarecvbase[\nmh]
    }
    \ell_{\nph}(t)
    \\
    =
    \tarecvbase[\nph]\ell_{\nph}(t)
    +
    \tarecvbase[\nmh]\ell_{\nmh}(t)
    =
    \lintrecv(t)
    ,
  \end{multline}
  and using the fact that $\ellop\tarecu=\discellop[n]\feta[n]u$ hence
  we obtain
  \begin{equation}
    \quatrecv(t)
    -
    \lintrecv(t)
    =
    \frac{
      \discellop[n+1]\feta[n+1]u
      -
      \discellop[n-1]\feta[n-1]u
    }2
    q_n(t)\timestep,
  \end{equation}
  as claimed.

  Similarly, owing to \crefnosort{eq:def:quadratic-time-reconstructions,%
    eq:def:staggered-linear-time-reconstructions}
  we have
  \begin{equation}
    \begin{split}
      \quatrecu(t)
      -
      \tarecu[n-1]
      &
      =
      \int_{\tnmo}^t\slintrecv(s)\d s
      +(t-\tnmo)\resv
      \\
      &
      =
      \int_{\tnmo}^t
      \lintrecv[n-1]\ell_{n-1}(s)
      +
      \lintrecv[n]\ell_n(s)
      +
      \resv
      \d s
      \\
      &
      =
      \int_{\tnmo}^t
      \frac{\tarecvbase[\nmth]+\tarecvbase[\nmh]}2
      \ell_{n-1}(s)
      +
      \frac{\tarecvbase+\tarecvbase[\nph]}2
      \ell_n(s)
      \\
      &
      \phantom{=\int_{\tnmo}^t}
      -\frac14\qp{\tarecvbase[\nmth]-2\tarecvbase[\nmh]+\tarecvbase[\nph]}
      \\
      &
      \phantom{=\int_{\tnmo}^t}
      +
      \left[ \ellrecop[n-\frac{1}{2}]\fetav[n-] - \tarecvbase[\nmh] \right]
      \\
      &
      \phantom{=\int_{\tnmo}^t}
      +
      \resvextmeshchange
      \d s.
    \end{split}
  \end{equation}
  Using the facts that $\ell_n(t)=\int_{\tnmo}^t\frac{\d s}{\timestep}$
  and {$\ell_{n-1}+\ell_{n}=1$}, and recalling
  \cref{eq:LTSLF:system-form:variable-fespace%
    ,eq:def:elliptic-reconstructor%
    ,eq:def:linear-time-reconstructions}
  yields
  \begin{equation}
    \begin{split}
      \quatrecu(t)
      -
      \tarecu[n-1]
      &
       =
      \int_{\tnmo}^t
      \frac{\tarecvbase[\nmth]-\tarecvbase[\nph]}2\qp{\ell_{n-1}(s)-\frac12}\d s
      \\
      &
      \phantom{=\int_{\tnmo}^t}
      +
      \qp{
        \ellrecop[n]\qb{
          \fetav[n-]
          \timestep
          +
          \passop[n]
          \feta[n-1]u
        }
        -
        \ellrecop[n-1]\feta[n-1]u
      }\ell_n(t)
      \\
      &
      =
      \frac{\tarecvbase[\nmth]-\tarecvbase[\nph]}2q_{\nmh}(t)\timestep
      +
      \lintrecu(t)
      -
      \tarecu[n-1],
    \end{split}
  \end{equation}
  which implies \cref{eq:quadratic-time-reconstruction-residual:recu}
  and concludes the proof.
\end{Proof}
\begin{Lem}[piecewise linear time-reconstruction residual]
  \label{lem:piecewise-linear-time-reconstruction-residual}
  \ \\For each
  $n\integerbetween0N$,
  if $\tnmh\leq t\leq\tnph$ we have
  \begin{equation}
    \label{eq:linear-time-reconstrucion-residual-velocity}
      \slintrecv(t)-\lintrecv(t)
      =
      \frac12
      \qp{
        \secodiff\tarecvbase[\nmh]\ell_{n-1}(t)
        +
        \secodiff\tarecvbase[\nph]\ell_{n+1}(t)
      }
      \qppow{\timestep}2
      ,
  \end{equation}
  and if $\tnmo\leq t\leq\tn$ we have
  \begin{equation}
    \label{eq:linear-time-reconstrucion-residual-recu}
    \slintrecu(t)-\lintrecu(t)
    =
    \frac12
    \qp{
      \secodiff\tarecu[n-1]\ell_{\nmth}(t)
      +
      \secodiff\tarecu[n]\ell_{\nph}(t)
    }
    \qppow\timestep2
    .
  \end{equation}
\end{Lem}
\begin{Proof}
  Suppose that $\tnmh\leq t\leq\tn$ then
  \begin{equation}
    \begin{split}
      \slintrecv(t)-\lintrecv(t)
      &
      =
      \lintrecv[n-1]\ell_{n-1}(t)+\lintrecv[n]\ell_n(t)
      -
      \qp{
        \tarecvbase[\nmh]\ell_{\nmh}(t)
        +
        \tarecvbase[\nph]\ell_{\nph}(t)
      }
      \\
      &
      =
      \frac{\tarecvbase[\nmth]+\tarecvbase[\nmh]}2
      \ell_{n-1}(t)
      +
      \frac{\tarecvbase[\nmh]+\tarecvbase[\nph]}2
      \ell_{n}(t)
      \\
      &
      \phantom=
      -
      \tarecvbase[\nmh]\ell_{\nmh}(t)
      -
      \tarecvbase[\nph]\ell_{\nph}(t)
      \\
      &
      =
      \frac{\tarecvbase[\nmth]}2
      \ell_{n-1}(t)
      +
      \frac{\tarecvbase[\nmh]}2
      +
      \frac{\tarecvbase[\nph]}2
      \ell_{n}(t)
      \\
      &
      \phantom=
      -
      \tarecvbase[\nmh]\ell_{\nmh}(t)
      -
      \tarecvbase[\nph]\ell_{\nph}(t)
      \\
      &
      =
      \frac{\tarecvbase[\nmth]}2
      \ell_{n-1}(t)
      +
      \frac{\tarecvbase[\nmh]}2
      -
      \tarecvbase[\nmh]\ell_{\nmh}(t)
      \\
      &\phantom=
      +
      \frac{\tarecvbase[\nph]}2
      \ell_{n}(t)
      -
      \tarecvbase[\nph]\ell_{\nph}(t)
      \\
      &
      =
      \frac12\Big(
        \tarecvbase[\nmth]
        \ell_{n-1}(t)
        +
        \tarecvbase[\nmh]
        \qp{1-2\ell_{\nmh}(t)}
        \\
        &
        \phantom{=\frac12\Big( }
        +
        \tarecvbase[\nph]
        \qp{\ell_n(t)-2\ell_{\nph}(t)}
      \Big)
    \end{split}
  \end{equation}
  Noting that
  \begin{equation}
    \begin{gathered}
      1-2\ell_{\nmh}(t)=-2\ell_{n-1}(t)
      \\
      \ell_{n}(t)-2\ell_{\nph}(t)
      =
      \ell_{n-1}(t)
    \end{gathered}
  \end{equation}
  and using definition
  \cref{eq:def:centered-second-difference}
  we obtain
  \begin{equation}
    \slintrecv(t)-\lintrecv(t)
    =
    \frac12\qp{\tarecvbase[\nmth]-2\tarecvbase[\nmh]+\tarecvbase[\nph]}\ell_{n-1}(t)
    =
    \frac12
    \secodiff\tarecvbase[\nmh]\ell_{n-1}(t)
    \qppow{\timestep}2
    .
  \end{equation}
  Similarly if $\tn\leq t\leq\tnph$ we get
  \begin{equation}
    \slintrecv(t)-\lintrecv(t)
    =
    \frac12
    \secodiff\tarecvbase[\nph]\ell_{n+1}(t)
    \qppow\timestep2
    .
  \end{equation}
  Therefore
  \begin{equation}
    \slintrecv(t)-\lintrecv(t)
    =
    \begin{cases}
      \frac12\qp{\tarecvbase[\nmth]-2\tarecvbase[\nmh]+\tarecvbase[\nph]}
      \ell_{n-1}(t)
      &
      \text{ for }\tnmh\leq t\leq\tn
      \\
      \frac12\qp{\tarecvbase[\nmh]-2\tarecvbase[\nph]+\tarecvbase[\npth]}
      \ell_{n+1}(t)
      &
      \tfor\tn\leq t\leq\tnph.
    \end{cases}
  \end{equation}
  Owing to the empty common support of $\ell_{n-1}$ and $\ell_{n+1}$
  we sum up to deduce
  \cref{eq:linear-time-reconstrucion-residual-velocity}.

  Showing \cref{eq:linear-time-reconstrucion-residual-recu} is similar,
  for $\tnmh\leq t\leq\tn$ we have
  \begin{equation}
    \begin{split}
      \slintrecu(t)-\lintrecu(t)
      =
      \slintrecu[\nmh]\ell_{\nmh}(t)
      +
      \slintrecu[\nph]\ell_{\nph}(t)
      -
      \lintrecu[n-1]\ell_{n-1}(t)
      -
      \lintrecu[n]\ell_n(t)
      \\
      =
      \frac{\tarecu[n-1]+\tarecu[n]}2
      \qp{1-\ell_{\nph}(t)}
      +
      \frac{\tarecu[n]+\tarecu[n+1]}2
      \ell_{\nph}(t)
      -
      \lintrecu[n-1]\ell_{n-1}(t)
      -
      \lintrecu[n]\ell_n(t)
      \\
      =
      \frac{\tarecu[n-1]+\tarecu[n]}2
      \qp{\ell_{n-1}(t)+\ell_n(t)}
      -
      \lintrecu[n-1]\ell_{n-1}(t)
      -
      \lintrecu[n]\ell_n(t)
      \\
      \phantom=
      +
      \qp{\lintrecu[n+1]-\lintrecu[n-1]}\ell_{\nph}(t)
    \end{split}
  \end{equation}
\end{Proof}
\begin{The}[full time-reconstruction residual]
  \label{the:full-time-recontstruction-residual}
  Using the time-functions $\ell_\nu$ and $q_\nu$ defined in
  \cref{sec:time-discretization}, for each $n\integerbetween1N$,
  we have
  \begin{equation}
    \label{eq:linear-minus-quad-time-reconstrucion-residual-velocity}
    \slintrecv(t)-\quatrecv(t)
    =
    \qp{
      \frac12
      \qp{
        \secodiff\tarecvbase[\nmh]\ell_{n-1}(t)
        +
        \secodiff\tarecvbase[\nph]\ell_{n+1}(t)
      }
      -
      \centdiff\qb{\discellop[n]\feta[n]u}q_n(t)
    }
    \qppow{\timestep}2
  \end{equation}
  if $\tnmh\leq t\leq\tnph$,
  and
  \begin{equation}
    \label{eq:linear-minus-quad-time-reconstrucion-residual-recu}
    \slintrecu(t)-\quatrecu(t)
    =
    \qp{
      \frac12
      \qp{
        \secodiff\tarecu[n-1]\ell_{\nmth}(t)
        +
        \secodiff\tarecu[n]\ell_{\nph}(t)
      }
      -
      \centdiff{\tarecvbase[\nmh]}q_{\nmh}(t)
    }
    \qppow{\timestep}2
  \end{equation}
  if $\tnmo\leq t\leq\tn$.
\end{The}
\begin{Proof}
Subtracting \cref{eq:linear-time-reconstrucion-residual-velocity} from \cref{eq:quadratic-time-reconstruction-residual:V} gives us 
\begin{align*}
\slintrecv(t)-\quatrecv(t) &=\slintrecv(t)-\lintrecv(t)-(\quatrecv(t)-\lintrecv(t)) \\
&= \frac12
\qp{
  \secodiff\tarecvbase[\nmh]\ell_{n-1}(t)
  +
  \secodiff\tarecvbase[\nph]\ell_{n+1}(t)
}
\qppow{\timestep}2 \\
&\phantom{=} - \centdiff\qb{\discellop[n]\feta[n]u}q_n(t)
\qppow{\timestep}2  \\
&= \qp{
  \frac12
  \qp{
    \secodiff\tarecvbase[\nmh]\ell_{n-1}(t)
    +
    \secodiff\tarecvbase[\nph]\ell_{n+1}(t)
  }
  -
  \centdiff\qb{\discellop[n]\feta[n]u}q_n(t)
}
\qppow{\timestep}2.
\end{align*}
Similarly, if we subtracting \cref{eq:linear-time-reconstrucion-residual-recu} from \cref{eq:quadratic-time-reconstruction-residual:recu} gives us 
\begin{align*}
 \slintrecu(t)-\quatrecu(t) 
    &= \slintrecu(t)-\lintrecu(t)-(\quatrecu(t)-\lintrecu(t)) \\
    &=  \frac12
    \qp{
      \secodiff\tarecu[n-1]\ell_{\nmth}(t)
      +
      \secodiff\tarecu[n]\ell_{\nph}(t)
    }
    \qppow\timestep2 \\
    &\phantom{=} - \centdiff{\fetav[n-]}q_{\nmh}(t)
    \qppow{\timestep}2 \\
    &= \qp{
      \frac12
      \qp{
        \secodiff\tarecu[n-1]\ell_{\nmth}(t)
        +
        \secodiff\tarecu[n]\ell_{\nph}(t)
      }
      -
      \centdiff{\tarecvbase[\nmh]}q_{\nmh}(t)
    }
    \qppow{\timestep}2.
\end{align*}
\end{Proof}
\pathword{sec-aposteriori-error-analysis.tex}%
\conword\par%
\section{\Aposteriori error analysis}
\label{sec:aposteriori-error-analysis}
We now present the main analytical result of this paper in the form of
\cref{the:full-error analysis}.  The starting point of the analysis is
given by the error-residual PDE for the error between the
reconstruction of the discrete solution and the exact solution in
\cref{sec:error-residual-PDE}.  We use this PDE to prove
\cref{the:reconstruction-exact-error-residual-estimate}.
In \cref{def:error-indicators} we introduce all the error indicators
needed to state and prove the main result.
\subsection{The reconstruction--exact error--residual PDE}
\label{sec:error-residual-PDE}
The rationale behind the definitions in
\cref{sec:reconstruction} is that differentiation in time
and \cref{eq:def:quadratic-time-reconstructions}
yield
\begin{equation}
  \begin{split}
    \pdt\quatrecu(t)
    -\quatrecv(t)
    &
    =
    \slintrecv(t)
    -
    \quatrecv(t)
    +
    \resvt(t)
    \\
    \pdt\quatrecv(t)
    +
    \ellop \quatrecu(t)
    &=
    \ellop\qb{
      \quatrecu(t)
      -
      \slintrecu(t)
    }
    +
    \resrecut(t)
    +
    \overline{\projf{}}(t)
    ,
  \end{split}
\end{equation}
where $\extprojf$ is the piecewise constant time-extension of the
$\projf n$ over the half-grid:
\begin{equation}
  \indexma[Fbar]{\extprojf}(t)
  =
  \sumifromto n0N\projf n\charfun{\timestaginter}(t),
  \Foreach
  t\in\clinter{\tn[-\fracl12]}{\tnph[N]}
  .
\end{equation}

This allows comparison with the wave equation in system form
\begin{equation}
  \begin{split}
    \pdt u(t)
    - v(t)
    &
    =
    0
    \\
    \pdt v(t)
    +
    \ellop u(t)
    &
    =
    f(t)
  \end{split}
\end{equation}
which, upon interpreting the residuals and referring to
\cref{eq:def:time-extended-residuals}, gives
\begin{equation}
  \label{eq:def:full-residual-vector}
  \begin{split}
    \pdt\qb{\quatrecu-u}
    -\qp{\quatrecv-v}
    &
    =
    \indexma[r0]{\reszero}
    :=
    \slintrecv
    -
    \quatrecv
    +
    \resvt
    \\
    \pdt\qb{\quatrecv-v}
    +
    \ellop\qb{\quatrecu-u}
    &
    =
    \indexma[r1]{\resone}
    :=
    \ellop\qb{\quatrecu-\slintrecu}
    +
    \resrecut
    +
      \extprojf
      -f
  \end{split}
\end{equation}
that is the \indexemph{error-residual partial differential equation}
\begin{equation}
  \label{eq:error-residual-PDE}
\pdt
\discolvectwo{\errrecu}{\errrecv}
+
\dismattwo0{-1}{\ellop}0
\discolvectwo{\errrecu}{\errrecv}
=
\discolvecitwoz r
\end{equation}
with the
\indexemph[reconstruction--exact error]{reconstruction--exact error}
for $\pair uv$
\begin{equation}
  \indexma[rho 0]{\errrecu}
  :=
  \quatrecu-u
  \tand
  \indexma[rho 1]{\errrecv}
  :=
  \quatrecv-v.
\end{equation}
In what follows we respectively denote
the pairs $\pair\errrecu\errrecv$ and $\pair\reszero\resone$ as the
(column) vectors $\errrec$ and $\resvec$.
\begin{The}[reconstruction--exact error--residual estimate]
  \label{the:reconstruction-exact-error-residual-estimate}
  With the notation introduced in
  \cref{sec:error-residual-PDE}
  we have
  \begin{equation}
    \label{eq:main-reconstruction-exact-error-estimate}
  \supergnorm{\errrec}0T
  \leq
  \ergnorm{\errrec(0)}
  +
  2
  \lebergnorm{\resvec}10T
    .
  \end{equation}
\end{The}
\begin{Proof}
  Testing the error--residual PDE \cref{eq:error-residual-PDE}
  with the reconstruction--exact error vector, 
  with $\pdt\errrecu\in\elldom$
  and $\pdt\errrecv\in\ellran$,
  \begin{equation}
    \label{eqn:wave-energy-reconstruction-differential-estimate}
    \begin{split}
      \frac12\ddt
      \qb{
        \ergnorm{
          \errrec
        }
        ^2
      }
      &=
      \frac12\ddt\qb[big]{\duality{\linop A\errrecu}\errrecu+\ltwop\errrecv\errrecv}
      \\
      &=
      \duality{\linop A\pdt\errrecu}\errrecu
      +
      \duality{\pdt\errrecv}\errrecv
      \\
      &=
      \duality[\big]{\linop A\qb{\errrecv
          +
          \reszero
      }}\errrecu
      +
      \duality[\big]{-\linop A\errrecu
        +
        \resone
      }\errrecv
      \\
      &=
      \duality[\big]{\linop A{
          \reszero
      }}\errrecu
      +
      \duality[\big]{
        \resone
      }\errrecv
      =
      \ergprod{
        \resvec
      }{
        \errrec
      }
      \\
      &\leq
      \ergnorm{
        \resvec
      }
      \ergnorm{
        \errrec
      }
    \end{split}
  \end{equation}
  Noting that $\ergnorm{\errrec(t)}$ is piecewise uniformly continuous in $t$
  over $\clinter0T$ the partition $\ltidotsfromto t0N$, there must exist
  a $T^*\in\clinter0T$ such that
  \begin{equation}
    \ergnorm{\errrec(T^*)}=\max_{\clinter0T}\ergnorm{\errrec}
  \end{equation}
  Integrating both sides of
  \cref{eqn:wave-energy-reconstruction-differential-estimate} over the
  time interval $\clinter0{T^*}$ and using the fact that $T\geq T^*$ we
  obtain
  \begin{equation}
    \begin{split}
      \supergnorm{\errrec}0T^2
      &
      :=
      \Normsupon{\ergnorm{\errrec}}{0,T}^2
      =
      \ergnorm{\errrec(T^*)}^2
      \\
      &
      \leq
      \ergnorm{\errrec(0)}^2
      +
      2
      \supergnorm{\errrec}0T
      \lebergnorm{\resvec}10T
    \end{split}
  \end{equation}
  Using the following elementary fact
  \begin{equation}
    {a,b,c\geq0\tand a^2\leq c^2+2ab}\implies{a\leq c+2b}
  \end{equation}
  we conclude that
  \begin{equation}
  \supergnorm{\errrec}0T
  \leq
  \ergnorm{\errrec(0)}
  +
  2
  \lebergnorm{\resvec}10T
    .
  \end{equation}
\end{Proof}
\begin{Def}[error indicators]
  \label{def:error-indicators}
  Let us now introduce the error indicators that appear in the
  \aposteriori error analysis and that we will implement in the
  numerical experiments \cref{sec:numerical-results}:
  \begin{description}
  \item[\indexen{mesh-change indicator}s]
    (nonzero only when the mesh changes)
    \begin{equation}
      \label{eq:def:mesh-change-indicator}
      \begin{split}
        \indexma[mu]{\indmeshchangezero}
        &
        :=
        \qp{
          \potenorm{
            \qb{\passop[n]-\Id}\feta[n-1]u
          }
          +
          \cE\qb[big]{
            \qb{\passop[n]-\Id}\feta[n-1]u
            ,
            \vespace\meet\vespace[n+1]
            ,
            \ellop
          }
        }
        \inverse\timestep
        ,\\
        \indexma[mu 1 n]{\indmeshchangeone}
        &
        :=
        \Big(
          \pivotnorm[big]{\qb{\passop[n+1]-\Id}\fetav[n-]}
          \\
          &
          \phantom{:=\Big(}%
          +
          \cE\qb[big]{\qb{\passop[n+1]-\Id}\fetav[n-],\vespace\meet\vespace[n+1],\pivot}
        \Big)
        \inverse\timestep
        ,
        \\
        \indexma[mu 2]{\indmeshchangelts}
        &
        :=
        \pivotnorm{\qb{\Id-\passop[n+1]}\tildefeop[n]a\feta[n]u}
        +
        \cE[\qb[big]{\Id-\passop[n+1]}\tildefeop[n]a\feta[n]u,\vespace[n+1],\pivot]
        ;
      \end{split}
    \end{equation}
  \item[\indexen{LTS error indicator}s]
    (due to using $\tildefeop[n]a$ in scheme instead of $\discellop$)
    \begin{equation}
    \label{def:LTS-error-indicators}
      \begin{split}
        \indexma[alpha 0]{\alpha_0^n}
        &
        :=
        \pivotnorm{\qb{\discellop-\tildefeop[n]a}\feta[n]u}
        ,\\
        \indexma[alpha 1]{\alpha_1^n}
        &
        :=
        \cE[\tildefeop[n]a\feta[n]u,\vespace[n+1],\pivot]
        ,\\
        \indexma[alpha]{\alpha^n}
        &
        :=
        \alpha_0^n+\alpha_1^n+\indmeshchangelts
        ;
      \end{split}
    \end{equation}
  \item[\indexen{time-error indicator}s]
    (mainly due to time discretization)
    \begin{equation}
      \label{eq:def:time-error-indicator-0}
      \begin{split}
        \indexma[hz0]{\indtimefunzero[n](t)}
        &
        :=
        \timestep^2
        \begin{cases}
          \potenorm{
            \secodiff\tarecv[n-]
            \frac{\ell_{n}(t)-1}2
            -
            \centdiff\qb{\discellop[n-1]\feta[n-1]u}
            q_{n-1}(t)
          },
          &
          t\in\timehalfinter[\nmh]
          ,\\[1mm]
          \potenorm{
            \secodiff\tarecv[n-]
            \frac{\ell_{n}(t)-1}2
            -
            \centdiff\qb{\discellop[n]\feta[n]u}
            q_{n}(t)
          },
          &
          t\in\timehalfinter[n]
          ,
        \end{cases}
        \\
        \indexma[hz1]{\indtimefunone[n](t)}
        &:=
        \timestep^2
        \begin{cases}
          \pivotnorm{
            \ellop\qb{
              \frac12
              \secodiff\tarecu[n]
              \ell_n(t)
              -
              \centdiff{\tarecvbase[\nmh]}q_{\nmh}(t)
            }
          }
          ,
          &
          t\in\timehalfinter[n]
          ,
          \\
          \pivotnorm{
            \ellop\qb{
              \frac12
              \secodiff\tarecu[n]
              \ell_n(t)
              -
              \centdiff{\tarecvbase[\nph]}q_{\nph}(t)
            }
          }
          ,
          &
          t\in\timehalfinter[\nph]
          ;
        \end{cases}
      \end{split}
    \end{equation}
  \item[\indexen{data approximation indicator}]
    (due to a possibly nonzero source)
    \begin{equation}
      \indexma[delta]{\inddatafun(t)}
      :=
      \pivotnorm{
        \projf n
        -
        f(t)
      }
      ;
    \end{equation}
  \item[\indexen{elliptic error indicator}s]
    (the ``standard'' error indicators depending on the residual functional discussed
    in \cref{def:elliptic-estimators} )
    \begin{equation}
    \label{elliptic error indicator}
      \begin{split}
        \indexma[e0]{\indellzero}
        &:=
        \cE[\feta[n]u,\vespace[n],\ellop]
        ,
        \\
        \indexma[e1]{\indellone}
        &:=
        \cE[\feta[\nmh]v,\vespace[n],\pivot]
        ;
      \end{split}
    \end{equation}
  \item[\indexen{time accumulation indicator}s]
    \begin{equation}
      \indtotnospace
      :=
      \int_{t_{\frac{{m}-1}2}}^{t_{\fracl{m}2}}
      \powqpsqrt{
        \qppow{
          \indmeshchangezero
          +
          \indtimefunzero(t)
        }2
        +
        \qppow{\indoperator
          +
          \indmeshchangeone
          +
          \inddatafun(t)
          +
          \indtimefunone(t)
        }2
      }
      \d t
    \end{equation}
    for $n=\ceil{2{m}}$ and ${m}\integerbetween1{2N}$.
  \end{description}
\end{Def}
\begin{The}[full-error analysis]
  \label{the:full-error analysis}
  With the notation introduced in \cref{def:error-indicators}
  we have the following error estimates
  \begin{equation}
  \label{estimateU}
    \max_{0\leq n\leq N}
    \potenorm{%
      \feta[n]u
      -
      u\attime n
    }
    \leq
    \max_{1\leq n\leq N}\indellzero
    +
    \ergnorm{\errvec(0)}
    +
    2
    \sumifromto m1{2N}\indtotnospace
    ,
  \end{equation}  
  and
  \begin{equation}
  \label{estimateV}
    \max_{1\leq n\leq N}
    \pivotnorm{%
      \fetav
      -
      v\attime\nmh
    }
    \leq
    \max_{1\leq n\leq N}\indellone
    +
    \ergnorm{\errvec(0)}
    +
    2
    \sumifromto m1{2N}\indtotnospace
    .
  \end{equation}
\end{The}
\begin{Proof}
  Using the facts that $\quatrecu[n]=\tarecu[n]$ and
  $\quatrecv[n-]=\tarecv[n-]$ for $n\integerbetween0N$, we can
  decompose the full discretization errors as follows
  \begin{equation}
    \begin{split}
      \indexma[e0]{\erroru[n]}
      &
      :=
      \feta[n]u
      -
      \tarecu[n]
      +
      \quatrecu[n]
      -
      u\attime n
      =:
      {\errellu\attime n}
      +
      \errrecu\attime n%
      \\
      \indexma[e1]{\errorv[n-]}
      &
      =
      \fetav
      -
      \tarecv[n-]
      +
      \quatrecv[n-]
      -
      v\attime\nmh
      =:
      \indexma[epsilon 1 n 12]{\errellv\attime{\nmh}}
      +
      \errrecv\attime{\nmh}%
      ,
    \end{split}
  \end{equation}
  where this defines the staggered components of the
  \indexemph{full error} $\errvec$ and its splitting into
  elliptic part $\errell$ and time-dependent part $\errrec$.

Thanks to the \aposteriori error estimators
discussed in \cref{def:elliptic-estimators} and the
equivalence between $\elldom$'s norm and the potential energy norm we
have
\begin{align}
  \potenorm{\errellu\attime n}
  &
  =
  \potenorm{
    \feta[n]u
    -
    \tarecu[n]
  }
  \leq
  \indellzero
  \intertext{and }
  \pivotnorm[big]{\errellv\attime{\nmh}}
  &
  =
  \pivotnorm{
    \feta[\nmh]v
    -
    \tarecv[n-]
  }
  \leq
  \indellone
  .
\end{align}
From \cref{eq:main-reconstruction-exact-error-estimate}
we also have
\begin{multline}
  \max_{1\leq n\leq N}%
  \maxi{\potenorm{\erc[0]\attime n}}{\pivotnorm{\erc[1]\attime{\nmh}}}
  \\
  \leq
  \supergnorm{\errrec}0T
  \leq
  \ergnorm{\errrec(0)}
  +
  2
  \lebergnorm{\resvec}10T
  .
\end{multline}
With definition \cref{eq:def:full-residual-vector} in mind we may write
\begin{equation}
    \ergnorm{\resvec}^2
    =
    \potenorm\reszero^2
    +
    \pivotnorm\resone^2
\end{equation}
and proceed to bound both terms separately.

Owing to
\cref{eq:linear-minus-quad-time-reconstrucion-residual-velocity}
and
\cref{eq:def:time-extended-residuals}
we see that when $n\integerbetween0N$
and $t\in\timeinter$
\begin{equation}
  \begin{split}
    \reszero(t)
    &
    =
    \slintrecv(t)
    -
    \quatrecv(t)
    +
    \resvt
    =
    \\
    &
    =
    \resvexttime
    +
    \resvextmeshchange
    \\
    &\phantom=
    +
    \timestep^2
    \begin{cases}
      \qp{
        \qp{
          \frac12%
          \ell_n(t)
          -
          \frac14}
        \secodiff\tarecvbase[\nmh]
        -
        \centdiff\qb{\discellop[n-1]\feta[n-1]u}
        q_{n-1}(t)
      }
      &
      \text{ for }
      t\leq\tnmh
      \\
      \qp{
        \qp{
          \frac12%
          \ell_{n-1}(t)
          -
          \frac14}
        \secodiff\tarecvbase[\nmh]
        -
        \centdiff\qb{\discellop[n]\feta[n]u}
        q_{n}(t)
      }
      &\text{ for }
      \tnmh<t
    \end{cases}
    \\
    &=
    \resvextmeshchange
    \\
    &\phantom=
    +
    \timestep^2
    \secodiff\tarecv[n-]
    \frac{\ell_{n}(t)-1}2
    -
    \timestep^2
    \begin{cases}
      \centdiff\qb{\discellop[n-1]\feta[n-1]u}
      q_{n-1}(t)
      &
      \text{ for }
      t\leq\tnmh
      \\
      \centdiff\qb{\discellop[n]\feta[n]u}
      q_{n}(t)
      &\text{ for }
      \tnmh<t    
    \end{cases}
  \end{split}
\end{equation}
By definitions \cref{eq:def:potential-energy-norm},
\cref{eq:def:mesh-change-indicator}
and Lemma \ref{lem:reconstrution-and-two-spaces}
we have the following bound 
\begin{equation}
  \potenorm{\resvextmeshchange}
  \\
  \leq
  {\indmeshchangezero}
  .
\end{equation}
Recalling \cref{eq:def:time-error-indicator-0} we obtain the
following bound, for all
$t\in\timeinter$ with $n=\ceil t$,
\begin{equation}
  \potenorm{\reszero(t)}
  \leq
  \indmeshchangezero
  +
  \indtimefunzero[n](t)
  .
\end{equation}

Next, we bound the residual $\resone$ which, thanks to 
\cref{eq:linear-minus-quad-time-reconstrucion-residual-recu}
and
\cref{eq:def:alt:displacement-residual}
can be written as
\begin{equation}
  \label{eq:residual-one}
  \begin{split}
    \resone(t)
    &=
    \resrecut
    +
    {\overline{\projf{}}-f}
    +
    \ellop\qb{\quatrecu-\slintrecu}
    \\
    &
    =
    \resrecuextoperators
    +
    \resrecuextmeshchange
    \\
    &\phantom{=\Bigg(}
    +
    \projf n
    -
    f
    +
    \frac14\qp{
      \discellop[n+1]\feta[n+1]u
      -
      2\discellop[n]\feta[n]u
      +
      \discellop[n-1]\feta[n-1]u
    }
    \\
    &\phantom{=\Bigg(}
    -
    \ellop
    \timestep^2
    \begin{cases}
      \frac12
      \secodiff\tarecu[n]\ell_{\nph}(t)
      -
      \centdiff{\tarecvbase[\nmh]}q_{\nmh}(t)
      &\text{ for }t\in\timehalfinter[n]
      \\
      \frac12
      \secodiff\tarecu[n]\ell_{\nmh}(t)
      -
      \centdiff{\tarecvbase[\nph]}q_{\nph}(t)
      &\text{ for }t\in\timehalfinter[\nph]
    \end{cases}
   \end{split}
\end{equation}
for all $n\integerbetween0N$ and $t\in\timestaginter$.

The first term on the right-hand side of \cref{eq:residual-one}
can be decomposed as follows
\begin{equation}
  \begin{split}
    \pivotnorm{
      \resrecuextoperators
    }
    &
    \leq
    \pivotnorm{\qb{\discellop-\tildefeop[n]a}
      \feta[n]u}
    \\
    &
    \phantom\leq
    +
    \pivotnorm{\qb{\Id-\ellrecop[n+1]}
      \tildefeop[n]a\feta[n]u}
    \\
    &
    \phantom\leq
    +
    \pivotnorm{\ellrecop[n+1]\qb{\Id-\passop[n+1]}
      \tildefeop[n]a\feta[n]u}
    \\
    &
    \leq
    \alpha_0^n
    +
    \alpha_1^n
    +
    \indmeshchangelts
    =
    \indoperator
    .
  \end{split}
\end{equation}
Here we have used \cref{def:error-indicators}
and Lemmas \ref{lem:residual-error-aposteriori}
and \ref{lem:reconstrution-and-two-spaces}.

To bound the second term in \cref{eq:residual-one} we use Lemma
\ref{lem:reconstrution-and-two-spaces} and definition
\cref{eq:def:mesh-change-indicator} to obtain
\begin{equation}
  \begin{split}
    \inverse\timestep
    &
    \pivotnorm{\ellrecop[n+1]\passop[n+1]\fetav[n-]-\tarecv[n-]}
    \\
    &=
    \inverse\timestep
    \pivotnorm{\qb{\ellrecop[n+1]\passop[n+1]-\ellrecop[n]}\fetav[n-]}
    \leq
    \indmeshchangeone
    .
  \end{split}
\end{equation}
\textnote{Add something about how to get the $\indtimefunone(t)$ term.}
Definitions in \cref{def:error-indicators} lead to the following bound
\begin{equation}
  \pivotnorm{\resone(t)}
  \leq
  \indoperator
  +
  \indmeshchangeone
  +
  \inddatafun(t)
  +
  \indtimefunone(t).
\end{equation}
Summing up we have
\begin{equation}
  \begin{split}
    \int_0^T
    \ergnorm{\resvec(t)}
    \d t
    =
    \sumifromto{m}1{2N}
    \indtotnospace
  \end{split}
\end{equation}
where $\indtotnospace$ is defined in \cref{def:error-indicators}.

Noting that with the discrete initial data taken as the
Ritz/$\leb2$ projections of $u(0)$ and $v(0)$,
\begin{equation}
  \ergnorm{\errrec(0)}
  \leq
  \ergnorm{\errvec(0)}
\end{equation}
we have thus
\begin{equation}
  \begin{split}
    \max_{0\leq n\leq N}
      \potenorm{\erroru[n]}
    \leq
    \max_{1\leq n\leq N}\qp{
      \indellzero
      +
      \maxi{\potenorm{\erc[0]\attime n}}{\pivotnorm{\erc[1]\attime{\nmh}}}
    }
    \\
    \leq
    \max_{1\leq n\leq N}\indellzero
    +
    \ergnorm{\errvec(0)}
    +
    2
    \sumifromto m1{2N}\indtotnospace
    .
  \end{split}
\end{equation}
Similarly
\begin{equation}
  \max_{1\leq n\leq N}
  \pivotnorm{\errorv[n-]}
  \leq
  \max_{1\leq n\leq N}\indellone
  +
  \ergnorm{\errvec(0)}
  +
  2
  \sumifromto m1{2N}\indtotnospace
  .
\end{equation}
\end{Proof}
\pathword{sec-numerical-results.tex}
\conword\par%
\section{Numerical results}
\label{sec:numerical-results}

We now provide a numerical example involving a \emph{time-varying mesh}
and the \indexemph{Gaussian beam} as solution for the exact problem.
\subsection{Set-up}
Consider the one-dimensional wave equation
\cref{eq:def:wave-equation:concrete} in $ \Omega = (-10, 10)$ with
homogeneous Dirichlet boundary conditions, i.e. $ \Gamma = \Gamma_D, c
\equiv 1$, and zero source, $f(x,t) = 0$. The exact solution is a
right-moving Gaussian pulse centered about $x = 1$ and $t = 0$:
\begin{equation}
  u(x,t) = \expp{-4(x-1-t)^2}.
\end{equation}
\par%
\changes{%
  For the numerical solution, we use piecewise linear
  $\sobh1$-conforming finite elements on a nonuniform mesh with
  mass-lumping in space and the leapfrog-based local time-stepping (LF-LTS)
  method with global time-step $\timestep$ without
  stabilization \citep[see][for details]{GroteMichelSauter:21:article:Stabilized}.
  The estimator functionals $\cE$
  are realised the residual Babuška--Rheinboldt on compatible meshes,
  discussed in detail in \cref{sec:residual-esimators}.
}%
  
  At any discrete time $n$ the mesh $\mesh[n]m$, which partitions the
  domain $\W$, is subdivided into a coarse part $\coarsemesh[n]m$ of
  mesh-size $\coarsemeshsize=h$ and a fine part $\finemesh[n]m$ of
  mesh-size $\finemeshsize=\coarsemeshsize/2$ (note that
  $h=\coarsemeshsize$ and $\finemeshsize$ themselves does not depend
  on time).  The initial coarse mesh $\mesh[0]m$ covers the subset
  $\Omega^\coarse_0=[-10, -1.9]\cup [3.9,10]$, while the initial fine
  mesh covers the interval $\Omega^\fine_0 = [-1.9, 3.9]$, inside each
  of which we use an equidistant mesh with respective mesh-sizes
  $\coarsemeshsize$ or $\finemeshsize$.  Hence inside
  $\Omega^\fine_n$, the LF-LTS method takes two local time-steps of
  size $\timestep/2$ for each global time-step of size $\timestep$
  inside $\Omega^\coarse_n$.

\begin{figure}
  \begin{subfigure}[t]{.49\textwidth}
    \centering
    \includegraphics[width=.8\linewidth]{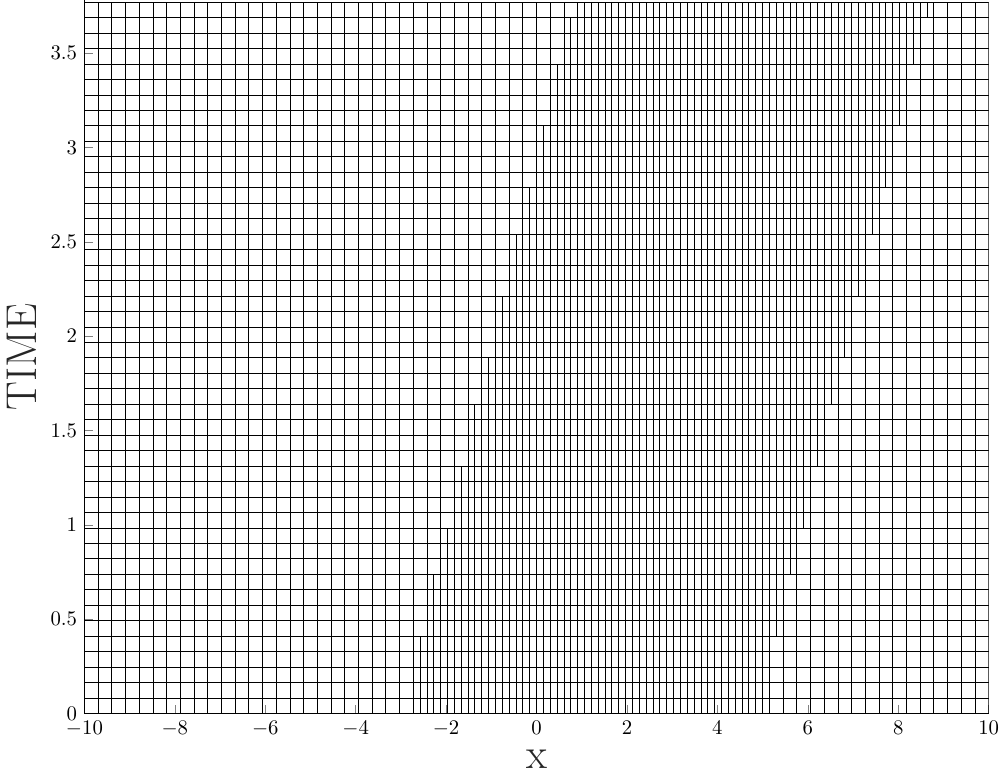}  
    \caption{Time-evolving mesh.}
    \label{fig:mesh}
  \end{subfigure}
  \hfill
  \begin{subfigure}[t]{.49\textwidth}
    \centering
    \includegraphics[width=.8\linewidth]{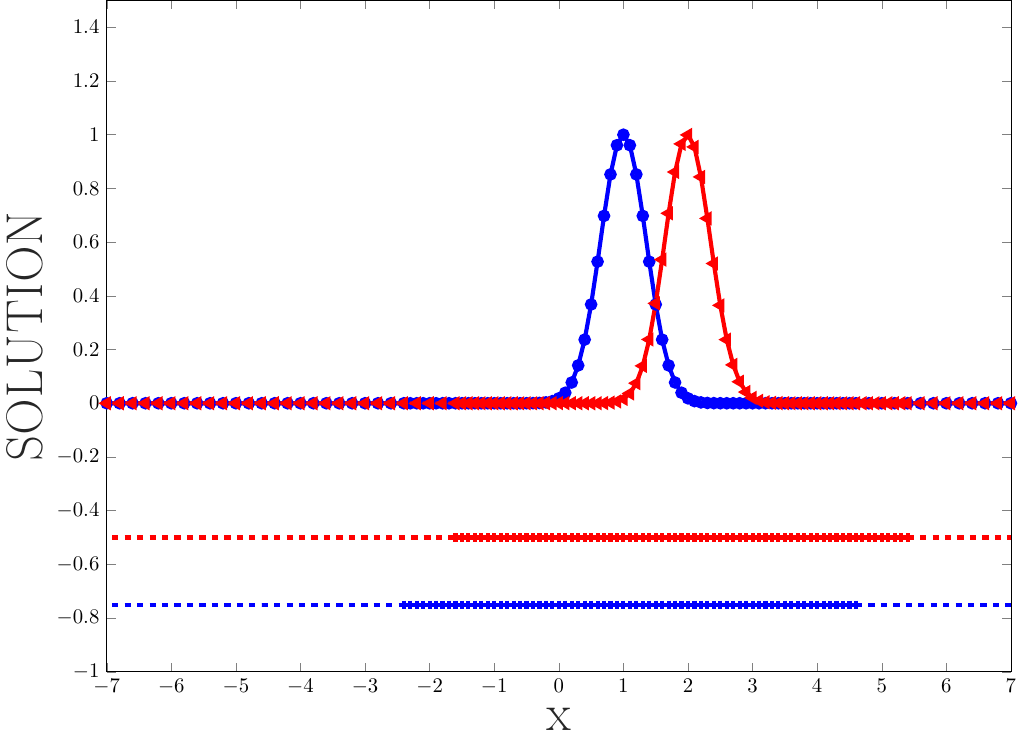}  
    \caption{Numerical solution and refined mesh at time $t=0$ (blue) and $t = 1$ (red).}
    \label{fig:solution}
  \end{subfigure}
\end{figure}

  The fine part, $\finemesh m$, of the mesh $\mesh[n]m$, which has all
  elements length $\finemeshsize$, ``follows'' the peak of Gaussian
  pulse as this propagates rightward across $\Omega$.  The mesh (and
  hence the associated FE space $\vespace$) changes whenever the
  elapsed time from the previous mesh change is greater then the
  coarse mesh-size $\coarsemeshsize$. Hence the fine mesh
  $\finemesh{m}$ moves to the right, as $n$ grows, with the same unit
  wave speed as the pulse, while two subsequent meshes $\vespace$ and
  $\vespace[n+1]$ always remain compatible (see \cref{compatible
    meshes}) during any mesh change.  The resulting space-time mesh
  is plotted in \cref{fig:mesh}. On newly created elements by
  refinement, the FE solution is interpolated on the finer mesh; hence
  no additional discretization error occurs.  Inside coarse elements produced
  by merging two fine elements, however, the removal of the node common
  to those to fine elements introduces an additional discretization error.

  Finally we take the global time-step to be $\timestep:=0.52 h$,
  to ensure it lies just under the CFL stability limit of
  a uniform mesh with mesh-size $h$ (which equals $\coarsemeshsize$ for our
  nonuniform meshes).

\subsection{Discussion}
\begin{figure}
\begin{subfigure}[t]{.475\textwidth}
  \centering
  \includegraphics[width=.8\linewidth]{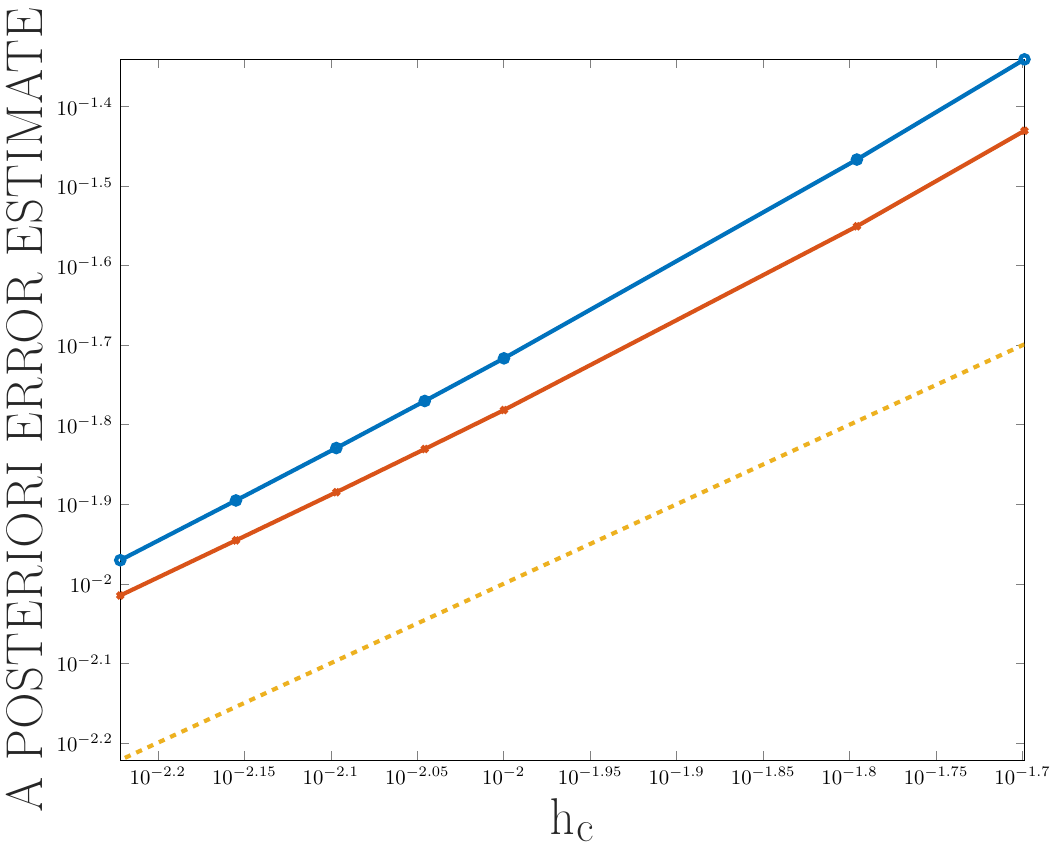}  
  \caption{%
    Convergence rate of the \aposteriori error estimate
    \cref{estimateU} (blue), \cref{estimateV} (red) and
    $\Oh(h)$ (yellow dash-dot).}
  \label{fig:aposteriori}
\end{subfigure}
\hfill
\begin{subfigure}[t]{.475\textwidth}
\centering
  \includegraphics[width=.8\linewidth]{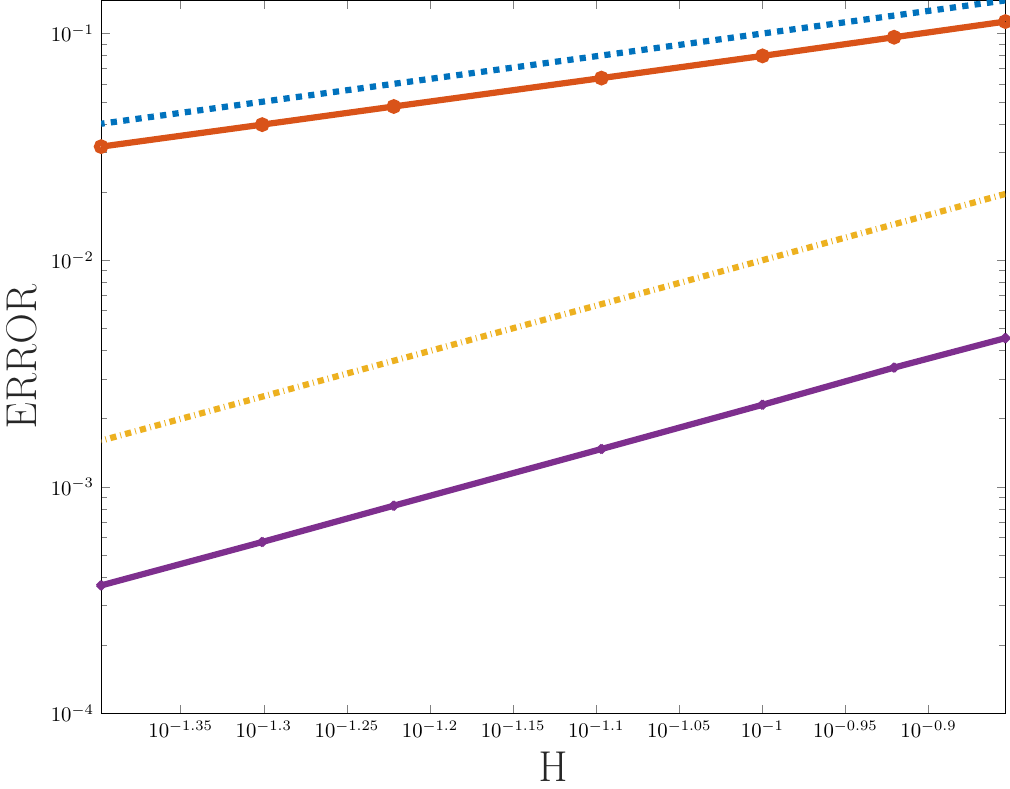}  
  \caption{LF-LTS-FEM convergence on a time evolving mesh. Relative
    energy-norm error (solid red) and $\leb2(\W)$-norm error (solid purple) and
    rates $\Oh(h)$ (blue dash-dot) and $\Oh(h^2)$ (yellow dash-dot).}
  \label{fig:LF-LTS-convergence}
\end{subfigure}
\newline
\begin{subfigure}[t]{.475\textwidth}
  \centering
  \includegraphics[width=.775\linewidth]{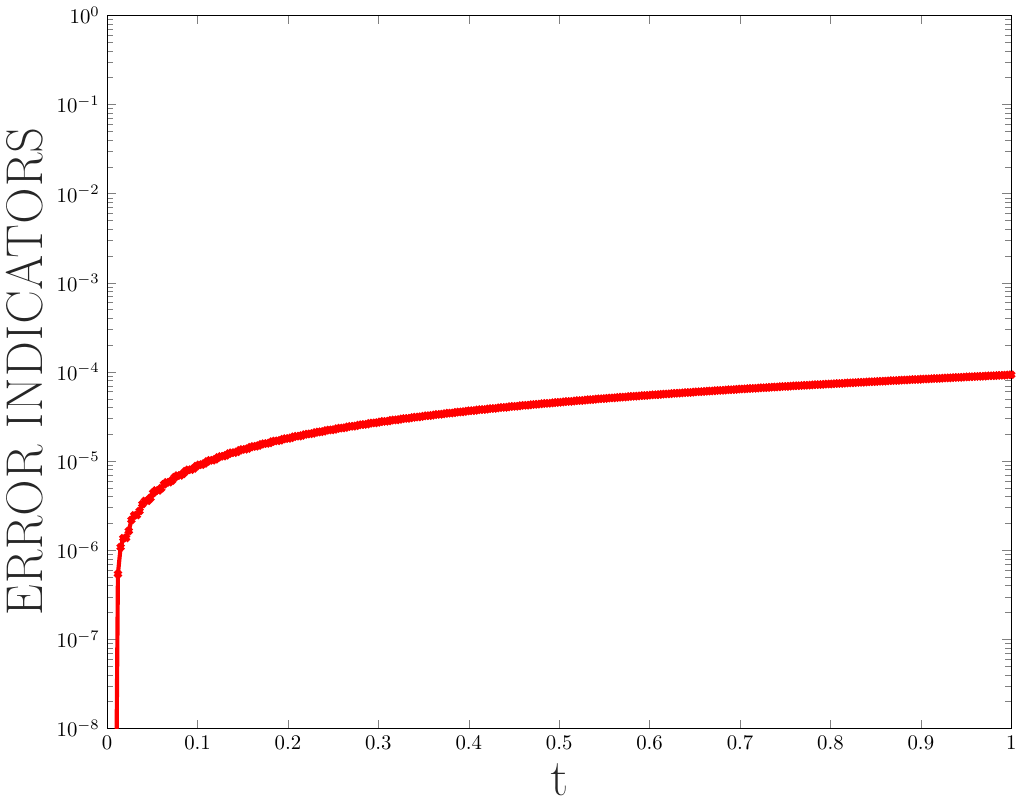}  
 \caption{Time evolution of the error indicator $\vartheta_0^n$ in \cref{eq:def:time-error-indicator-0}.}
  \label{fig:r0-terms_accu}
\end{subfigure}
\hfill
\begin{subfigure}[t]{.475\textwidth}
  \centering
  \includegraphics[width=.775\linewidth]{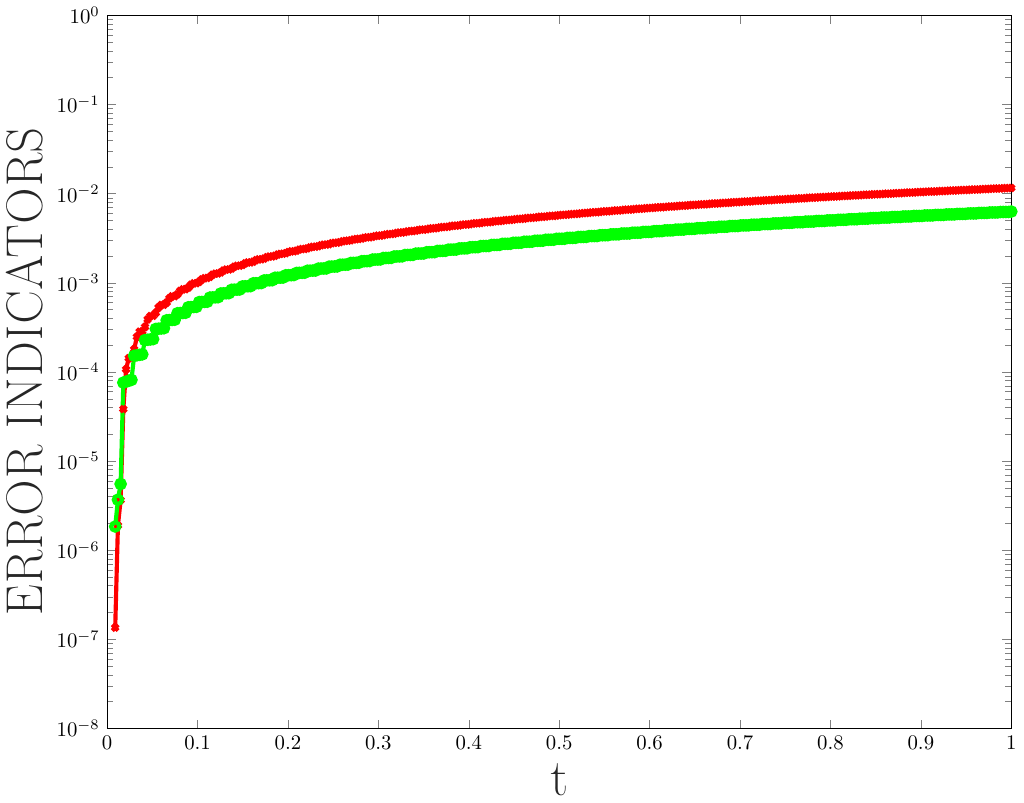}  
    \caption{Time evolution of the time error indicator $\vartheta_1^n$ (red) in \cref{eq:def:time-error-indicator-0} and the LTS error indicator $\alpha^n$ (green) in \cref{def:LTS-error-indicators}.}
    \label{fig:r1-terms_accu}
\end{subfigure}
\newline
\begin{subfigure}[t]{.475\textwidth}
  \centering
  \includegraphics[width=.775\linewidth]{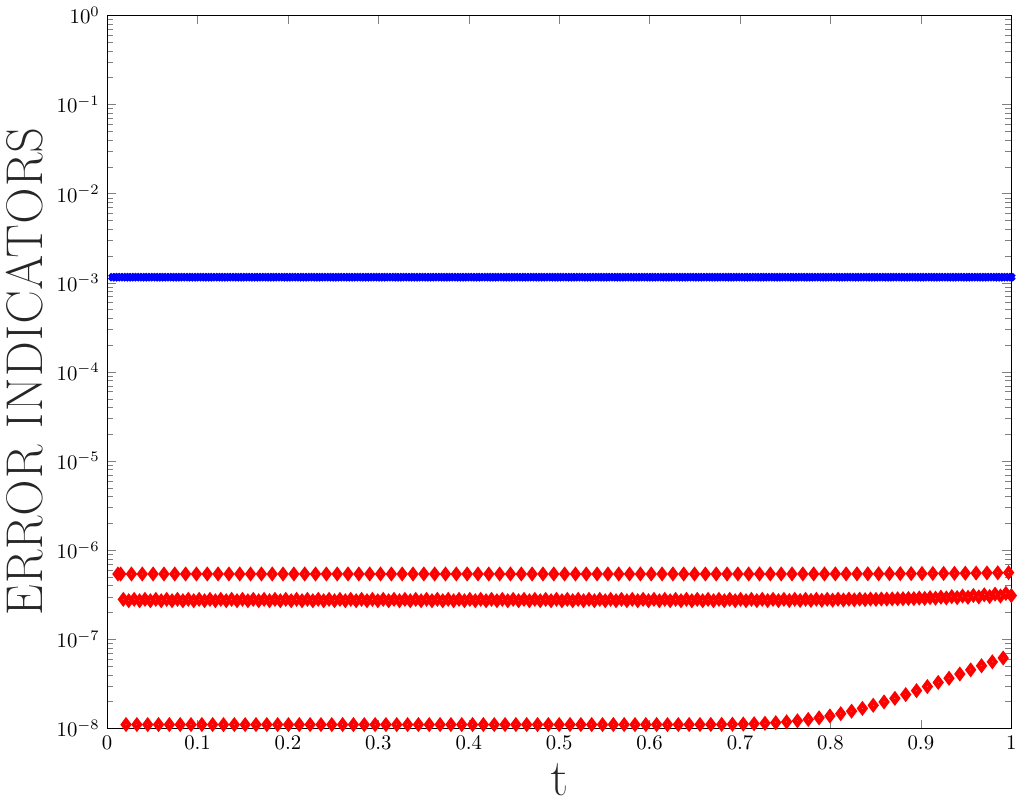}  
 \caption{Elliptic error indicator $\varepsilon_0^n$ in \cref{elliptic error indicator} (blue) and time error indicator $\vartheta_0^n$ in \cref{eq:def:time-error-indicator-0} (red) vs. time without time accumulation.}
  \label{fig:r0-terms}
\end{subfigure}
\hfill
\begin{subfigure}[t]{.475\textwidth}
  \centering
  \includegraphics[width=.775\linewidth]{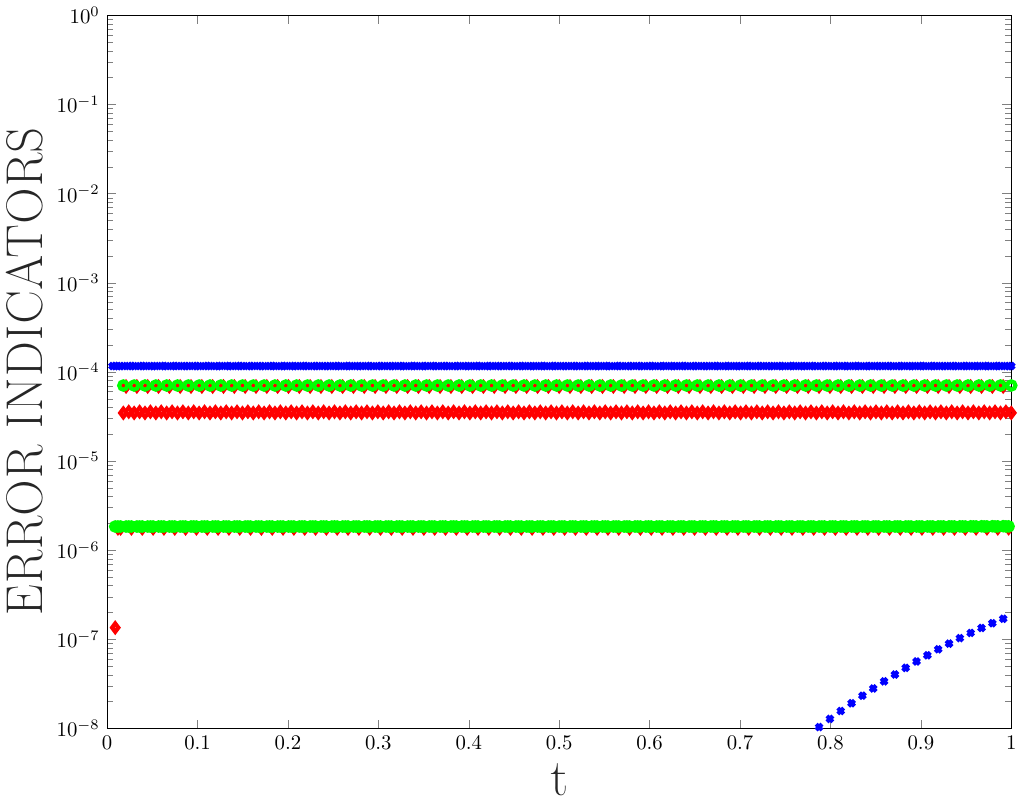}  
    \caption{Elliptic error indicator  $\varepsilon_1^n$ in \cref{elliptic error indicator} (blue), time error indicator $\vartheta_1^n$ (red) in \cref{eq:def:time-error-indicator-0}, and LTS error indicator $\alpha^n$ (green) in \cref{def:LTS-error-indicators} vs. time without time accumulation.}
    \label{fig:r1-terms}
\end{subfigure}
\end{figure}

In \cref{fig:solution}, we display the numerical solutions and the
underlying meshes for $h=0.3$ at initial time $0$ and when
time is $1$. The entire space-time time-evolving mesh
with $\coarsemeshsize = 0.3$ is shown in \cref{fig:mesh}. The refined
part moves to the right with the same unit speed as the Gaussian
pulse. \Cref{fig:LF-LTS-convergence} confirms that the numerical
method \cref{eq:LTSLF:two-step-form:variable-fespace}, including
local time-stepping and a time-evolving mesh, achieves the optimal
convergence rates $\Oh(h)$ and $\Oh(h^2)$ with
respect to the $\sobh1(\Omega)$- and $\leb2(\Omega)$-norm, respectively.

In \cref{fig:aposteriori} the convergence rates of the \aposteriori
error estimates introduced in \cref{the:full-error analysis} are
displayed. We observe that estimate \cref{estimateV} is slightly
smaller then estimate \cref{estimateU}, but both achieve a convergence
rate of $\Oh(h)$. In \cref{fig:r0-terms_accu} and
\cref{fig:r1-terms_accu} the individual indicators in
\cref{def:error-indicators} accumulated over time are displayed. The
behavior of the LTS error indicator $\alpha^n$ in
\cref{def:LTS-error-indicators} and time-error indicators
$\indtimefunzero[n](t)$ and $\indtimefunone[n](t)$ together with the
elliptic error indicators $\indellzero$ and $\indellone$ in
\cref{elliptic error indicator} are shown in \cref{fig:r0-terms} and
\cref{fig:r1-terms} vs. time without accumulation. Note that the
elliptic error indicators $\varepsilon_0^n$ and $\varepsilon_1^n$ in
\cref{elliptic error indicator} are equal to zero whenever no mesh
change occurs. The mesh-change indicators $\indmeshchangezero$ and $
\indmeshchangeone$ \cref{eq:def:mesh-change-indicator} are not
displayed here, as mesh coarsening/refinement occurs only in regions
where the solution is nearly zero. Since the source $f$ is identically
zero, the data approximation indicator $\inddatafun(t)$ also remains
identically zero in this example.

\section{Conclusion}
\label{sec:conclusion}
Building on the time-discrete analysis
\citet{GeorgoulisLakkisMakridakisVirtanen:16:article:A-Posteriori} we
have derived rigorous a posteriori error bounds for a fully discrete
Galerkin formulation of the wave equation with explicit leapfrog time
integration and mesh change. Moreover, our error bounds also accommodate the use
of leapfrog based local time-stepping methods
\citet{%
  DiazGrote:09:article:Energy-Conserving,%
  GroteMitkova:10:article:Explicit-Local,%
  GroteMehlinSauter:18:article:Convergence,%
  GroteMichelSauter:21:article:Stabilized} %
which overcome the stringent CFL stability condition imposed on
explicit time integrators by local mesh refinement.

The fully discrete a posteriori error bounds for the displacement in
the energy norm and for the velocity in the $\leb2$-norm are given in
\cref{the:full-error analysis}.  All the error indicators in the two
upper error bounds \eqref{estimateU}, \eqref{estimateV} are fully
computable while our numerical results in \cref{sec:numerical-results}
confirm their expected optimal rates of convergence with mesh
refinement.

By monitoring local contributions from the error indicators,
algorithms for automatic space-time mesh adaptation can be devised for
computational wave propagation without sacrificing the explicitness of
time integration.  Thus our fully discrete \aposteriori error
estimates for the wave equation pave the way for incorporating
adaptivity with mesh change into explicit time integration while
retaining its ease of use, efficiency and inherent parallelism.
\clearpage
\appendix
\longonly{
  \section{Residual estimators} We discuss now an application
  with one possible choice for the elliptic error estimators $\cE$
  introduced in \S\ref{def:elliptic-estimators}.
}%
\label{sec:residual-esimators}
\subsection{Compatible meshes}
\label{compatible meshes}
\longonly{
  For what concerns mesh refinement, coarsening and management, we
  follow the ideas described in \citet{SchmidtSiebert:05:book:Design},
  where the original simplicial mesh subdivision algorithms of
  \citet{Mitchell:89:article:A-Comparison-of-Adaptive,Kossaczky:94:article:A-recursive}
  are discussed and to which we refer the reader for the details.
  
  We assume that the domain $\Omega$ is a polytope and that it can be
  partitioned into simplices exactly with the coarsest mesh, $\mesh m$
  called the \indexemph{macro triangulation} where every element of
  $\mesh m$ is ``ready'' to be bisected (following the newest vertex
  bisection algorithm in $2=d$ and the Kossaczký algorithm in $3=d$),
  then we are provided with a finite forest of infinite complete (or
  perfect) binary trees, $\frak T=\setofsuch{\frak t^M}{M\in\mesh
    m}$. For each $M\in\mesh{m}$ each node in $\frak t^M$ represents a
  subsimplex of $M$ and its two children represent the subsimplices at
  the next refinement level (see
  \cref{fig:triangle-refinement-by-bisection} for an example). A
  mesh $\mesh k$ that is obtained via refinement by bisection of
  $\mesh m$ is represented by a forest of full \emph{finite} binary
  trees, where each tree is that is a full finite subtree of one of the
  $\frak t^M$, where each leaf (i.e., a tree node that has no children)
  corresponds to an element of $\mesh k$.
\begin{figure}
  \begin{center}
   \begin{tikzpicture}
     [level distance=24pt,
       every node/.style={draw,minimum size=20pt,circle,inner sep=1pt},
       level 1/.style={sibling distance=128pt},
       level 2/.style={sibling distance=64pt},
       level 3/.style={sibling distance=32pt}]
     \node {$K_1$}
     child {node {\(K_2\)}
       child {node {\(K_4\)}
         child {node {\(K_8\)}}
         child {node {\(K_9\)}}
       }
       child {node {\(K_5\)}
         child {node {\(K_{10}\)}}
         child {node {\(K_{11}\)}}
       }
     }
     child {node {\(K_3\)}
       child {node {\(K_6\)}
         child {node {\(K_{12}\)}}
         child {node {\(K_{13}\)}}
       }
       child {node {\(K_7\)}
         child {node {\(K_{14}\)}}
         child {node {\(K_{15}\)}}
       }
     }
     ;
   \end{tikzpicture}
   \begin{tikzpicture}
     \coordinate (A1) at (0,0);
     \coordinate (A1B1) at (.5,2);
     \coordinate (B1) at ($(A1)+(A1B1)$);
     \coordinate (B1C1) at (1.5,-1);
     \coordinate (C1) at ($(B1)+(B1C1)$);
     \coordinate (Z1) at ($(A1)!.5!(B1)$);
     \coordinate (C1Z1) at ($(Z1)-(C1)$);
     \coordinate (Z2) at ($(C1)!.5!(A1)$);
     \coordinate (Z1Z2) at ($(Z2)-(Z1)$);
     \coordinate (Z3) at ($(A1B1)!.5!(C1)$);
     \coordinate (Z1Z3) at ($(Z3)-(Z1)$);
     \coordinate (Z4) at ($(A1)!.5!(Z1)$);
     \coordinate (Z5) at ($(C1)!.5!(Z1)$);
     \coordinate (Z6) at (Z5);
     \coordinate (Z7) at ($(B1)!.5!(Z1)$);
     \coordinate (Z2Z4) at ($(Z4)-(Z2)$);
     \coordinate (Z2Z5) at ($(Z5)-(Z2)$);
     \coordinate (Z3Z6) at ($(Z6)-(Z3)$);
     \coordinate (Z3Z7) at ($(Z7)-(Z3)$);
     \node[left] at (A1) {\tiny ${\vec a}^1$};
     \node[left] at (A1B1) {\tiny ${\vec b}^1$};
     \node[left] at (Z1) {\tiny ${\vec z}^1$};
     \node[right] at (C1) {\tiny ${\vec c}^1$};
     \fill (C1) circle(1pt);
     \fill (Z1) circle(1pt);
     \foreach \oltikzi in {0,3,6,9}{
       \draw (\oltikzi,0)
       --++(A1) coordinate (A1\oltikzi)
       --++(A1B1) coordinate (B1\oltikzi)
       --++(B1C1) coordinate (C1\oltikzi)
       --cycle
       ;
     }
     \foreach \oltikzi in {3,6,9}{
       \draw (C1\oltikzi)--++(C1Z1) coordinate (Z1\oltikzi);
     }
     \fill[color=a] (Z13) circle (1pt);
     \node[above,rotate=76] at (Z13) {\color a\tiny ${\vec z}^1={\vec c}^2={\vec c}^3$};
     \node[left] at (A13) {\tiny ${\vec a}^2$};
     \fill[color=b] (C13) circle (1pt);
     \node[rotate=76,below] at (C13) {\color b\tiny ${\vec b}^2={\vec a}^3$};
     \node[left] at (B13) {\tiny ${\vec b}^3$};
     \foreach \oltikzi in {6,9}{
       \foreach \oltikzj in {2,3}
       \draw (Z1\oltikzi)--++(Z1Z\oltikzj) coordinate (Z\oltikzj\oltikzi);
     }
     \node[left] at (A16) {\tiny ${\vec a}^4$};
     \node[left] at (B16) {\tiny ${\vec b}^7$};
     \fill[color=a] (Z26) circle(1pt);
     \node[rotate=27,below] at (Z26) {\color a\tiny ${\vec z}^2={\vec c}^4={\vec c}^5$};
     \fill[color=a] (Z36) circle(1pt);
     \node[rotate=327,above] at (Z36) {\color a\tiny ${\vec z}^3={\vec c}^6={\vec c}^7$};
     \foreach \oltikzj in {4,5}
     \draw (Z29)--++(Z2Z\oltikzj) coordinate (Z\oltikzj9);
     \foreach \oltikzj in {6,7}
     \draw (Z39)--++(Z3Z\oltikzj) coordinate (Z\oltikzj9);
     \node[left] at (A19) {\tiny ${\vec a}^8$};
     \node[left] at (B19) {\tiny ${\vec b}^{15}$};
     \coordinate (K1) at (barycentric cs:A10=1,B10=1,C10=1);
     \coordinate (K2) at (barycentric cs:A13=1,C13=1,Z13=1);
     \coordinate (K3) at (barycentric cs:B13=1,C13=1,Z13=1);
     \coordinate (K4) at (barycentric cs:A16=1,Z26=1,Z16=1);
     \coordinate (K5) at (barycentric cs:C16=1,Z26=1.2,Z16=1);
     \coordinate (K6) at (barycentric cs:C16=1,Z16=1.2,Z36=1);
     \coordinate (K7) at (barycentric cs:B16=1,Z16=1,Z36=1);
     \coordinate (K8) at (barycentric cs:A19=0.8,Z29=1,Z49=1.2);
     \coordinate (K9) at (barycentric cs:Z19=0.8,Z29=1,Z49=1.2);
     \coordinate (K10) at (barycentric cs:Z19=1.05,Z29=0.9,Z59=1.05);
     \coordinate (K11) at (barycentric cs:C19=0.9,Z29=0.8,Z59=1.3);
     \coordinate (K12) at (barycentric cs:C19=0.9,Z39=0.8,Z59=1.3);
     \coordinate (K13) at (barycentric cs:Z19=1,Z39=1,Z59=1.5);
     \coordinate (K14) at (barycentric cs:Z19=1,Z39=1,Z79=1.5);
     \coordinate (K15) at (barycentric cs:B19=1,Z39=1,Z79=1.5);
     \foreach \oltikzj in {1,2,3,4,5,6,7,8,9,10,11,12,13,14,15}
     \node at (K\oltikzj) {\tiny $K_{\oltikzj}$};
     \coordinate (T) at (3,-1);
     \node at (T)  {$K_1$}
     child {node[circle,draw,fill=a!12.5!g] {$K_2$}}
     child {node {$K_3$}
       child {node[circle,draw,fill=b!18.75!g] {$K_6$}}
       child {node[circle,draw,fill=c!12.5!g] {$K_7$}}
     };
     \coordinate (A1S) at (6,-4);
     \draw (A1S)--++(A1B1) coordinate (B1S)--++(B1C1) coordinate (C1S)--cycle;
     \draw[fill=a!12.5!g,line join=bevel]
     (C1S)--++(C1Z1) coordinate (Z1S)--(A1S)--cycle;
     \draw[fill=b!18.75!g,line join=bevel]
     (Z1S)--++(Z1Z3) coordinate (Z3S)--(C1S)--cycle;
     \draw[fill=c!12.5!g,line join=bevel]
     (Z1S)--(Z3S)--(B1S)--cycle;
     \coordinate (SK2) at (barycentric cs:A1S=1,C1S=1,Z1S=1);
     \coordinate (SK6) at (barycentric cs:C1S=1,Z1S=1.2,Z3S=1);
     \coordinate (SK7) at (barycentric cs:B1S=1,Z1S=1,Z3S=1);
     \foreach \oltikzj in {2,6,7}
     \node at (SK\oltikzj) {\tiny $K_\oltikzj$};
     \node at ($(A1S)+(1.25,.25)$) {mesh $\mesh k$};
   
     \coordinate (A1R) at (9,-4);
     \draw (A1R)--++(A1B1) coordinate (B1R)--++(B1C1) coordinate (C1R)--cycle;
     \draw[line join=bevel]
     (B1R)--++($(Z2)-(B1)$) coordinate (Z2R)--(A1R)--cycle;
     \draw[line join=bevel]
     (Z2R)--++($(Z3)-(Z2)$) coordinate (Z3R)--(C1R)--cycle;  
     \node at ($(A1R)+(1.25,.25)$) {mesh $\mesh l$};
   \end{tikzpicture}
  \end{center}
  \caption{
    Three successive refinements by bisection of a macro element
    $K_1=\simplex\qb{{\vec a}^1,{\vec b}^1,{\vec c}^1}$ where ${\vec c}^1$ is the
    refinement vertex. For each $i$ new node ${\vec z}^i$ on triangle $K_i$
    is created as the midpoint of segment $\simplex\qb{{\vec a}^i,{\vec b}^i}$.\\
    The triangle $K_i$, where $i:=2^l+j$ for some $l\geq0$ and
    $j\integerbetween0{2^l-1}$ (i.e., $l:=\floor{\log_2 i}$ and
    $j:=i-2^l$), is then split into triangles $K_{2i}$ (left child, or
    child $0$) and $K_{2i+1}$ (right child, or child $1$). The
    vertices of the new triangles are named by the following rules:\\
    \(  \begin{matrix}
        {\vec a}^{2i}:={\vec a}^i
        &
        {\vec a}^{2i+1}:={\vec c}^i
        \\
        {\vec b}^{2i}:={\vec c}^i
        &
        {\vec b}^{2i+1}:={\vec b}^i
        \\
        {\vec c}^{2i}:={\vec z}^i
        &
        {\vec c}^{2i+1}:={\vec z}^i
      \end{matrix}
    \)
    where the ${\vec c}$ vertices and $\simplex\qb{{\vec a},{\vec b}}$
    edges are always the one to be refined at the next bisection.\\
    Note that the numbering of simplices is only for notational
    convenience and is generally not used in practice (objects and
    references such as pointers to data structures are used).\\
    The last row shows a bisection tree with the corresponding
    triangulation.  The leaves of the tree correspond to the mesh
    elements $\mesh k=\setofsuch{K_i}{i=2,3,6}$.  The second mesh,
    $\mesh l$, (also obtained by successive bisections of $K_1$, but
    starting from a different labelling) is \emph{not compatible} with
    $\mesh k$.
  }
  \label{fig:triangle-refinement-by-bisection}
\end{figure}

  It is worth noting that not all simplicial partitions of $\Omega$ can be
  represented by such forests, but when two meshes are generated by the
  same macro partition we say that they form a
  \indexemph{compatible} pair of meshes.
}%
In this \namecref{sec:residual-estimators}, we consider
given a compatible pair $\mesh k$ and $\mesh l$ of $\Omega$. It
can be seen that in this that if $K\in\mesh k$ either
\begin{enumerate}[(a)\ ]
\item
    \label{item:K-is-finer-than-L}
    for some element
    $L_K\in\mesh l$ we have $\closure K\subsetneq\closure L_K$\\
    or
  \item
    for some submesh
    $\mesh l_K$ we have $\closure K=\unions L{\mesh l_K}\closure L$.
  \end{enumerate}
  If \ref{item:K-is-finer-than-L} occurs for all $K\in\mesh k$ we say
  that $\mesh k$ is strictly coarser than $\mesh l$ or that $\mesh l$ is
  strictly finer than $\mesh k$.  This induces a partial ordering and a
  Boolean structure on the forest of $\frak T$.
  
  We write also write $\sidesofmesh k$ for the set of sides of $\mesh k$ and
  denote the union of such sides with
  \begin{equation}
    \Sigma_{\mesh k} = \unions S{\sidesofmesh k}\closure S.
  \end{equation}
  
  If $E$ is an element of $\mesh k$ or $\sidesofmesh k$, we denote its
  diameter by $\meshsizeof E$.  The \indexemph{mesh-size} of the mesh
  $\mesh k$ is the piecewise constant function defined by
  \begin{equation}
    \meshsizemesh k(\vec x)=
    \begin{cases}
      \meshsizeof K\tif 
      \vec x\in\interior K\text{ (interior of $K$)}
      \Forsome K\in\mesh k,
      \\
      \meshsizeof S\tif
      \vec x\in S
      \Forsome
      S\in\sidesofmesh k.
    \end{cases}
  \end{equation}

In the rest of this \namecref{sec:residual-estimators} we will consider a pair of compatible
meshes $\mesh k$ and $\mesh l$ upon which we build the conforming finite
element spaces
\begin{equation}
  \label{eq:def:compatible-finite-element-spaces-VW}
  \fespace w:=\poly k(\mesh k)\meet\elldom
  \tand
  \fespace v=\poly k(\mesh l)\meet\elldom,
\end{equation}
where $\elldom:=\sobhz[\Gamma_0]1(\Omega)$ and $\pivot:=\leb2(\Omega)$.
\longonly[\par]{\subsection{Residual estimators}}%
For $\fe w\in\fespace w$, noting that $\ellop\fe w$ belongs to the
dual space $\ellran$ but generally not to the pivot space $\pivot$, In
fact, the distribution $\ellop\fe{w}$ can be decomposed into a
\indexemph{regular part} and a singular \indexemph{jump part}
\begin{equation}
  \begin{split}
    \ellopmesh k\fe{w}
    :=
    -
    \sums K{\mesh k}
    \charfun K
    \divof{c\grad\fe{w}}
    \Aein
    \W%
    \\
    \ellop[\sidesofmesh k]\fe{w}
    :=
    \sums S{\sidesofmesh K}
    \charfun S
    \jump[S]{c\grad\fe{w}}
    \Aeon[\area]
    {\Sigma_{\mesh k}}
    \\
    \text{ where }
    \jump[S]{\vec\psi(\vec x)}
    :=
    \sum_{\substack{K\in\mesh k\\\closure K\supseteq\closure S}}
    \restriction{\vec\psi}K(\vec x)\inner\normalto K(\vec x)
    \tand
    \restriction{\vec\psi}K(\vec x)
    :=
    \lim_{\theta\to0}
    \vec\psi(\vec x-\theta\normalto K(\vec x))
    ,
  \end{split}
\end{equation}
with $\normalto K$ the outer boundary normal to $K$ and
$\vec{w}\in\cont0(\mesh{k})$.

\longonly{%
The following hold
\begin{equation}
  \begin{split}
    \ellopmesh k\fe{w}
    \in\leb2(\W)
    \tand
    \\
    \ltwop{\ellopmesh k\fe{w}}{\vec\chi}
    =
    -
    \sum_{K\in\mesh k}
    \int_K{
      \divof{
        c(\vec x)
        \grad\fe{w}(\vec x)
      }
    }{
      \chi(\vec x)
    }\d\vec x
    \Forall\chi\in\leb2(\W),
  \end{split}
\end{equation}
as well as
\begin{equation}
  \begin{split}
    \ellop_{\sidesofmesh k}\in\leb2(\Sigma_{\mesh k})
    \tand
    \\
    \ltwopon{\ellop_{\sidesofmesh k}\fe{w}}{\chi}{\Sigma_{\mesh k}}
    =
    \sum_{S\in\sidesofmesh k}
    \int_S{\jump{c(\vec x)\grad\fe{w}(\vec x)}}{\chi(\vec x)}\ds(\vec x)
    \Forall
    \chi\in\leb2(\Sigma_{\mesh k}).
  \end{split}
\end{equation}
To summarize we have for each $\fe W\in\fespace W$ and $\phi\in\elldom$
\begin{equation}
  \label{eq:regular-singular-decomposition}
  \duality{\ellop\fe W}{\phi}
  =
  \ltwop{\ellop_{\mesh k}\fe W}{\phi}
  +
  \ltwopon{\ellop_{\sidesofmesh k}\fe W}{\phi}{\Sigma_{\mesh k}}
\end{equation}
where the $\phi$ on $\Sigma_{\mesh k}$ is understood as the trace
of $\phi$.}%

The associated Babuška--Rheinboldt \aposteriori error estimator 
\citep{BabuskaRheinboldt:78:article:Error}
\begin{equation}
  \label{eq:babuska-rheinboldt}
  \cE_{\operatorname{BR}}[\fe{w},\fespace v,\linspace z]
  :=
  \Normonleb[big]{
    \qppow{\meshsizemesh l}\sigma
    \qp{\discellop[\fespace v]\fe w-\ellopmesh k\fe{w}}}2\W
  +
  \Normonleb[big]{
    \powqp{\sigma-1/2}{\meshsizeof{\sidesofmesh l}}
    \ellop_{\sidesofmesh k}\fe W}2
  {\Sigma_{\mesh k}}
\end{equation}
where $\sigma=1$ if $\linspace z=\elldom$ and $\sigma=2$ if $\linspace z=\leb2(\W)$.
\subsection{Discrete elliptic operators and elliptic reconstructors}
Given a conforming finite element space, say $\fespace w\subseteq\elldom$,
we define the corresponding \indexemph{discrete elliptic operator}
\begin{equation}
  \dfunkmapsto{\discellop[\fespace w]}w{\elldom}{\discellop[\fespace w]w}{\fespace w}
\end{equation}
defined (thanks to Riesz representation) by
\begin{equation}
  \ltwop{\discellop[\fespace w]w}{\fe\Phi}=\abil w{\fe\Phi}
  \Foreach\fe\Phi\in\fespace w.
\end{equation}
Alternatively we can think of $\discellop[\fespace w]=\lproj[\fespace w]\ellop$,
where $\funk{\lproj[\fespace w]}{\dualspace v}{\fespace w}$ is the $\leb2$ projector
onto $\fespace w$.

Denote by $\ellrecop[\fespace w]$ the \indexemph[elliptic reconstruction]{elliptic
reconstruction with respect to $\fespace w$}, defined by
\begin{equation}
  \ellrecop[\fespace w]
  =
  \inverse\ellop
  \discellop[\fespace w]
  =
  \inverse\ellop
  \lproj[\fespace w]\ellop
  .
\end{equation}
Note that $\funk{\ellrecop[\fespace w]}{\elldom}{\elldom}$ has finite dimensional range.
We can now state\longonly[only, and omit the proof]{} the three auxiliary's results
needed to use the elliptic residual estimators in the time-dependent problems with
time-varying meshes.
\begin{Lem}[two-space residual \aposteriori error estimate]
  \label{lem:residual-error-aposteriori}
  Suppose $\fespace v\subseteq\fespace w$, and $\linspace z$ one of
  $\pivot$ or $\elldom$, then for all $\fe W\in\fespace W$ we have
  \begin{equation}
    \Normonspace{\ellrecop[\fespace v]\fe{w}-\fe{w}}z
    \leq
    \cE_{\operatorname{BR}}[\fe{w},\fespace v,\linspace z].
  \end{equation}
\end{Lem}
\longonly{
  \begin{Proof}
    Let $\fe W\in\fespace w$.
    The key observation is that
    \begin{equation}
      \ellrecop[{\fespace v}]\fe {w}-\fe {w}
      \orthogonalto[\ellop]
      \vespace[].
    \end{equation}
    This results immediately from the definition
    \begin{equation}
      \abil{\ellrecop[\fespace v]\fe {w} }{\fe\Phi}
      =
      \ltwop{\discellop[\fespace v]\fe {w} }{\fe\Phi}
      =
      \abil{\fe {w} }{\fe\Phi}
      \Foreach
      \fe\Phi\in\fespace v.
    \end{equation}
    Respectively writing $\mesh K$ and $\mesh L$ for the finite element
    mesh underpinning $\fespace w$ and $\fespace v$, and
    $\funk{\clementop[\fespace v]}{\leb2(\W)}{\fespace v}$ for the
    Clément--Scott--Zhang interpolant onto $\fespace v$, it follows that
    for any $\phi\in\elldom$
    \begin{equation}
      \begin{split}
        \abil{
          \qb{\ellrecop[\fespace v]\fe{w}-\fe{w}}
        }{
          \phi
        }
        &
        =
        \abil{
          \qb{\ellrecop[\fespace v]\fe{w}-\fe{w}}
        }{
          \phi-\clementop[\fespace v]\phi
        }
        %
        %
        %
        %
        %
        %
        %
        %
        %
        %
        %
        %
        %
        %
        %
        %
        %
        %
        %
        %
        %
        %
        %
        %
        %
        %
        %
        %
        %
        %
        %
        %
        %
        %
        %
        %
        %
        %
        %
        %
        %
        %
        %
        %
        %
        %
        %
        %
        %
        %
        %
        %
        \\
        &
        =
        \ltwop{
          \qp{\discellop[\fespace v]\fe W-\ellopmesh K\fe W}\meshsizemesh l
        }{
          \qp{\meshsizemesh l}^{-1}\qp{\phi-\clementop[\fespace v]\phi}
        }
        \\
        &
        \phantom=
        +
        \ltwopon{
          \ellop_{\sidesofmesh k}\fe W\qp{\meshsizeof{\sidesofmesh l}}^{\fracl12}
        }{
          \qp{\meshsizemesh l}^{-\fracl12}\qp{\phi-\clementop[\fespace v]\phi}          
        }{\Sigma_{\mesh k}}  
      \end{split}
    \end{equation}
    where we used that
    \begin{equation}
      \ellop\ellrecop[\fespace v]\fe w=\discellop[\fespace v]\fe w
    \end{equation}
    and the decomposition (\ref{eq:regular-singular-decomposition})
    of $\ellop\fe w$ into regular and singular part.

    Taking $\phi=\ellrecop[\fespace v]W-W$ we obtain
    \begin{equation}
      \begin{split}
        \potenorm{\ellrecop[\fespace v]\fe{w}-\fe{w}}^2
        &=
        \abil{
          \qb{\ellrecop[\fespace v]\fe{w}-\fe{w}}
        }{
          \phi
        }
        \\
        &
        =
        \ltwop{
          \qp{\discellop[\fespace v]\fe W-\ellopmesh K\fe W}\meshsizemesh l
        }{
          \qp{\meshsizemesh l}^{-1}\qp{\phi-\clementop[\fespace v]\phi}
        }
        \\
        &
        \phantom=
        +
        \ltwopon{
          \ellop_{\sidesofmesh k}\fe W\qp{\meshsizeof{\sidesofmesh l}}^{\fracl12}
        }{
          \qp{\meshsizemesh l}^{-\fracl12}\qp{\phi-\clementop[\fespace v]\phi}          
        }{\Sigma_{\mesh k}}.  
      \end{split}
    \end{equation}
    By \CBS, the Clément--Scott--Zhang inequalities on compatible meshes
    \citep{LakkisMakridakis:06:article:Elliptic} and the coercivity
    of $\ellop$ we have
    \begin{equation}
      \begin{split}
        \potenorm{\ellrecop[\fespace v]\fe{w}-\fe{w}}^2
        &
        \leq
        \constext{CSZ}
        \Normonleb{\gradof{\ellrecop[\fespace v]\fe{w}-\fe{w}}}2{\W}
        \\
        &
        \qp{
          \Normonleb[big]{
            \qp{\discellop[\fespace v]\fe W-\ellopmesh K\fe W}\meshsizemesh l
          }2{\Omega}
          +
          \Normonleb[big]{
            \ellop_{\sidesofmesh k}\fe W\qp{\meshsizeof{\sidesofmesh l}}^{\fracl12}
          }2{\Sigma_{\mesh k}}}
        \\
        &
        \leq
        \constext{CSZ}\constref[\flat]{eqn:ellop-Lax-Milgram}
        \potenorm{\ellrecop[\fespace v]\fe{w}-\fe{w}}
        \cE_{\operatorname{BR}}[\fe{w},\fespace v,\ellop]
        .
      \end{split}
    \end{equation}
    The proof for the case $\linspace z=\pivot$ follows the same line, albeit by
    testing  with $\inverse\ellop\qb{\ellrecop[\fespace v]\fe w-\fe w}$;
    details are found in \citet[\S2.4]{AinsworthOden:00:book:A-posteriori}.
  \end{Proof}
}
\begin{Lem}[reconstructions on two different spaces]
  \label{lem:reconstrution-and-two-spaces}
  Let $\fespace v$ and $\fespace w$ be two compatible conforming
  finite element spaces, $\linspace z=\pivot$ or
  $\elldom$. Respectively denote by $\ellrecop[\fespace w]$ and
  $\ellrecop[\fespace v]$ the elliptic reconstructors with respect
  to $\fespace w$ and $\fespace v$, then for each $\fe v\in\fespace v$
  and $\fe w\in\fespace w$ we have
  \begin{equation}
    \Normonspace{\ellrecop[\fespace w]\fe w+\ellrecop[\fespace v]\fe v}z
    \leq
    \Normonspace{\fe w+\fe v}z
    +
    \cE[\fe w+\fe v,\fespace w\meet\fespace v,\linspace z].
  \end{equation}
\end{Lem}
\longonly{
\begin{Proof}
  Note that
  \begin{equation}
    \begin{split}
      \Normonspace{\ellrecop[\fespace w]\fe w+\ellrecop[\fespace v]\fe v}z
      \leq
      \Normonspace{\fe w+\fe v}z
      +      
      \Normonspace{
        \ellrecop[\fespace w]\fe w-\fe w
        +
        \ellrecop[\fespace v]\fe v-\fe v
      }z
      .
    \end{split}
  \end{equation}
  But
  \begin{equation}
    \ellrecop[\fespace w]\fe w-\fe w
    \orthogonalto[\ellop]\fespace w
    \tand
    \ellrecop[\fespace v]\fe v-\fe v
    \orthogonalto[\ellop]\fespace v
  \end{equation}
  imply that
  \begin{equation}
    \ellrecop[\fespace w]\fe w-\fe w
    +
    \ellrecop[\fespace v]\fe v-\fe v
    \orthogonalto[\ellop]
    \fespace w\meet\fespace v.
  \end{equation}
  Following the proof of Lemma \ref{lem:residual-error-aposteriori}
  yields the result.
\end{Proof}
}%
\begin{Lem}[reconstruction on the coarser space]
  \label{lem:reconstruction-on-coarser-space}
  Let $\fespace v\subseteq\fespace w$ be two compatible conforming
  finite element spaces, $\linspace z=\leb2(\W)$ or
  $\sobhz1(\W)$. Denote by $\ellrecop[\fespace v]$ the elliptic
  reconstructor with respect to $\fespace v$ and $\cE$ the error
  estimator functional, then for each $\fe {w}\in\fespace w$ we have
  that
  \begin{equation}
    \Normonspace{\ellrecop[\fespace v]\fe {w}}z
    \leq
    \cE[\fe {w},{\vespace[]},\linspace z]
    +
    \Normonspace{\fe {w}}z
    .
  \end{equation}
\end{Lem}
\longonly{\begin{Proof}
  Using the triangle inequality and Lemma \ref{lem:residual-error-aposteriori}
  we write
  \begin{equation}
    \begin{split}
      \Normonspace{\ellrecop[\fespace v]\fe {w}}z
      \leq
      \Normonspace{\ellrecop[\fespace v]\fe {w}-\fe {w}}z
      +
      \Normonspace{\fe {w}}z
      \\
      \leq
      \cE[\fe {w},{\vespace[]},\linspace z]
      +
      \Normonspace{\fe {w}}z
      .
    \end{split}
  \end{equation}
\end{Proof}
}
\bibliographystyle{abbrvnat}

\ifthenelse{\boolean{shownotes}}{
  \color{b!50!g}
  \printindex
  \color i
}{}
\end{document}